\documentclass{amsart}
\usepackage{graphicx} % Required for inserting images
\usepackage[a4paper, total={6.5in, 9in}]{geometry}
\newcommand{\Fa}{\mathrm{Fa}}
\newcommand{\Fi}{\mathrm{Fi}}
\usepackage{amsmath,amssymb,amsthm,mathtools}
\usepackage{booktabs,longtable,array,tabularx}
\usepackage{setspace}
\usepackage{hyperref}
\usepackage{comment}
\usepackage{mathtools}
\usepackage{relsize}
\usepackage{enumitem}
\usepackage{subcaption}
\usepackage{bm}
\usepackage{natbib}
\usepackage{tikz,tikz-cd}
\usepackage{hhline}
\newcommand{\N}{\mathbb{N}}
\newcommand{\Z}{\mathbb{Z}}

\newcommand{\R}{\mathbb{R}}

\newcommand{\lt}{\left}
\newcommand{\rt}{\right}
\newcommand{\half}{\frac{1}{2}}
\newcommand{\trm}{\textrm}
\newcommand{\inv}[1]{\frac{1}{#1}}
\newcommand{\op}{\operatorname}
\newcommand{\poch}[3]{\left({#1};{#2}\right)_{#3}}

\newcommand{\qmultinom}[3]{\genfrac{[}{]}{0pt}{}{{#1}}{{#2}}_{{#3}}}

\numberwithin{equation}{section}

\newtheorem{theorem}{Theorem}[section]
\newtheorem{conjecture}[theorem]{Conjecture}
\newtheorem{proposition}[theorem]{Proposition}
\newtheorem{lemma}[theorem]{Lemma}
\newtheorem{corollary}[theorem]{Corollary}
\newtheorem{question}[theorem]{Question}
\theoremstyle{definition}
\newtheorem{definition}[theorem]{Definition}
\newtheorem{remark}[theorem]{Remark}
\newtheorem*{remark*}{Remark}
\newtheorem{example}[theorem]{Example}

\title{Rational Knots, Quivers, and Tails of Colored HOMFLY-PT Polynomials}
\author{Elijah Berry \and Rongjiang Guo \and Jonathan A. Higgins \and Eduardo Perez \and Henry Zaslow}
\address{Department of Mathematics, University of Illinois Urbana-Champaign, 1409 W Green st, Urbana, IL, 61801, USA}
\date{October 2026}
\thanks{This project was supported by NSF grant DMS-2405302 to Jacob Rasmussen and NSF grant DMS-2246537 to Richard Laugesen.}

\begin{document}

\begin{abstract}
    We use the knots-quivers correspondence to prove that the colored HOMFLY-PT polynomial of a rational knot has a tail and give an explicit formula for calculating the tail in terms of the quiver.
\end{abstract}

\maketitle

\section{Introduction}
\label{IntroSection}

The colored Jones polynomial is a well-known invariant of knots which associates a sequence of Laurent polynomials to each knot. One of the interesting phenomena it exhibits is that for certain classes of knots, the coefficients of highest and lowest $q$-powers stabilize as the color increases. This stabilization was first observed by Dasbach and Lin \cite{dasbach2006headtailcoloredjones}, where they proved that the first three leading coefficients of the colored Jones polynomial for an alternating knot stabilized up to a sign, and they conjectured that this stabilization would continue for all leading coefficients, forming a ``tail." Later, the tail of the colored Jones polynomial of an alternating link $L$ was proven to exist by Garoufalidis and L\^{e} \cite{garoufalidis2012nahmsumsstabilitycolored} using generalized Nahm sums. The tail was also independently proven to exist for a more general class of knots called adequate knots by Armond \cite{Armond_2013}. Finally, Beirne and Osburn \cite{Beirne_2017} ask whether a tail exists in some sense for generalizations of the colored Jones polynomial.

The colored HOMFLY-PT polynomial of a knot is a sequence of Laurent polynomials in $q$ with coefficients in $\mathbb{Z}[a^{\pm 1}]$. After a substitution, it specializes to the colored Jones polynomial. Thus, it is natural to wonder if there is a similar stability phenomenon in the coefficients of the colored HOMFLY-PT polynomial. 

In this paper, we use the links-quivers correspondence \cite{KRSS19} to prove that the normalized colored HOMFLY-PT polynomials of every rational knot converge coefficient-wise to a formal power series $\Psi_K(a,q)$, which we call the tail of the colored HOMFLY-PT polynomial. In addition, we provide an explicit formula for the tail in terms of the quivers for the knots.  

\subsection{Limits of a sequence of polynomials}\label{subsec:lims}

The colored Jones polynomial can be thought of as a sequence of Laurent polynomials with integer coefficients. 
Consider the $j$th colored Jones polynomial of $K$, $J_K^j(q)$, normalized so that $J_K^j(q)=1$ for $K$ the unknot. Note that for $j=1$ this is the ordinary Jones polynomial. By dividing $J_K^j(q)$ by its monomial with lowest power of $q$, we can write the normalized $j$th colored Jones polynomial of $K$ to be
\[
\hat{J}_K^j(q):=\sum_{i=0}^\infty c_{i}^jq^{2i}
\]
for $c_{i}^j \in \mathbb{Z}$. We then say that $\hat{J}_K^j(q)$ converges coefficient-wise to $\Phi_K(q) \in \mathbb{Z}[[q]]$ as a formal power series if for every $i \ge 0$, there exists an $m \ge 0$ such that the coefficient of $q^{2i}$ is independent of $j$ for all $j \ge m$. In this case, we say that the tail of the colored Jones polynomial $\Phi_K(q)$ exists. In fact, for all knots in which the tail exists, it is the case that the coefficient of $q^{2i}$ stabilizes when $m=i$ \cite{Armond_2013}.

Likewise, one can also normalize by the monomial with the greatest power of $q$ and obtain the ``head" of the colored Jones polynomial. Note that the head of $K$ is determined by the tail of $-K$, the mirror knot of $K$. 

Now we define the tail of the colored HOMFLY-PT polynomial analogously to the colored Jones polynomial. Recall that for a knot $K$ the $j$th (antisymmetric) colored HOMFLY-PT polynomial $P^j_K(a,q) $ is a Laurent polynomial in two variables and that one can obtain the colored Jones polynomial of $K$ from the specialization $q \to q^{-1}, a\to q^2$. Of the monomials of $P^j_K(a,q)$ with minimum $q$-power, choose the monomial with the greatest $a$-power (there may be multiple monomials with minimum $q$-power, so it is by convention that we choose the one with greatest $a$-power). Then by dividing $P^j_K(a,q)$ by this monomial, we can define the normalized $j$th colored HOMFLY-PT polynomial of $K$ to be
\[
\hat{P}^j_K(a,q):=\sum_{s\in \Z,r\ge0}c_{s,r}^ja^{2s}q^{2r}
\]
for $c_{s,r}^j \in \mathbb{Z}$. We then say that $\hat{P}_K^j(a,q)$ converges coefficient-wise to $\Psi_K(a,q) \in \mathbb{Z}[a^{\pm1}][[q]]$ as a formal power series if for every $i \ge 0$, there exists an $m \ge 0$ such that the coefficient of $q^{2i}$ is independent of $j$ for all $j \ge m$. In this case, we say that the tail of the colored HOMFLY-PT polynomial $\Psi_K(a,q)$ exists.

Analogously to the Jones case, if we instead normalized $P^j_K(a,q)$ by the monomial with maximum $q$-power and lowest $a$-power, then we would obtain the head of the colored HOMFLY-PT polynomial.

\subsection{Knots-Quivers Correspondence}
    
 A knot $K$ satisfies the knots-quivers correspondence \cite{KRSS17, KRSS19} if there exist linear forms $R, S, A$ and a quadratic form $Q$ defined on a free $\mathbb{Z}$-module $V 
 \cong\mathbb{Z}^u$, such that the generating function for the colored HOMFLY-PT polynomials can be written as
\begin{equation}\label{eq:homgen}
    \sum_{r\ge0}P_K^r(a, q)x^r=\sum_{d\in\N^u}(-1)^{R\cdot d}q^{S\cdot d}a^{A\cdot d}q^{d^TQd}\genfrac{[}{]}{0pt}{}{|d|}{d_1,...,d_u}_{q^2}x^{|d|}.
\end{equation}

In \citep{Sto_i__2021}, Sto\v si\'{c} and Wedrich prove that the knots-quivers correspondence holds for all arborescent knots, extending their previous proof for rational knots \citep{Sto_i__2019}. They detail how to calculate these forms algebraically for rational knots. Additionally, the third author gives a geometric method for calculating these quivers \cite{higgins2}. Note that these methods may give different quivers, but they produce the same generating function and are equivalent up to a permutation of off-diagonal entries in a manner allowable by \cite{JKLN21}. We will thus refer to quivers obtained by the third author's method as \textit{geometric quivers} and those obtained by Sto\v si\'{c} and Wedrich's method as \textit{algebraic quivers}. The quivers we give as examples in this paper will be geometric quivers, but some of the proofs about properties of quivers will be in terms of algebraic quivers. 

\subsection{Main Result}

In order to prove the existence of the tail of the colored HOMFLY-PT polynomial, we consider a distinguished submatrix of the quiver $Q$ for a knot $K$ which controls which vertices contribute to any fixed coefficient of the limiting power series. Terms involving vertices outside of these submatrices will have arbitrarily high degree of $q$ as the color increases. In order to guarantee these submatrices exist and can produce a tail, the quivers must have certain properties. Thus, we give the following definition.

\begin{definition}\label{def:classq}
    Let $K$ be a knot which satisfies the knots-quivers correspondence. Then $K$ is of class $\mathcal{Q}$ if it admits quiver data $Q, S, A$ such that:
    \begin{enumerate}
        \item $Q$ has a unique minimum on its diagonal
        \item the minimum of diag$(Q)$ is the minimum (not necessarily unique) of the whole quiver.
        \item $\text{diag}(Q)=S-2A$
    \end{enumerate} 
\end{definition}

Now, the main theorem can be stated as follows.

\begin{theorem}\label{thm:homfly-tail}
    The tail of the colored HOMFLY-PT polynomial exists for knots $K$ of class $\mathcal{Q}$ and is determined by the restriction of the quiver data to quadratic and linear forms $Q^t, S^t, A^t, R^t$ defined on a certain submodule $V^t \cong \mathbb{Z}^m$ of $V$ with a formula given by
    \begin{equation}\label{eq:hom}
        \Psi_K(a,q)=\sum_{d\in\N^{m}}(-1)^{R^t\cdot d}\frac{a^{A^t \cdot d}q^{d^TQ^td+S^t \cdot d}}{\prod_{i=1}^{m}(q^2;q^2)_{d_i}}.
    \end{equation}
    
\end{theorem}

The precise definitions of these restricted subquivers will be given in Section \ref{sec:homfly-tail}. In Section \ref{sec:Q} we will prove that rational knots are of class $\mathcal{Q}$. Consequently, Theorem \ref{thm:homfly-tail} holds for all rational knots. 

If we know the tail of the colored HOMFLY-PT polynomial, one can obtain the head of the colored Jones polynomial (i.e., the Jones tail of the mirror knot) by substituting $q=q^{-1}$ followed by $a=q^2$. More generally, by instead specializing $a\to q^N$, one obtains a head of the colored $\mathfrak{sl}_N$ polynomials. 

Furthermore, we will see that the tail of the colored HOMFLY-PT polynomial is a stronger invariant than the tail of the colored Jones polynomial. This will be proven in Section \ref{sec:ex} by explicitly computing the HOMFLY-PT tails of the torus knots and twist knots. In Proposition \ref{diff_tails} we will show that these two families have different HOMFLY-PT tails yet specialize to the same Jones heads. 

\subsection{Quiver Computation Program}

\begin{figure}
    \centering
    \includegraphics[width=0.5\linewidth]{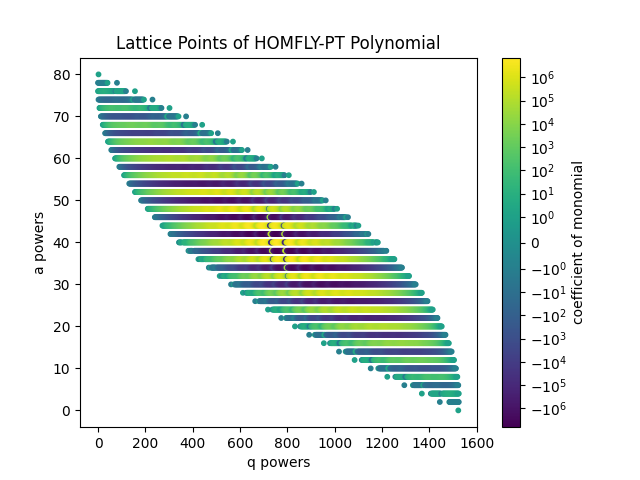}
    \caption{The 20th-Colored HOMFLY-PT Polynomial of $K_{5/2}$}
    \label{fig:K5/2}
\end{figure}

This paper originated from a UIUC summer IML project supported by NSF grant DMS-2405302.
As a part of the project, we wrote a Python program to compute the quivers and colored HOMFLY-PT polynomials of rational knots utilizing the third author's geometric method of computing the quivers \cite{higgins2}. This approach let us compute these knot invariants efficiently, which let us investigate if the coefficients stabilized. For example, our program can produce graphs, as shown in Figure \ref{fig:K5/2}, for visualizing the shape of the colored HOMFLY-PT polynomials. After seeing that tails were forming, we used the program to investigate the structure of the quivers, from which the data suggested that all of the rational knots are of class $\mathcal{Q}$. The program, figures, and associated data are available at \cite{homfly-tails}.

\subsection{Structure of the Paper}

Our paper proceeds as follows. In Section \ref{sec:background}, we give necessary background on rational knots, the colored HOMFLY-PT polynomial, quantum algebra, and the knots-quivers correspondence. In Section \ref{sec:homfly-tail}, we prove Theorem \ref{thm:homfly-tail} for knots of class $\mathcal{Q}$. In Section \ref{sec:properties}, we explore various other properties of the quiver and colored HOMFLY-PT polynomials. In Section \ref{sec:Q}, we prove that rational knots are of class $\mathcal{Q}$ by showing that the necessary properties hold after a sequence of top and right twists. In Section \ref{sec:ex}, we compute the tail of the colored HOMFLY-PT polynomial for various knots, such as the torus knots, twist knots, and all rational knots up to numerator $30$. Finally, in Section \ref{sec:questions}, we conclude with some questions and conjectures we have concerning the quivers and colored HOMFLY-PT polynomials.

\subsection{Acknowledgments}

We would like to thank Jake Rasmussen for leading this project and for his support. Additionally, we wish to thank Christine Ruey Shan Lee for introducing us to the tail of the colored Jones polynomial.

\subsection{AI Disclosure}

The idea for the proof that rational knots are of class $\mathcal{Q}$ given in Section 5 was suggested by ChatGPT. It attempted a proof, which, although incorrect, was useful in giving us a starting point in proving the claim. Additionally, we used ChatGPT and Claude to find various quantum algebra identities that we used in the computations in several of the HOMFLY-PT tails. Lastly, we were assisted by ChatGPT and Claude in writing the code for our program.

\section{Background}\label{sec:background}

\subsection{Rational Knots}
\label{RKsection}

The rational knots are a family of alternating knots which can be defined inductively by finite sequences of simple twist rules. Even though rational knots are invertible, we will use the results of \cite{higgins2,Sto_i__2019}, where rational knots are taken as closures of oriented rational tangles; consequently, we will define them as such here. 

\begin{figure}
\centering
    \begin{tikzpicture}
        \node at (0,0) {\includegraphics[height=3cm]{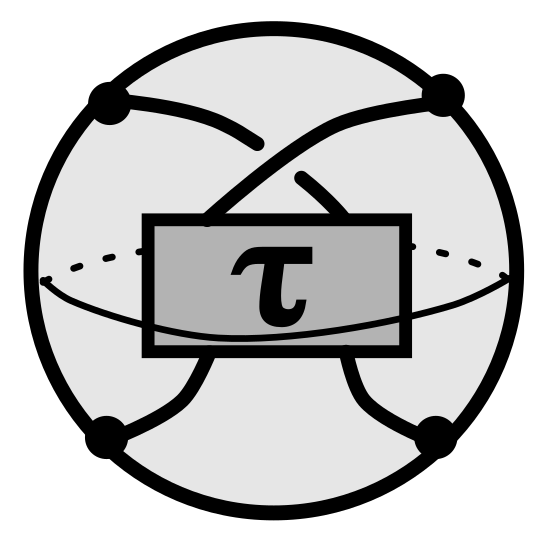}};
        \node at (4.5,0) {\includegraphics[height=3cm]{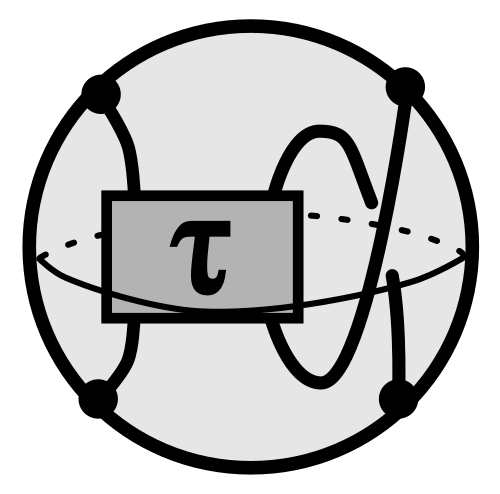}};
        \node at (9,0) {\includegraphics[height=3cm]{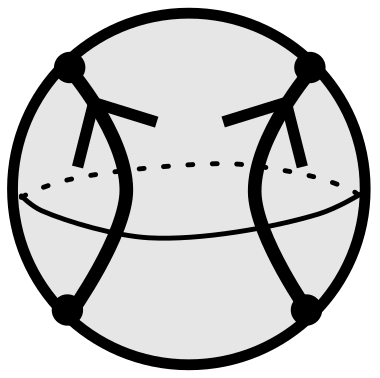}};
        \node at (0,-1.8) {Top Twist: $T$};
        \node at (4.5,-1.8) {Right Twist: $R$};
        \node at (9,-1.8) {Trivial Tangle: $\tau_{0/1}$};
    \end{tikzpicture}
    \caption{The two twists and trivial tangle used to build positive rational tangles}
    \label{twistsfig}
\end{figure}

\begin{definition}
    The \textit{(positive) rational tangles} are the infinite family of oriented 4-ended tangles defined inductively by applying finite sequences of top twists ($T$) and right twists ($R$), as shown in Figure \ref{twistsfig}, to the \textit{trivial tangle}, which we call $\tau_{0/1}$.
\end{definition}

Negative rational tangles can be defined similarly by taking mirrors, but we will focus our discussion on positive rational tangles, following \cite{higgins2,Sto_i__2019}.
%\newpage

As Figure \ref{twistsfig} suggests, it is natural to think of rational tangles as pairs of oriented embedded arcs in the $3$-ball $B^3$ with boundary on $\partial B^3=S^2$. Furthermore, rational tangles correspond to rational numbers by the sequence of twists used to construct them. Given $u/v\in \mathbb{Q}$, we will write $\tau_{u/v}$ to denote the rational tangle corresponding to this rational number. There is a topological way to understand this correspondence via Dehn twists and double branched covers (see \cite{higgins1} for more details), but we will only provide the combinatorial description here, which follows from the rules:
\[
T\tau_{u/v}=\tau_{(u+v)/v}\qquad \qquad \text{and}\qquad\qquad R\tau_{u/v}=\tau_{u/(u+v)}.
\]

Observe first that $R\tau_{0/1}$, the tangle we get by applying a right twist to the trivial tangle, is equal to $\tau_{0/1}$ since the two tangles are related by a Reidemeister I move. Consequently, any rational tangle can be built by first applying a top twist to $\tau_{0/1}$. Then, given a sequence of twists $T^{a_1}R^{a_2}T^{a_3}...R^{a_{n-1}}T^{a_n}\tau_{0/1}$ ending in a top twist gives the rational tangle $\tau_{u/v}$, and a sequence $R^{b_1}T^{b_2}...R^{b_{n-1}}T^{b_n}\tau_{0/1}$ ending in a right twist gives $\tau_{u'/v'}$, where $u/v$ and $u'/v'$ are given by the continued fraction expansions
\[
u/v=a_1+\frac{1}{a_2+\frac{1}{a_3+...+\frac{1}{a_{n-1}+\frac{1}{a_n}}}}
\qquad\qquad\text{and}\qquad\qquad
u'/v'=\frac{1}{b_1+\frac{1}{b_2+...+\frac{1}{b_{n-1}+\frac{1}{b_n}}}},
\]
respectively. The compositions of twists should be read from right to left (like function composition), and the exponents (the $a_i$ and $b_i$) are all positive integers. 

\begin{remark}
    In the first case, we have that $n$ is odd and $u/v\geq1$; in the second, $n$ is even and $u'/v'<1$.
\end{remark}

Positive rational tangles are oriented objects, but there are only three different possible orientations of the arcs at the boundary. We adopt the following convention for naming the orientations from \cite{W16}.

\begin{figure}
    \centering
    \begin{tikzpicture}
        \node at (0,0) {\includegraphics[height=3cm]{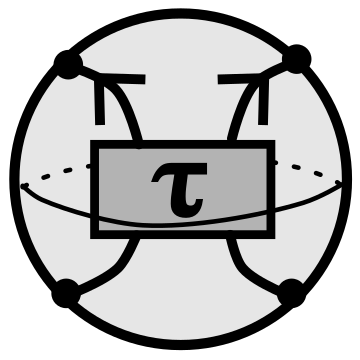}};
        \node at (4.5,0) {\includegraphics[height=3cm]{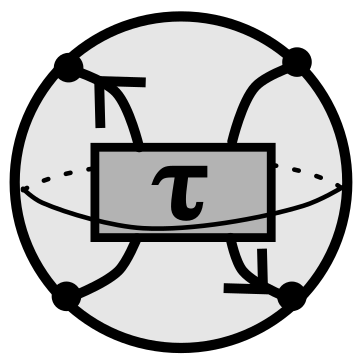}};
        \node at (9,0) {\includegraphics[height=3cm]{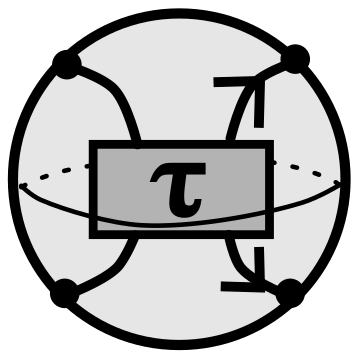}};
        \node at (0,-1.8) {$UP$};
        \node at (4.5,-1.8) {$OP$};
        \node at (9,-1.8) {$RI$};
    \end{tikzpicture}
    \caption{The three possible orientations for positive rational tangles}
    \label{orientationsfig}
\end{figure}
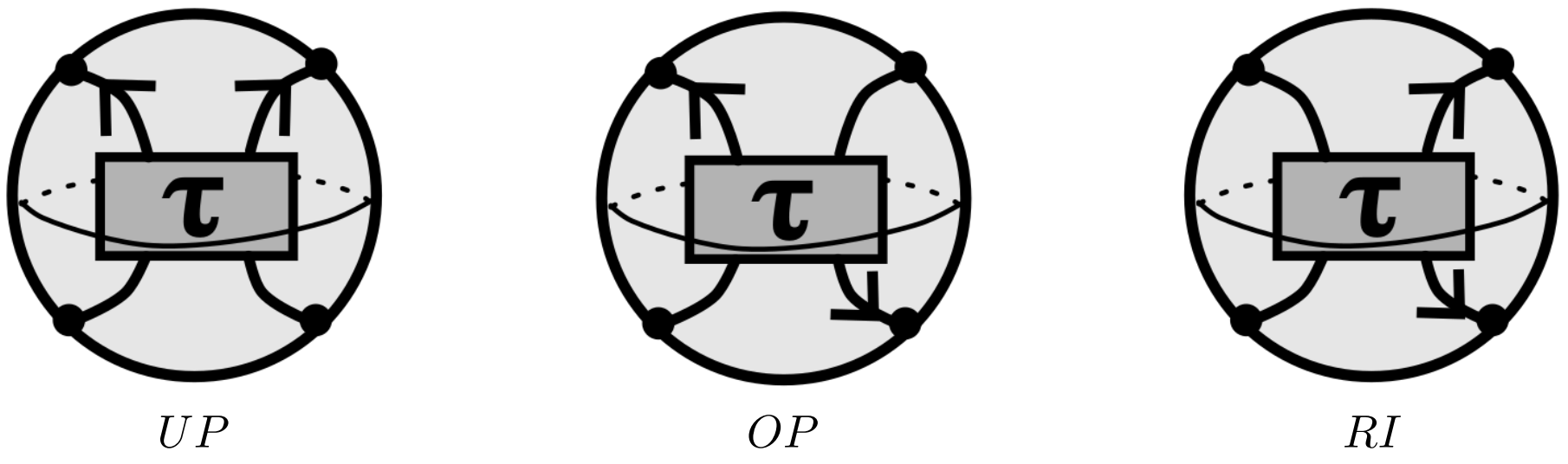

\begin{definition}
    Let $UP,OP,$ and $RI$ be the names of the three possible orientations for positive rational tangles, as shown in Figure \ref{orientationsfig}.
\end{definition}

There are two types of closure operations for (positive) rational tangles: the numerator and denominator closures. However, we will only be concerned with numerator closures of rational tangles $\tau_{u/v}$ with $u/v\geq 1$, as we will justify below.
%, as all rational knots (or their mirrors) can be obtained by taking the numerator closure of one of these tangles by Lemma 2.6 of \cite{Sto_i__2019}. Furthermore, the same lemma says that t is a knot if and only if the numerator $u$ is odd

\begin{definition}
    The \textit{numerator closure} of a rational tangle $\tau_{u/v}$ is the link obtained by identifying the top left (resp. right) endpoints with the bottom left (resp. right) endpoints, if the orientation allows. We will use $Cl(\tau_{u/v})$ to denote the numerator closure of $\tau_{u/v}$. The \textit{denominator closure} can be defined similarly by identifying the top and bottom endpoints.
\end{definition}

A couple of observations must be made. We are only taking the numerator closure of rational tangles with $u/v\geq 1$ because doing so to $\tau_{u/v}$ with $u/v<1$ results in nugatory crossings that can be removed by Reidemeister I moves. Additionally, the numerator closure is not compatible with all three orientations. In particular, we can only take the numerator closure of tangles with $UP$ or $OP$ orientation. 

As was previously stated, rational knots are invertible, which means that we lose nothing by ignoring orientations. However, \textit{rational links}, or the links obtained by taking the closures of positive or negative rational tangles, may have one or two components; in the case of $2$-component rational links, there are two possible orientation choices which result in different colored HOMFLY-PT polynomials. Nevertheless, we restrict our attention to closures of positive rational tangles that yield single-component rational knots.

 Naturally, one should take the numerator closure of $\tau_{u/v}$ when $|u/v|\geq 1$ and the denominator closure when $|u/v|<1$. It is well-known, however, that it is sufficient to only consider numerator closures of rational tangles $\tau_{u/v}$ with $u/v\geq 1$, and that the closure is a knot if and only if $u$ is odd (otherwise, it is a $2$-component link).

\begin{definition}
    Given $\tau_{u/v}$ as an unoriented tangle with $u/v\geq 1$ and $u$ odd, let $K_{u/v}=Cl(\tau_{u/v})$ be the numerator closure of $\tau_{u/v}$. This is the $u/v$\textit{-rational knot}.
\end{definition}

Next, we state some facts about rational knots (for additional background, see \cite{Ku96}).

\begin{remark}
    Two rational knots $K_{u/v}$ and $K_{u'/v'}$ are equal if and only if $u=u'$ and either $v=v'$ or $vv'\equiv 1 \,(\text{mod}\, u)$. Furthermore, $K_{u/v}$ and $K_{u'/v'}$ are mirrors of each other if and only if $u=u'$ and either $vv'\equiv -1 \,(\text{mod}\, u)$ or $v'=u-v$. Consequently, one criterion for $K_{u/v}$ being amphichiral is $v^2\equiv -1 \,(\text{mod}\,u)$.
\end{remark}

The following lemma shows that, when taking orientation into account as in \cite{higgins2,Sto_i__2019}, we obtain all rational knots (up to mirrors).

\begin{lemma}[Lemma 2.6 of \cite{Sto_i__2019}]
\label{RKClLemma}
 Every rational knot or its mirror can be obtained by taking the numerator closure of a positive rational tangle $\tau_{u/v}$ having  $UP$ or $OP$ orientation and $u/v\geq 1$. Furthermore, in this case, $\tau_{u/v}$ has $UP$ orientation when $v$ is odd and $OP$ orientation when $v$ is even. 
\end{lemma}

This lemma dramatically reduces the number of rational tangle closures we need to consider. Later, we will discuss computations related to the closures of positive rational tangles as in Lemma \ref{RKClLemma}; the corresponding data for the mirror knots can easily be obtained from this.

We conclude the section with an example.

\begin{example}
    Using the continued fraction decomposition $\frac{5}{2}=2+\frac{1}{1+\frac{1}{1}}$, we see that $\tau_{5/2}=T^2RT\tau_{0/1}$ and that it has $OP$ orientation. Thus, we have (ignoring the embeddings in $B^3$):
    \[
    \begin{tikzpicture}
        \node at (0,0) {\includegraphics[height=2.5cm]{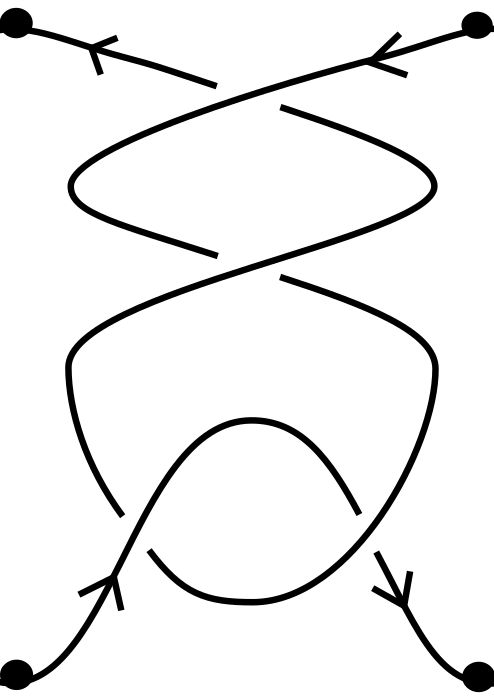}};
        \node at (2,0) {$\xrightarrow{Cl(-)}$};
        \node at (5,0) {\includegraphics[height=2.8cm]{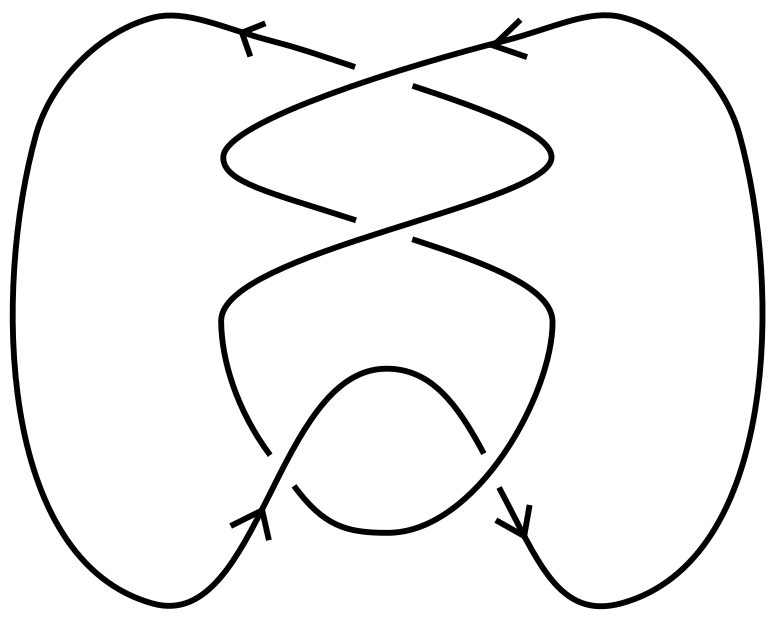}};
        \node at (0,-1.7) {$\tau_{5/2}$};
        \node at (5,-1.7) {$4_1=K_{5/2}$};
    \end{tikzpicture}
    \]
\end{example}

\subsection{Colored HOMFLY-PT Polynomials and Their Specializations}

The Reshetikhin-Turaev invariants, originally defined in \cite{RT90}, give link polynomials that are computed by coloring the components in a diagram for the link with irreducible representations of a simple Lie algebra (or, more accurately, the quantum group corresponding to the Lie algebra). More generally, invariants for tangles are obtained in the same way. In this paper, however, we are only concerned with Type-A invariants of knots, that is, the colored $\mathfrak{sl}_N$ polynomials and colored HOMFLY-PT polynomials, which can be thought of as the $N\to\infty$ limit of the $\mathfrak{sl}_N$ polynomials. We will mainly be interested in the anti-symmetric (single column) colored HOMFLY-PT polynomials and the symmetric (single row) colored $\mathfrak{sl}_N$ polynomials, which are related by simple variable substitutions. In the latter case, the symmetric-colored $\mathfrak{sl}_2$ polynomials are precisely the colored Jones polynomials.

Both the colored HOMFLY-PT and $\mathfrak{sl}_N$ polynomials can be computed via skein theories, as in \cite{MOY98}. Given a tangle $\tau$, one can compute its colored $\mathfrak{sl}_N$/HOMFLY-PT skein module evaluation as a linear combination of MOY-type webs. In the case of rational tangles, Wedrich used quantum skew Howe duality in \cite{W16} to provide explicit formulas for the effects of top and right twists on the skein module evaluations of colored rational tangles. The skein relations then allow one to obtain the desired polynomial invariant for $Cl(\tau_{u/v})$. We will not say anything further about direct calculations of these invariants, but the curious reader is encouraged to consult the sources listed above. For now, our goal is to establish notation, conventions, and list a few important facts that will be used later.

\begin{definition}
  Given a knot $K$, let $P^j_K(a,q)\in\mathbb{Z}[q^{\pm1},a^{\pm1}]$ denote its $\bigwedge^j$-colored ($j$th antisymmetric-colored) reduced HOMFLY-PT polynomial.  
\end{definition}

Observe that the reduced HOMFLY-PT polynomial of a knot is actually a Laurent polynomial in $q$ and $a$. The unreduced polynomial, however, is only a rational expression in the two variables (which is generally the case for (un)reduced HOMFLY-PT polynomials of links with multiple components).

\begin{remark}
    The conventions used in this paper differ from those found in some others by replacing the variables with their squares (i.e., $q\mapsto q^2$ and $a\mapsto a^2$). For example, the uncolored HOMFLY-PT polynomial of the positive trefoil, $3_1=K_{3/1}$, should be 
    \[
    P_{K_{3/1}}^1(a,q)=a^2q^2-a^4+a^2q^{-2}
    \]
    with respect to our conventions. These conventions eliminate the need for fractional exponents and are consistent with grading conventions in categorification.
\end{remark}

Next, we establish notation for the colored $\mathfrak{sl}_N$ polynomials.

\begin{definition}
    Given a knot $K$ and $N\geq 3$, let $F^{j,N}_K(q)\in \mathbb{Z}[q^{\pm1}]$ denote its $Sym^j$-colored ($j$th symmetric-colored) reduced $\mathfrak{sl}_N$ polynomial. In the case where $N=2$, we will write $J^j_K(q)\in\mathbb{Z}[q^{\pm1}]$ for the $j$th colored Jones polynomial.
\end{definition}

%If we were instead to consider the $\bigwedge^j$-colored $\mathfrak{sl}_N$ polynomial $F_K^{\bigwedge^j,N}(q)$, we would have
%\begin{equation}
%F_K^{\bigwedge^j,N}(q)=P^j_K(a,q)|_{a\mapsto q^N}=P^j_K(q^N,q).
%\end{equation}
%However, these polynomials are not ideal for studying as an infinite family since $P_K^{\bigwedge^j,N}(q)=0$ for $j>N$. Consequently, we study the $Sym^j$-colored $\mathfrak{sl}_N$ polynomials, 
The symmetric-colored $\mathfrak{sl}_N$ polynomials are related to the antisymmetric-colored HOMFLY-PT polynomials by
\begin{equation}
F^{j,N}_K(q)=P^j_K(a,q)|_{q\mapsto q^{-1},a \mapsto q^N}=P^j_K(q^N,q^{-1}).
\end{equation}
(See \cite{TVW19} for a proof.)

As discussed in Section \ref{RKsection}, we will need to work with mirrors of knots. Given a knot $K$ and $-K$ its mirror, we have
\begin{equation}\label{mirror-spec}
    P^j_{-K}(a,q)=P^j_K(a,q)|_{q\mapsto q^{-1},a \mapsto a^{-1}}=P^j_K(a^{-1},q^{-1}).
\end{equation}
One consequence of this is 
\begin{equation}\label{sln-spec}
    F^{j,N}_{-K}(q)=P^j_K(a,q)|_{a\mapsto q^{-N}}=P^j_K(q^{-N},q).
\end{equation}

\subsection{Quantum Algebra}
We review some facts about quantum algebra which will be useful in our computations. Recall that the \textit{quantum Pochhammer symbol} is defined as \[
\poch{x}{q}{n} = \prod_{k=0}^{n-1}\lt(1-xq^k\rt),\quad \poch{x}{q}{\infty} = \prod_{k=0}^{\infty}\lt(1-xq^k\rt),\quad (x;q)_0=1,
\]
and the \textit{q-analog} of a number is given by \[
[n]_q=\frac{1-q^n}{1-q}=1+q+q^2+\dots+q^{n-1}.
\]
We say \[
[n]_q!=[n]_q\cdot[n-1]_q\cdots[1]_q, \quad[0]_q!=1,
\]
and for $d_1,...,d_k \ge 0$, $n=d_1+\cdots+d_k$, the \textit{quantum multinomial} is defined as \[
\qmultinom{n}{d_1,d_2,\dots,d_k}{q}=\frac{[n]_q!}{[d_1]_q!\cdots [d_k]_q!}=\frac{(q;q)_n}{\prod_{i=1}^k (q;q)_{d_i}}.
\]
We also use the convention that for $n < 0$,  $\frac{1}{(q;q)_n} = 0$.

Many of our calculations in computing the tail of the colored HOMFLY-PT polynomial and Jones polynomial will rely on various quantum algebra identities, most of which can be found in Appendix II of \cite{GR04}. Some of the most important are the Jacobi Triple Product:
\begin{equation}\label{jtp}
    \poch{q}{q}{\infty}\poch{-\tfrac{1}{z}}{q}{\infty}\poch{-zq}{q}{\infty}=\sum_{n\in\Z} q^{\frac{n(n+1)}{2}}z^n
\end{equation}
and the Rogers-Ramanujan identities (see chapter 7 of \cite{A98}):
\begin{equation}\label{eq:rr1}
    \frac{1}{(q;q^5)_\infty(q^4;q^5)_\infty}=\sum_{n\ge0} \frac{q^{n^2}}{(q;q)_n}
\end{equation}
\begin{equation}\label{eq:rr2}
    \frac{1}{(q^2;q^5)_\infty(q^3;q^5)_\infty}=\sum_{n\ge 0} \frac{q^{n^2+n}}{(q;q)_n}.
\end{equation}
The $q$-binomial theorem, which can be attributed to Heinrich August Rothe \cite{bressoud1981identities}, is given as 

\begin{equation}\label{eq:qbinom}
    \frac{(ax;q)_\infty}{(x;q)_\infty}=\sum_{n=0}^\infty \frac{(a;q)_n}{(q;q)_n}x^n,
\end{equation}
from which one can derive the following $q$-series identities
\begin{equation}\label{eq:qbinom3}
    (x;q)_\infty=\sum_{n=0}^\infty \frac{(-1)^nq^{n(n-1)/2}x^n}{(q;q)_n}
\end{equation}
\begin{equation}\label{eq:t/q}
    \frac{(t;q)_K}{(q;q)_K}=\sum_{n=0}^\infty \frac{(-1)^nt^nq^{n(n-1)/2}}{(q;q)_n(q;q)_{K-n}}.
\end{equation}
For finite Pochhammer symbols, we have the following identity in terms of the $q$-binomial:
\begin{equation}\label{eq:binom}
    (x;q)_n= \sum_{i=0}^n(-1)^iq^{i(i-1)/2}x^i\qmultinom{n}{i}{q}.
\end{equation}
One also gets the following by splitting the products:
\begin{equation}\label{eq:tele}
    \poch{x}{q}{m+n}=\poch{x}{q}{m}\poch{xq^m}{q}{n}
\end{equation}
\begin{equation}\label{eq:qinfty}
    (x;q)_\infty=(x;q)_n(xq^n;q)_\infty.
\end{equation}
Another useful identity for Pochhammer symbols is
\begin{equation}\label{eq:aqpoch}
    (a;q)_n=(-a)^nq^{(n^2-n)/2}(a^{-1}q^{1-n};q)_n,
\end{equation}
which one can easily verify from the definition. Then, substituting $a=xq^{-k}$ gives:
\begin{equation}\label{eq:qpochinverse}
    (xq^{-k};q)_n=(-1)^nx^nq^{(n^2-n)/2-kn}(x^{-1}q^{k+1-n};q)_n.
\end{equation}

Now we give some lemmas that will be used in various proofs later in the paper.

\noindent
\begin{lemma}\label{multinomlim}
    Given a vector ${\bf a}\in \N^r$, let $|{\bf a} |=\sum_{i=1}^ra_i$. Then
    \begin{equation}
    \lim_{N\to \infty}\qmultinom{N}{a_1,\dots,a_r, N-|{\bf a}|}{q}=\frac{1}{\prod_{i=1}^r\poch{q}{q}{a_i}}
    \end{equation}
    as a coefficient-wise limit in $\mathbb{Z}[[q]]$.
\end{lemma}
\begin{proof}
For $N \ge |\bf a|$, we have
\begin{align*}
\qmultinom{N}{a_1,\dots,a_r, N-|{\bf a}|}{q}
&= \frac{(q;q)_N}{(q;q)_{N-|\bf a|} \prod_{i=1}^r (q;q)_{a_i}}\\
&= \frac{\prod_{n=1}^N(1-q^n)}{\prod_{n=1}^{N-\bf |a|}(1-q^n)\prod_{i=1}^r (q;q)_{a_i}}\\
&=\frac{\prod_{n=N-|{\bf a}|+1}^N (1-q^n)}{\prod_{i=1}^r\poch{q}{q}{a_i}}.
\end{align*}
Then, for any $n>0$, the coefficient of $q^n$ in the numerator will be $0$ for sufficiently large $N$. Hence
\[\lim_{N\to\infty}\prod_{n=N-|{\bf a}|+1}^N (1-q^n) = 1
\]
and
\[
\lim_{N\to \infty}\qmultinom{N}{a_1,\dots,a_r, N-|{\bf a}|}{q}=\frac{1}{\prod_{i=1}^r\poch{q}{q}{a_i}}\]

\end{proof}

Next, we consider the effect of inverting $q$ on quantum multinomials.

\begin{lemma}\label{multinomswap} Suppose $N=\sum_{i=1}^kd_i$. Then, we have
    \begin{equation}
        \qmultinom{N}{d_1,\dots , d_k}{q} = q^{\sum\limits_{1\le i<j\le k}d_id_j}\qmultinom{N}{d_1,\dots,d_k}{q^{-1}}.
    \end{equation}
\end{lemma}
\begin{proof}
    Using $[n]_q=q^{n-1}[n]_{q^{-1}}$, it is easily shown that \([n]_q!=q^{n(n-1)/2}[n]_{q^{-1}}!\). This gives the identity \[
    \prod_{i=1}^{k}[d_i]_{q}!=q^{\sum_{i=1}^kd_i(d_i-1)/2}\prod_{i=1}^{k}[d_i]_{q^{-1}}!,
    \]
    so we can rewrite the quantum multinomial as \[
    \qmultinom{N}{d_1,d_2,\dots,d_k}{q}=\frac{[N]_q!}{[d_1]_q!\cdots [d_k]_q!}=q^{N(N-1)/2-\sum_{i=1}^kd_i(d_i-1)/2}\frac{[N]_{q^{-1}}!}{[d_1]_{q^{-1}}!\cdots [d_k]_{q^{-1}}!}.
    \]
    The exponent in the final expression can be simplified by
    \[N(N-1)/2-\sum_{i=1}^kd_i(d_i-1)/2=\half\lt(\sum_{i,j=1}^kd_id_j - \sum_{i=1}^kd_i-\sum_{i=1}^k(d_i^2-d_i)\rt)=\sum\limits_{1\le i<j\le k}d_id_j,\]
    which completes the proof.
\end{proof}

Now, we give two additional quantum algebra identities, which both involve sums over multiple indices.

\begin{lemma}[Andrews \cite{Andrews1974}, Theorem 1]
    \begin{equation} \label{ag}
\sum_{n_1, n_2, \dotsc, n_{k-1} \geq 0} \frac{q^{N_1^2 + N_2^2 + \cdots + N_{k-1}^2 + N_i + N_{i+1} + \cdots + N_{k-1}}}{(q;q)_{n_1} (q;q)_{n_2} \cdots (q;q)_{n_{k-1}}} = \frac{(q^i; q^{2k+1})_{\infty} (q^{2k+1-i}; q^{2k+1})_{\infty} (q^{2k+1}; q^{2k+1})_{\infty}}{(q;q)_{\infty}}
\end{equation}
where $k \ge 2$ and $N_{j}=n_j + n_{j+1} + \cdots + n_{k-1}$.
\end{lemma}
This identity is written differently in \cite{Andrews1974}, but the form shown above can be found in \cite{Keilthy_2016}. Additionally, the following gives a formula for the false theta functions $h_b$, which we shall define in Section \ref{sec:ex}.
\begin{lemma}
\[
\sum_{n_1,\dots, n_{k-1}\ge 0}(-1)^{\sum_{j=1}^{k-1}N_j}\frac{q^{\half\sum_{j=1}^{k-1}N_j(N_j+1)}}{(q;q)_{n_1}\cdots (q;q)_{n_{k-1}}}=h_{k+1}(q)
\]
where $k \ge 2$ and $N_{j}=n_j + n_{j+1} + \cdots + n_{k-1}$.
\end{lemma}
\begin{proof}
    The proof for the cases $k$ even and $k$ odd may be found on pages 20 and 22 of \cite{Keilthy_2016} respectively, after the substitution $n_i\to i_{k-i}$.
\end{proof}

\subsection{Knots-Quivers Correspondence}
\label{KQC}

The knots-quivers correspondence of Kucharski, Reineke, Sto\v si\'c, and Su\l kowski \cite{KRSS17,KRSS19} predicts that the generating function for the (antisymmetric or symmetric) colored HOMFLY-PT polynomials for a rational knot should be related to a generating function associated with a symmetric quiver, after a change of variables. 

\begin{definition}
    Given a knot $K$, let $P_K(a,q)$ denote the generating function of its $\bigwedge^j$-colored HOMFLY-PT polynomials, i.e.
    \begin{equation}
        P_K(a,q)=\sum_{j\geq0}P^j_K(a,q)x^j.
    \end{equation}
\end{definition}

\begin{comment}
For a symmetric quiver $Q$ with $|Q|=m$ vertices, the generating function for its cohomological Hall algebra is of the form
\begin{equation}
P_Q(x_1,...,x_m) = \sum_{\textbf{d}=(d_1,..,d_m)\in\mathbb{N}^m} (-q)^{-\langle \textbf{d},\textbf{d}\rangle_Q}\prod_{i=1}^m \frac{x_1^{d_1}...x_m^{d_m}}{(1-q^{-2})...(1-q^{-2d_i})},
\end{equation}
where we also use $Q$ to denote its adjacency matrix, $\langle \textbf{d},\textbf{e}\rangle_Q = \sum_{i,j}(\delta_{ij}- Q_{ij})d_ie_j$ is its Euler form, and $x_1,...,x_m$ are formal variables. Given a knot $K$, the knots-quivers correspondence claims that one can find a symmetric quiver $Q_K$ such that 
\begin{equation}
    \widehat{P}(K):=\sum_{j\geq 0}(q^2;q^2)_j^{-1}P^j_K(a,q)x^j=P_{Q_K}(x_1,...,x_{|Q_K|})|_{x_i\mapsto (-1)^{Q_{ii-S_i}}q^{S_i-1}a^{A_i}x},
\end{equation}
where $S=(S_1,...,S_{|Q_K|}), A=(A_1,...,A_{|Q_K|})\in \mathbb{Z}^{|Q_K|}$. 
\end{comment}
Practically, the most helpful way for us to think of the knots-quivers correspondence is in terms of $P_K(a,q)$ being in ``quiver form,'' as in \cite{higgins2,Sto_i__2019,Sto_i__2021}.

\begin{definition}
    We say $P_K(a,q)$ is in \textit{quiver form} if it is expressed as
    \begin{equation}
    \label{quivform0}
        P_K(a,q)=\sum_{\textbf{d}=(d_1,...,d_m)\in\mathbb{N}^m}(-1)^{R\cdot \textbf{d}}a^{A\cdot \textbf{d}}q^{S\cdot \textbf{d}+\textbf{d}^T\cdot Q\cdot \textbf{d}}\qmultinom{|\textbf{d}|}{d_1,d_2,\dots,d_m}{q^2}x^{|\textbf{d}|}
    \end{equation}
    where $R,A,S\in \mathbb{Z}^m$ and $Q$ is a symmetric $m\times m$ matrix. We call these forms the \textit{quiver data} of $K$. 
\end{definition}

For rational knots, we have $R=S$, so Equation (\ref{quivform0}) can be rewritten as
\begin{equation}
    \label{quivform}
        P_K(a,q)=\sum_{\textbf{d}=(d_1,...,d_m)\in\mathbb{N}^m}(-q)^{S\cdot \textbf{d}}a^{A\cdot \textbf{d}}q^{\textbf{d}^T\cdot Q\cdot \textbf{d}}\qmultinom{|\textbf{d}|}{d_1,d_2,\dots,d_m}{q^2}x^{|\textbf{d}|}.
    \end{equation}
Consequently, this is the form of $P_K(a,q)$ that we will call its quiver form. In \cite{Sto_i__2019} and \cite{Sto_i__2021}, Sto\v si\'c and Wedrich proved the knots-quivers correspondence for rational knots and arborescent knots, respectively, by showing that $P_K(a,q)$ can be put in quiver form for each case.

In proving the existence of the tail of the colored HOMFLY-PT polynomial, it will be convenient to modify the quiver. We give the following lemma to justify being able to shift the quiver data by a constant. 

\begin{lemma}\label{lem:shift}

Let $P^j_K(a,q)$ be the $j$-th colored HOMFLY-PT polynomial of a knot with quiver data $R,S,A,Q$. Then for $n_1,n_2,n_3,n_4 \in \mathbb{Z}$, the quiver data $R+n_1,S+n_2,A+n_3,Q+n_4$ produces the same polynomial as $P^j_K(a,q)$ up to an overall sign and a shift of $a$ and $q$ powers.
    
\end{lemma}

\begin{proof}
    Suppose $K'$ has quiver data $R+n_1,S+n_2,A+n_3,Q+n_4$. Then
    \begin{align*}
    P^j_{K'}(a,q) &= \sum_{|\textbf{d}|=j}(-1)^{(R+n_1)\cdot \textbf{d}}a^{(A+n_3)\cdot \textbf{d}}q^{(S+n_2) \cdot \textbf{d} + \textbf{d}^T\cdot (Q+n_4)\cdot \textbf{d}}\qmultinom{|\textbf{d}|}{d_1,d_2,\dots,d_m}{q^2} \\
    &= (-1)^{jn_1}a^{jn_3}q^{j^2n_4+jn_2}\sum_{|\textbf{d}|=j}(-1)^{R\cdot \textbf{d}}a^{A\cdot \textbf{d}}q^{S \cdot \textbf{d} + \textbf{d}^T\cdot Q\cdot \textbf{d}}\qmultinom{|\textbf{d}|}{d_1,d_2,\dots,d_m}{q^2}\\
    &= ((-1)^{jn_1}a^{jn_3}q^{j^2n_4+jn_2})P^j_K(a,q).
    \end{align*}
\end{proof}

It should be noted that, when the linear and quadratic forms are appropriately shifted to compute framed versions of the colored HOMFLY-PT polynomials, a framing shift of $+1$ is given by taking $(n_1,n_2,n_3,n_4)=(-1,-1,-1,1)$. This can be seen by combining Lemma \ref{lem:shift} with the fact that this framing shift corresponds to multiplying $P_K^j(a,q)$ by $(-q)^{-j}a^{-j}q^{j^2}$, as noted in Section 4.3 of \cite{Sto_i__2019}.

\begin{remark}
 One helpful way of thinking about the knots-quivers correspondence is that it is equivalent to associating with a knot $K$ two linear forms $S, A$ and one quadratic form $Q$ on a free $\mathbb{Z}$-module with preferred basis. In particular, the summation indices represent non-negative coordinate vectors with respect to this basis.
\end{remark}

The proof of the knots-quivers correspondence for rational knots in \cite{Sto_i__2019} was inductive in nature, as it involved assembling $S, A,$ and $Q$ sequentially according to the twists used to build the positive rational tangle $\tau_{u/v}$ and then taking the numerator closure. In \cite{higgins2}, the third author gave a direct geometric method for computing the linear and quadratic forms for positive rational knots ($K_{u/v}=Cl(\tau_{u/v})$ with $u/v\geq 1$ and tangle orientation $UP$ or $OP$). As \cite{JKLN21} discusses systematically, the quivers (or quadratic forms) associated with knots by the knots-quivers correspondence are highly non-unique, but \cite{higgins2} provides a way to compute a ``geometrically preferred'' one. We will now provide a brief description of this geometric method, as we will use it later, but the interested reader should consult \cite{higgins2} for further details.

\subsection{Geometric Computation of the Quiver}
\label{GeoQuiverSection}

In \cite{higgins2}, the third author provided a method for computing ``geometrically preferred'' linear and quadratic forms satisfying the Knots-Quivers Correspondence for rational knots via winding numbers of loops in the punctured plane and its second configuration space. In this paper, we will occasionally utilize these geometric methods, so we will briefly remind the reader of them and comment on notational differences, but the reader is encouraged to refer to \cite{higgins2} for a detailed discussion and clear definitions.

\begin{figure}
    \centering
    \begin{tikzpicture}
        \node at (0,0) {\includegraphics[height=4cm]{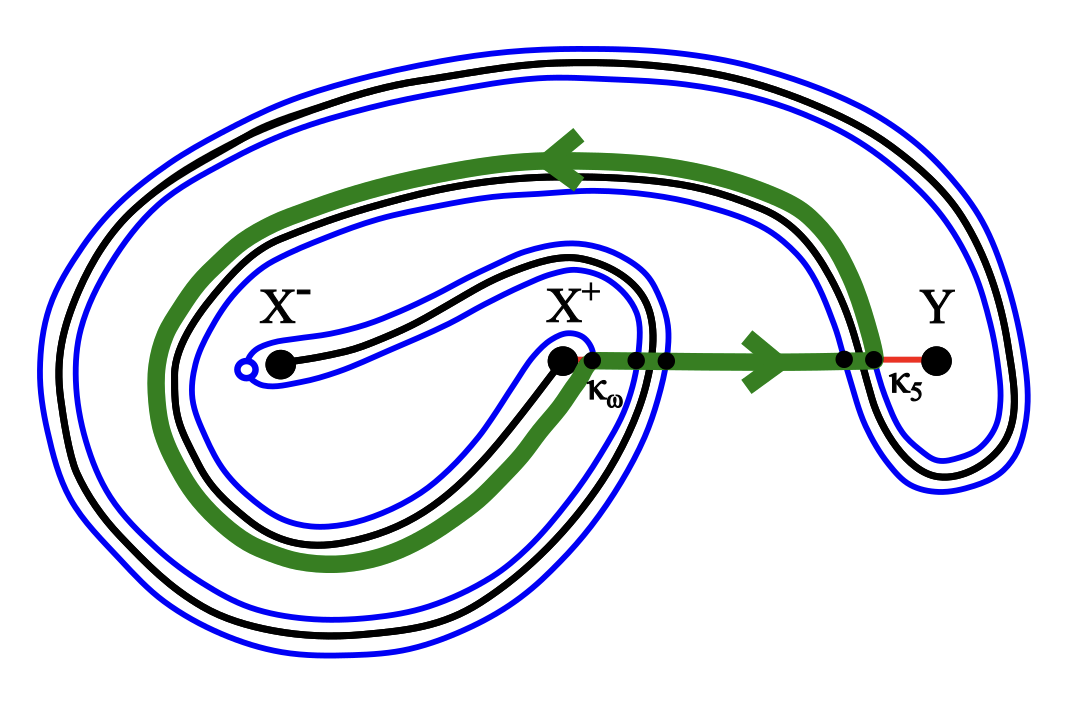}};
        \node at (7,0) {\includegraphics[height=4cm]{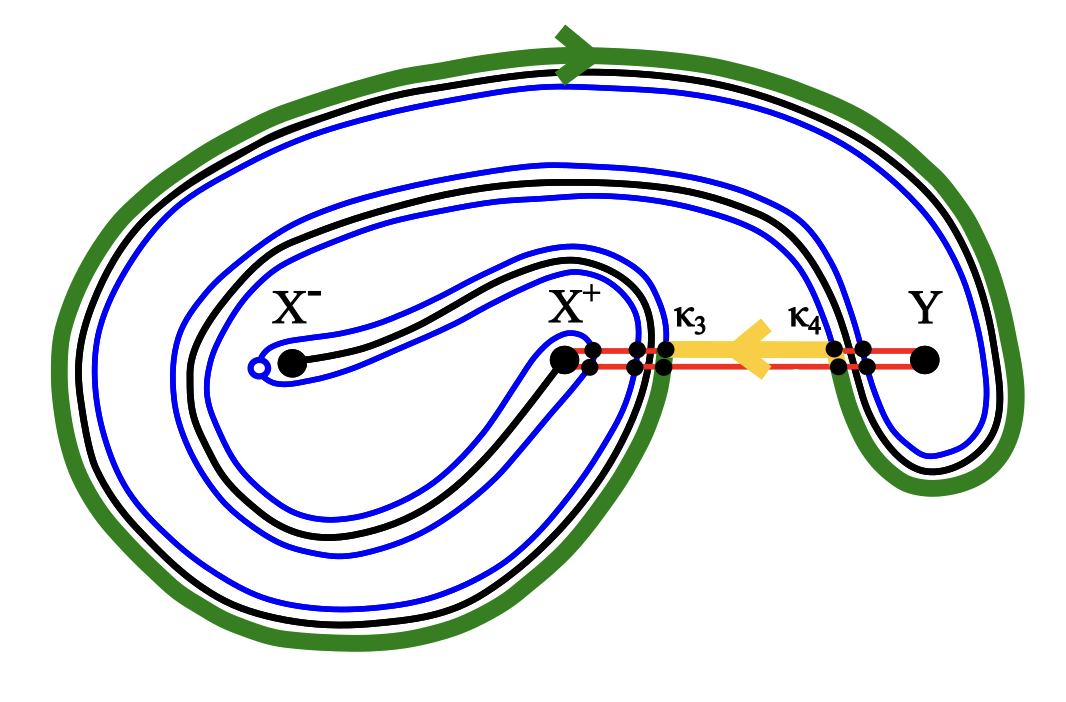}};
        \node at (0,-2) {(a)};
        \node at (7,-2) {(b)};
    \end{tikzpicture}
    \caption{Loops for the knot $K_{5/2}$. $\overline{\alpha_{5/2}}$ and $\mathbb{L}_{5/2}^2$ are drawn in blue; $\beta$ and $\mathbb{L}_\beta^2$ are drawn in red. (a) gives the loop $\gamma_5\subset M$; (b) gives the loop $\gamma_{3,4}^2\subset\text{Conf}_2(M)$.}
    \label{gammaloopsfig}
\end{figure}

Given a rational knot $K_{u/v}$, taken as the numerator closure of a rational tangle $\tau_{u/v}$ with $UP$ or $OP$ orientation, one can associate two types of Lagrangians, denoted $\overline{\alpha_{u/v}}$ and $\beta$, in the $3$-punctured plane. The intersection $\overline{\alpha_{u/v}}\cap \beta$ consists of $u$ points, which provide a preferred basis for the free $\mathbb{Z}$-module on which the linear and quadratic forms are defined. The three punctures are labeled $X^+,X^-,$ and $Y$, and they are permuted by the sequence of top and right twists used to construct the Lagrangians (in fact, there is an inductive procedure for building $\overline{\alpha_{u/v}}$ via this sequence of twists, which is discussed in \cite{higgins2}). The $u$ intersection points in $\overline{\alpha_{u/v}}\cap \beta$ are called $\kappa_1,...,\kappa_u$, which we order from left to right for convenience. Then, given $\kappa_i,\kappa_j\in \overline{\alpha_{u/v}}\cap \beta$, one can define the loop $\gamma_{i,j}$ that follows a specified path from $\kappa_i$ to $\kappa_j$ along $\overline{\alpha_{u/v}}$ (there are two options up to homotopy, but one is ruled out by placing a dot on $\overline{\alpha_{u/v}}$ near $X^-$), and then returning to $\kappa_i$ along $\beta$. There is an intersection point closest to $X^+$, denoted $\kappa_\omega$ in \cite{higgins2}, and many of the necessary loops are of the form $\gamma_i:=\gamma_{i,\omega}$. Such loops are sufficient for computing the linear forms and the diagonal of $Q$. See Figure \ref{gammaloopsfig}(a) for an example.

Letting $M=\mathbb{C}\setminus\{X^+,X^-,Y\}$ denote the $3$-punctured plane, we must now consider analogues of the two Lagrangians in $\text{Conf}_2(M)=\text{Sym}^2(M)\setminus \Delta$, where $\Delta$ is the diagonal subspace of $\text{Sym}^2(M)$. For $\overline{\alpha_{u/v}}$, we just take $\mathbb{L}_{u/v}^2:=\text{Conf}_2(\overline{\alpha_{u/v}})$, and for $\beta$, we define $\mathbb{L}_\beta^2$ to be given by two parallel copies of $\beta$ embedded in $\text{Conf}_2(M)$. Both $\mathbb{L}_{u/v}^2$ and $\mathbb{L}_\beta^2$ are Lagrangian submanifolds of $\text{Conf}_2(M)$, and their intersection consists of $u^2$ points, but we will only need some of these for defining loops. Associated with two points $\kappa_i$ and $\kappa_j$ is a loop $\gamma_{i,j}^2\subset \text{Conf}_2(M)$ defined by taking the appropriate basepoint in $\mathbb{L}_{u/v}^2\cap \mathbb{L}_\beta^2$ at $\kappa_i$ and $\kappa_j$, following a path $\eta_{i,j}$ in $\mathbb{L}_{u/v}^2$ to the point $(\kappa_j,\kappa_j)\in \mathbb{L}_{u/v}^2\cap \mathbb{L}_\beta^2$, and then returning to the basepoint along a path $\overline{\eta}_{j,i}$ in $\mathbb{L}_\beta^2$. An example of such a loop is provided in Figure \ref{gammaloopsfig}(b).

The key result of \cite{higgins2} is that the linear and quadratic forms may be computed in terms of winding numbers of these two types of loops about the punctures in $M$ and the diagonal $\Delta$ in $\text{Conf}_2(M)$. Given $Z\in\{X^+,X^-,Y\}$, let $w_Z(\gamma_{i,j})$ be the winding number of $\gamma_{i,j}$ about $Z$, and let $w_\Delta(\gamma_{i,j}^2)$ be the winding number of $\gamma_{i,j}^2$ about $\Delta$ (this is also discussed more carefully in \cite{higgins2}). Then, we have the following theorem.

\begin{theorem}[Theorem 4.6,\cite{higgins2}]\label{geothm}
     For a rational knot $K_{u/v}$, $P(K_{u/v})$ can be written as
    \begin{equation}
        P(K_{u/v})=\sum_{\textbf{d}=(d_1,...,d_u)\in \mathbb{N}^u}(-q)^{S\cdot \textbf{d}}a^{A\cdot \textbf{d}}q^{\textbf{d} \cdot Q\cdot \textbf{d}^t} \qmultinom{|\textbf{d}|}{d_1,d_2,\dots,d_u}{q^2}x^{|\textbf{d}|}
    \end{equation}
    where $S$, $A$, and $Q$ are defined on the free $\mathbb{Z}$-module with basis $\{\kappa_1,...,\kappa_u\}$, and they are computed by the following formulas:
    \[
    \begin{cases}
        S_i=\displaystyle\sum_{Z\in\{X^+,X^-,Y\}}w_Z(\gamma_i)+\mu_1(K_{u/v})\\
        A_i=2w_{X^+}(\gamma_i)+\mu_2(K_{u/v})\\
        \begin{cases}
            Q_{ii}=w_{X^-}(\gamma_i)+w_Y(\gamma_i)-3w_{X^+}(\gamma_i)+\mu_3(K_{u/v})\\
            Q_{ij}=Q_{ii}+w_\Delta(\gamma_{j,i}^2)-2w_{X^+}(\gamma_{j,i}),
        \end{cases}
    \end{cases}
    \]
    where the $\mu_i(K_{u/v})$ are integers that depend on the types of twists used to build the knot (formulas are given in Section 4.5 of \cite{higgins2}).
\end{theorem}

An important observation to make is that the notation in Theorem \ref{geothm} differs slightly from what is written in \cite{higgins2}. In particular, we use the $w_Z$ and $w_\Delta$ notation to emphasize that we are computing winding numbers. To make the comparison concrete, we have $w_Z(\gamma_{i,j})=\Psi_Z([\gamma_{i,j}])$ and $w_\Delta(\gamma_{i,j}^2)=\Phi([\gamma_{i,j}^2])$, with the $\Psi_Z$ and $\Phi$ notation coming from \cite{higgins2}.

\begin{definition}
    Given a rational knot $K_{u/v}$, we call the quiver $Q$ computed by Theorem \ref{geothm} the \textit{geometric quiver} for $K_{u/v}$.
\end{definition}

Next, we make two remarks before proceeding.

\begin{remark}
\label{fixedquiverrmk}
    By labeling the intersection points $\kappa_1,...,\kappa_u$ from left to right, we can think of $Q$ as a fixed matrix that defines a quadratic form on the free $\mathbb{Z}$-module with ordered basis. This ordering of the basis does not affect the computations of colored HOMFLY-PT polynomials and their tails, but we will require $Q$ to be defined by this ordering of the $\kappa_i$'s in Section \ref{subsec:skeleton} when discussing the ``skeletons" of quivers.
\end{remark}

\begin{remark}
    The code accompanying this paper (as discussed in Section \ref{IntroSection}) was written using Theorem \ref{geothm}. This is why we have chosen to use this method when discussing explicit quiver computations.
\end{remark}

\section{Existence of the Tail of the Colored HOMFLY-PT Polynomial}\label{sec:homfly-tail}

Recall Definition \ref{def:classq}, which says that a knot $K$ satisfying the Knots-Quivers Correspondence is of class $\mathcal{Q}$ if it admits quiver data $Q, S, A$ such that:
    \begin{enumerate}
        \item $Q$ has a unique minimum on its diagonal
        \item the minimum of diag$(Q)$ is the minimum (not necessarily unique) of the whole quiver.
        \item $\text{diag}(Q)=S-2A$
    \end{enumerate} 

It should be noted that, with the third condition, we are assuming that we are allowed to shift all linear and quadratic forms by integers, which we know gives the same colored HOMFLY-PT polynomials up to signs and shifts in $a$- and $q$-powers by Lemma \ref{lem:shift}. For example, diag$(Q)$ is related to $S-2A$ by the knot signature for the quivers computed for rational knots in \cite{higgins2} and \cite{Sto_i__2019}; we simply ignore this shift since we only consider colored HOMFLY-PT polynomials up to multiplication by monomials in $a$ and $q$. This third condition is added to rule out quiver data that yield vanishing tail formulas, which we consider in Example \ref{unknotexample} in the case of the unknot. 
    
\begin{definition}
\label{def:subquiverdata}
Let $K\in \mathcal{Q}$ be a knot with quiver data $Q,R,A,S$ of size $u$, with $Q$ satisfying the minimum diagonal conditions. We permute the quiver data such that $Q_{u,u}$ is the minimum of the diagonal, and normalize the quiver data such that $Q_{u,u}=R_u=A_u=S_u=0$, which is allowable by Lemma \ref{lem:shift}. Now let \[
I = \{i\in \{1,\dots,u-1\}:Q_{u,i}=0\}
\]
and define \[
Q^t = Q_{I,I},\quad R^t = R_I,\quad A^t=A_I,\quad S^t = S_I
\]
to be the restrictions of the quiver data to the indices in $I$. 
We call $Q^t, R^t, A^t, S^t$ the \textit{subquiver data} of the knot.
\end{definition}

Now we restate the definition of the tail of the HOMFLY-PT polynomial given in Section \ref{subsec:lims}.

\begin{definition}\label{def:homftail}
    Let $K$ be a knot and $\hat{P}_K^j(a,q) = \sum_{s\in \Z,r\ge0}c_{s,r}^ja^{2s}q^{2r}$ its normalized $j$th colored HOMFLY-PT polynomial (obtained by normalizing ${P}_K^j(a,q)$ by its monomial with lowest $q$-power and greatest $a$-power).
    
    We then say $\hat{P}_K^j(a,q)$ converges coefficient-wise to $\Psi_K(a,q) \in \mathbb{Z}[a^{\pm1}][[q]]$ as a formal power series if for every $i \ge 0$, there exists an $m \ge 0$ such that the coefficient of $q^{2i}$ is independent of $j$ for all $j \ge m$. In this case, we say that the \textit{tail of the colored HOMFLY-PT polynomial} $\Psi_K(a,q)$ exists.
\end{definition}

In subsequent formulas and examples, we will give the tail up to a monomial normalization. 
Now we can express the tail of the HOMFLY-PT polynomial purely in terms of the subquiver data. We restate Theorem \ref{thm:homfly-tail} as the following.

\begin{theorem}\label{thm:homflytail}
    Let $K \in \mathcal{Q}$ be a knot with quiver data satisfying Definition \ref{def:classq} and associated subquiver data $Q^t, R^t, A^t, S^t$ of size $k$. Then $\Psi_K$ exists and can be expressed as \begin{equation} \label{homftail}
\Psi_K(a,q)=\sum_{v\in\N^k}(-1)^{R^t\cdot v}\frac{a^{A^t\cdot v}q^{v^T\cdot Q^t\cdot v+S^t\cdot v}}{\prod_{i=1}^k \poch{q^2}{q^2}{v_i}}.
\end{equation}
\end{theorem}

Before we prove this theorem, we will give an example to show how one may explicitly compute the HOMFLY-PT tail given the quiver data.

\begin{example}\label{ex:5/2}
We compute the HOMFLY-PT tail of $K_{5/2}$. For the quiver data, we compute using \cite{homfly-tails} and find that up to a constant shift and permutation, \[
R=S= \begin{pmatrix}
    -2\\-1\\0\\-3\\-4
\end{pmatrix},\
A = \begin{pmatrix}
    -2\\-2\\0\\-2\\-4
\end{pmatrix},\
Q=\begin{pmatrix}
    2&2&1&1&2\\2&3&1&1&3\\1&1&0&0&2\\1&1&0&1&2\\2&3&2&2&4
\end{pmatrix}.
\]
Looking at the minimum diagonal entry of $Q$, we see it is at index 3, so we look for all other $0$ values in row/column 3 of $Q$. The only other $0$ entry is at index 4, so we can calculate the subquiver data:
$$S^t = \begin{pmatrix}-3\end{pmatrix},\ A^t = \begin{pmatrix}-2\end{pmatrix},\ Q^t = \begin{pmatrix}1\end{pmatrix}.$$
Then using \eqref{homftail}, we can compute the tail:
\begin{align*}
\Psi_{K_{5/2}}(a,q)&=\sum_{n\ge 0}(-1)^na^{-2n}q^{n^2-3n}\frac{1}{\poch{q^2}{q^2}{n}}\\
&=\poch{a^{-2}q^{-2}}{q^{2}}{\infty}&\text{by \eqref{eq:qbinom3}}.
\end{align*}
To recover the tail of the colored Jones polynomial for $K_{-5/2}=K_{5/2}$, we can specialize $a\to q^{-2}$, giving us $\Phi_{K_{-5/2}}(q)=(q^2;q^2)_\infty$, as originally computed in \cite{garoufalidis2012nahmsumsstabilitycolored}.
\end{example}

Now, we prove Theorem \ref{thm:homflytail}.

\begin{proof}[Proof of Theorem \ref{thm:homflytail}]
Up to constant shifts and conjugation by a permutation matrix, we may say \begin{equation}\label{quiverdatashifted}
Q_{u,u}=0, \qquad Q_{i,i}>0,\,i\ne u,\qquad Q_{i,j}\ge 0, \qquad R_u=S_u=A_u=0.
\end{equation}
Let $v\in \N^{u - 1}$, and write $|v|=\sum_{i=1}^{u-1}v_i$. For $|v|\leq j$, define \(d^v_j=\lt(v_1,v_2,\dots,v_{u-1},j-|v| \rt)\). Now, for each $j$, we can rewrite \eqref{eq:homgen} as \begin{equation}
\label{homreindex}
P_K^j(a,q)=\sum_{\substack{v\in\N^{u-1}\\|v|\le j}}(-1)^{R\cdot d^v_j}a^{A\cdot d^v_j}q^{(d^v_j)^T\cdot Q\cdot d^v_j+{S\cdot d^v_j}}\qmultinom{j}{d^v_j}{q^2}.
\end{equation}

By the conditions in (\ref{quiverdatashifted}), we get the inequality\[
(d^v_j)^T\cdot Q\cdot d^v_j+S\cdot d^v_j \ge \sum_{i=1}^{u-1} \left(Q_{i,i}v_i^2+S_iv_i\right),
\]
from which it follows that, for each $n\in \Z$, there exists $N(n)\in \N$ such that if $|v|>N(n)$, then \((d^v_j)^T\cdot Q\cdot d^v_j+S\cdot d^v_j>n\). Consequently, to find the coefficient of the $q^n$ term in $P^j_K(a, q)$, it suffices to consider only those $d^v_j$ for which $|v|\le N(n)$. Moreover, for sufficiently large $j$, if there exists $1\le k\le u-1$ such that $Q_{k,u}\ne 0$ and $v_k\ne 0$, then $(d^v_j)^T\cdot Q\cdot d^v_j+S\cdot d^v_j>n$.

It follows that, for sufficiently large $j$, the coefficient of the terms with $q$-power at most $n$ is determined by the vectors \[\Gamma_n=\lt\{\lt(v_1, \dots,v_{u-1} \rt)\in \N^{u-1}:
|v| \le N(n),Q_{k, u} \ne 0\implies v_k = 0\rt\}.
\]
Furthermore, if $v$ satisfies these conditions, we can define \[ Q(v):=(d^v_j)^T\cdot Q\cdot d_j^v, \quad A(v):=A(d^v_j), \quad S(v):=S(d^v_j), \quad R(v):=R(d^v_j)\] for any choice of $j\geq |v|$. Additionally, we have that when $j$ and $j'$ are sufficiently large, the coefficients in the binomials \[
\qmultinom{j}{d^v_j}{q^2} \trm{ and } \qmultinom{j'}{d^v_k}{q^2}
\] agree for terms with exponent less than $n$ by Lemma \ref{multinomlim}.
\begin{comment}
then, for all $j, k\ge \sum v_i$, \begin{align*}
    (d^v_j)^T\cdot Q\cdot d_j^v =(d^v_k)^T\cdot Q\cdot d_k^v =: Q(v),
    \ A(d^v_j)=A(d^v_k)=: A(v),
    \ S(d^v_j)=S(d^v_k)=: S(v),
    \ R(d^v_j)=R(d^v_k)=: R(v).
\end{align*}
By Lemma \ref{multinomlim}, for sufficiently large $j, k$, the coefficients of \[
\qmultinom{j}{d^v_j}{q^2} \trm{ and } \qmultinom{k}{d^v_k}{q^2}
\] below $n$ also agree. 
\end{comment}
Therefore, each coefficient of $q^n$ in $P^j_K$ stabilizes when $j$ is sufficiently large. Moreover, since the coefficient of $q^n$ is determined by $\Gamma_n$, which is finite, the coefficient is in $\mathbb{Z}[a^{\pm 1}]$. So $\Psi_K$ exists, and 
\[
\Psi_K(a,q) = \lim _{n\to \infty} \sum_{v\in \Gamma_n} (-1)^{R(v)}a^{A(v)}q^{Q(v)+S(v)}\inv{\prod_{i=1}^{u-1} \poch{q^2}{q^2}{v_i}}.
\]
As $\Gamma_{n}\subset \Gamma_{n+1}$ for all $n$, we can write
\[
\Psi_K(a,q)=\sum_{v\in\Gamma}(-1)^{R(v)}a^{A(v)}q^{Q(v)+S(v)}\inv{\prod_{i=1}^{u-1} \poch{q^2}{q^2}{v_i}},
\]
where 
\[
\Gamma\coloneqq \bigcup_{n\in \N}\Gamma_n=\lt\{\lt(v_1, \dots,v_{u-1} \rt)\in \N^{u-1}:
Q_{u, k} \ne 0\implies v_k = 0\rt\}.
\]

The set $\Gamma$ is simply the restriction of $\N^{u-1}$ to the indices $I = \{i\in \{1,\dots,u-1\}:Q_{u,i}=0\}$, and $Q(v)=v^T\cdot Q_{I,I}\cdot v$. Thus, we can write \[
Q^t = Q_{I,I},\quad R^t = R_I,\quad A^t=A_I,\quad S^t = S_I
\]
and \begin{equation*}
\Psi_K(a,q)=\sum_{v\in\N^k}(-1)^{R^t\cdot v}\frac{a^{A^t\cdot v}q^{v^T\cdot Q^t\cdot v+S^t\cdot v}}{\prod_{i=1}^k \poch{q^2}{q^2}{v_i}},
\end{equation*}
where $k=|I|$.

Now we will show that $\Psi_K(a,q)$ is nonzero. Consider the coefficient of $a^0q^0$. We can see that $v=(0,...,0)\in \mathbb{N}^k$ contributes a coefficient of $1$ to this term, so we want to check that no other vectors cancel out this term. Suppose $v \in \mathbb{N}^k$ is a nonzero vector in the summand that contributes to the term $a^0q^0$. Then $A^t \cdot v = 0$. Since $\text{diag}(Q^t)=S^t-2A^t$,
\begin{align*}
    v^T \cdot Q^t \cdot v + S^t \cdot v &\ge \sum_{i=1}^kQ^t_{i,i} v_i^2 + \text{diag}(Q^t)\cdot v +2A^t\cdot v\\
    &=\sum_{i=1}^kQ^t_{i,i} (v_i^2+v_i).
\end{align*}
This is strictly greater than $0$ since $v\not=0$ and each $Q^t_{i,i} > 0$. Then, since $\frac{1}{\prod_{i=1}^k \poch{q^2}{q^2}{v_i}}$ does not contribute any negative $q$-powers, such a $v$ cannot exist. Therefore, the coefficient of $a^0q^0$ is $1$ and $\Psi_K(a,q)$ is nonzero.

Lastly, to align with Definition \ref{def:homftail}, we show that the $q$-powers of $\Psi_K(a,q)$ are bounded below. We have shown that a term with $q^0$ exists, so we want to show that there are only finitely many negative $q$-powers. Like before, $\frac{1}{\prod_{i=1}^k \poch{q^2}{q^2}{v_i}} \in \mathbb{Z}[[q]]$ does not contribute to any negative $q$-powers for any $v \in \mathbb{N}^n$. 

Let $m =\min\{Q^t_{i,i}\}_{i=1}^k $. For $v \in \mathbb{N}^k$, 
\begin{equation}\label{eq:bounded}
v^T\cdot Q^t \cdot v + S^t \cdot v \ge \sum_{i=1}^k Q_{i,i}^tv_i^2 + S^t \cdot v \ge m||v||_2^2-||S^t||_2||v||_2.
\end{equation}
Then, since  $m > 0$, we have that $m||v||_2^2-||S^t||_2||v||_2 \to \infty$ as $||v||_2\to\infty$. Therefore, there are only finitely many $v \in \mathbb{N}^k$ that contribute to negative $q$-powers. It follows that the $q$-powers of $\Psi_K(a,q)$ are bounded below. We can then renormalize $\Psi_K(a,q)$ by the monomial with the lowest $q$-power and highest $a$-power, ensuring $\Psi_K(a,q) \in \mathbb{Z}[a^{\pm1}][[q]]$.
\end{proof}

For the remainder of the paper, it will be more convenient to consider the tail in the form given in Equation \eqref{homftail}, so we will not be renormalizing it, given that they are equivalent up to a shift of an $a$ and $q$ power.

One may notice that the third condition of the definition of class $\mathcal{Q}$ is only used to show that the tail is nonvanishing. We give the following example to show that there are quivers which satisfy the minimum diagonal conditions (hence they have well-defined subquivers), yet they have a vanishing tail. 

\begin{example}
\label{unknotexample}
    Let $K$ be the unknot. $K$ admits quiver data 
$$R=(0), \quad S=(0), \quad A =(0), \quad Q=(0).$$

We can verify this, as for all $j$: $$P^j_K(a,q)=(-1)^{0\cdot j}q^{0\cdot j}a^{0\cdot j}q^{j\cdot0\cdot j}\genfrac{[}{]}{0pt}{}{j}{j}_{q^2}=1.$$ 

Then the subquiver data will be empty, giving us $\Psi_K(a,q)=1$. 

Now suppose a knot $K'$ has the following quiver data:
$$R'=S' = \begin{pmatrix} 0 \\ -1 \\ 0 \end{pmatrix},
\quad A' = \begin{pmatrix} 0 \\ 0 \\ 0 \end{pmatrix}, 
\quad Q' = \begin{pmatrix} 0 & 0 & 0 \\ 0 & -1 & -2 \\ 0 & -2 & -2 \end{pmatrix}.$$

Calculating the $j$th colored HOMFLY-PT polynomial from this data gives:
\begin{align*}
    P^j_{K'}(a,q)&=\sum_{0 \le n+m\le j}(-1)^nq^{-n^2-2m^2-4nm-n}\qmultinom{j}{j-n-m,n,m}{q^2}\\
    &=\sum_{m=0}^j\sum_{n=0}^m (-1)^n q^{n^2-n-2m^2}\qmultinom{j}{m}{q^2} \qmultinom{m}{n}{q^2} & \text{shifting $m \to m - n$}\\
    &=\sum_{m=0}^j q^{-2m^2} (1;q^2)_m \qmultinom{j}{m}{q^2}& \text{by \eqref{eq:binom}}\\
    &=1 & \text{only the $m=0$ term is nonzero}.
\end{align*}

Hence this is also quiver data for the unknot. It satisfies the minimum diagonal conditions, but it does not satisfy that $\text{diag}(Q')=S'-2A'$. After a shift of $+2$ to $Q'$, we can calculate its subquiver data to be

$$R'^t=S'^t=(-1), \quad A'^t=(0), \quad Q'^t=(1).$$

Then if we attempt to use this subquiver data to calculate the tail, we get
\begin{align*}
    \Psi_{K'}(a,q)&= \sum_{n\ge0} (-1)^nq^{n^2-n}\frac{1}{(q^2;q^2)_n}\\
    &= (1;q^2)_\infty & \text{by \eqref{eq:qbinom3}}\\
    &= 0.
\end{align*}

Therefore, the unknot admits two quivers where one can be used to calculate the tail, and the other cannot.
\end{example}

\begin{remark}
    For rational knots, the $R$ vector is equal to the $S$ vector, so we will not consider it when providing examples in Section \ref{sec:ex}. Note, however, that a general knot in class $\mathcal{Q}$ may have an $R$ vector different from the $S$ vector.
\end{remark}

% We conclude this section by showing that the Jones tail exists for all knots in class $\mathcal{Q}$ and conjecture that the $\mathfrak{sl}_N$ tails exist as well for all $N > 2$.

We conclude this section by discussing specializations of the tail of the colored HOMFLY-PT polynomial. First we show that the HOMFLY-PT tail agrees with the Jones tail.

\begin{corollary}\label{cor:jones}
    If $K\in \mathcal{Q}$, then the tail of the colored Jones polynomial of $-K$ exists and is equal to \(\Psi_K(q^{-2},q)\).
\end{corollary}

\begin{proof}
    From Equation \eqref{homftail} we write 
        $$\Psi_K(q^{-2},q) =\sum_{v\in\N^k}(-1)^{R^t\cdot v}\frac{q^{v^T\cdot Q^t\cdot v+(S^t-2A^t)\cdot v}}{\prod_{i=1}^k \poch{q^2}{q^2}{v_i}}.$$
    Let $m=\min\{Q^t_{i,i}\}_{i=1}^k$. Following the same procedure as in the proof of Theorem \ref{thm:homflytail}, we get that for any fixed $v \in \mathbb N^k$,
    \begin{equation}
        v^T\cdot Q^t \cdot v + (S^t-2A^t) \cdot v \ge m||v||_2^2-||S^t-2A^t||_2||v||_2.
    \end{equation}
    Hence only finite many summation vectors can contribute to any fixed $q$-power, and the formula gives a well-defined coefficient-wise limit. 

    Now since $\text{diag}(Q)=S-2A$, we have that for any nonzero $v \in \mathbb N^k$,
    \begin{equation}
        v^T\cdot Q^t \cdot v + (S^t-2A^t) \cdot v \ge \sum_{i=1}^k Q^t_{i,i}(v_i^2+v_i) > 0.
    \end{equation}
    Hence $v=0$ is the only summation vector that contributes to $q^0$. That is, the $q^0$ coefficient is $1$ and $\Psi_K(q^{-2},q)$ is nonzero.

    Lastly, specializing $a\to q^{-2}$ for $P^j_K(a,q)$ computed from the shifted quiver data as in Equation \eqref{quiverdatashifted} gives us exactly $\hat J^j_{-K}(q)$, the normalized colored Jones polynomial for $-K$. 
    
    Therefore, $$\Psi_K(q^{-2},q) = \lim_{j\to\infty}\hat J^j_{-K}(q)=\Phi_{-K}(q).$$

\end{proof}

Unlike in the proof of Theorem \ref{thm:homflytail}, we do not need to renormalize $\Phi_{-K}(q)$, as the lowest $q$-power of $\Psi_K(q^{-2},q)$ (unnormalized) is exactly $q^0$. For this reason it is natural to let $\Psi_K(a,q)$ be the unnormalized HOMFLY-PT tail obtained from Equation \eqref{homftail}.

Observe that if $K\in \mathcal{Q}$ has a quiver $Q$ satisfying the class $\mathcal{Q}$ conditions (along with the rest of the quiver data), where $Q$ also has a unique maximum entry on its diagonal, which is also the unique maximum for the whole quiver, then $-K$ is also in class $\mathcal{Q}$. This will follow immediately from Proposition \ref{mirrorquiv}, and we will show in Section \ref{sec:Q} that this condition holds for all rational knots. In this case, Corollary \ref{cor:jones} gives
\begin{equation}
\Phi_K(q)=\Psi_{-K}(q^{-2},q).
\end{equation}

Given that the HOMFLY-PT tail has the nice property that specializing $a\to q^{-2}$ recovers the Jones tail, it is natural to consider the $a \to q^{-N}$ specializations for $N > 2$ as well. Rather than considering the coefficient-wise limit, we give a definition for the tail of the colored $\mathfrak{sl}_N$ polynomials in terms of this specialization.

\begin{definition}
    Let $K$ be a knot such that $-K$ is in class $\mathcal{Q}$. We define the tail of the colored $\mathfrak{sl}_N$ polynomial to be 
    \begin{equation}
        \Phi_K^{N}(q) := \Psi_{-K}(q^{-N},q)
    \end{equation}
    for $N \ge 2$.
\end{definition}

\begin{remark}
\label{rmk:slNtailnonvanishing}
    Repeating the same argument as in the proof of Corollary \ref{cor:jones}, one can show that $\Phi_K^N(q)$ is well-defined for knots in class $\mathcal{Q}$. However, we are not guaranteed that the tail is nonvanishing. 

    Suppose a knot $K$ has subquiver data $$R^t= S^t=(5), \quad A^t=(2), \quad Q^t=(1).$$
    This data satisfies the conditions of class $\mathcal{Q}$. One can calculate the HOMFLY-PT tail to be
    \begin{align*}
        \Psi_{K}(a,q)&=\sum_{n\ge 0}(-1)^na^{2n}q^{n^2+5n}\frac{1}{\poch{q^2}{q^2}{n}}\\
        &=\poch{a^2q^{6}}{q^{2}}{\infty}&\text{by \eqref{eq:qbinom3}}.
    \end{align*}

    Then the $\mathfrak{sl}_3$ tail of $-K$ is $$\Phi_{-K}^3(q)=\Psi_{K}(q^{-3},q)=(1;q^2)_\infty=0.$$
    This subquiver data likely does not appear for any knots, but this example shows that if we want the $\mathfrak{sl}_N$ tails to be nonvanishing and to agree with the coefficient-wise limit of the colored $\mathfrak{sl}_N$ polynomials (if it exists), we need to alter the conditions in class $\mathcal{Q}$.
    
    One sufficient condition that would imply nonvanishing would be $A^t \le 0$. Indeed, computational data suggests that this holds for all rational knots.
\end{remark}

%\begin{conjecture}
%    If $K\in \mathcal{Q}$, then for all $N \ge 2$, the tail of the colored $\mathfrak{sl}_N$ polynomial of $-K$ exists and is equal to \(\Psi_K(q^{-N},q)\).
%\end{conjecture}

%\begin{remark}
%    Note that repeating the same proof of Corollary \ref{cor:jones} for this conjecture will still result in coefficient-wise convergence of \(\Psi_K(q^{-N},q)\), but we are not guaranteed that the tail is nonvanishing. One sufficient condition that would imply nonvanishing would be $A^t \le 0$. Indeed, computational data suggests that this holds for all rational knots. We have verified this for all rational knots up to numerator $100$ \cite{homfly-tails}.
%\end{remark}

\section{Properties of Quivers and Tails}\label{sec:properties}

If a knot is of class $\mathcal Q$, then the tail of the HOMFLY-PT polynomial can be computed directly from the quiver. The following propositions describe certain properties of the map from quivers to tails of polynomials which will be useful in later computations.

\begin{proposition}\label{prop:block}
    Let $K$ be a knot of class $\mathcal{Q}$, and suppose $Q_K^t$ is block-diagonal; that is, \[
    Q_K^t=\begin{pmatrix}
        Q_1&0\\0&Q_2
    \end{pmatrix},\quad A_K^t = \begin{pmatrix}
        A_1\\A_2
    \end{pmatrix},\quad S^t_K=\begin{pmatrix}
        S_1\\S_2
    \end{pmatrix},\quad R_K^t = \begin{pmatrix}
        R_1\\R_2
    \end{pmatrix}
    \] 
    where $Q_1, A_1, S_1, R_1$ are size $m$ and $Q_2,A_2,S_2,R_2$ are size $n$.
    Then, \[
    \Psi_K(a,q) = \sum_{v\in\N^m}(-1)^{R_1\cdot v}\frac{a^{A_1\cdot v}q^{v^T\cdot Q_1\cdot v+S_1\cdot v}}{\prod_{i=1}^m \poch{q^2}{q^2}{v_i}}\times
    \sum_{v\in\N^n}(-1)^{R_2\cdot v}\frac{a^{A_2\cdot v}q^{v^T\cdot Q_2\cdot v+S_2\cdot v}}{\prod_{i=1}^n \poch{q^2}{q^2}{v_i}}.
    \]
\end{proposition}
\begin{proof}
    The proof is immediate upon noticing $v^T\cdot Q^t_K\cdot v = (v_1^T\cdot Q_1\cdot v_1) + (v_2^T\cdot Q_2\cdot v_2)$, where $v= \begin{pmatrix} v_1 \\ v_2\end{pmatrix}$ is analogous to the equations above. Then, we have 
    \begin{align*}
        \Psi_K(a,q)&=\sum_{v\in\N^{m+n}}(-1)^{R^t\cdot v}\frac{a^{A^t\cdot v}
        q^{v^T\cdot Q^t\cdot v+S^t\cdot v}}
        {\prod_{i=1}^{m+n} \poch{q^2}{q^2}{v_i}}\\
        &= \sum_{v\in\N^m}\sum_{w\in \N^n}(-1)^{R_1\cdot v+R_2\cdot w}\frac{a^{A_1\cdot v+A_2\cdot w}
        q^{(v^T\cdot Q_1\cdot v+S_1\cdot v)+(w^T\cdot Q_2\cdot w+S_2\cdot w)}}
        {\prod_{i=1}^m \poch{q^2}{q^2}{v_i}\times\prod_{i=1}^n \poch{q^2}{q^2}{v_i}},
    \end{align*}
    so the sum splits.
\end{proof}

Next, we consider how the linear and quadratic forms for a knot are related to the ones for its mirror.

\begin{proposition}\label{mirrorquiv}
    Let $K$ be a knot satisfying the knots-quivers correspondence, with $m\times m$ quiver $Q_K$ and vectors $R_K,A_K,S_K$. Then the mirror knot $-K$ satisfies the knots-quivers correspondence, with quiver $Q_{-K}=I-1-Q_K$ and vectors $R_{-K}=-R_K,\ A_{-K}=-A_K, \ S_{-K}= -S_K$.
    \begin{proof}
        Let $\mathbf{d}\in \N^m$. Then $\mathbf{d}\cdot (I-1-Q)\cdot \mathbf{d}=\mathbf{d}\cdot (-Q)\cdot \mathbf{d} - 2\sum_{1\le i<j\le m}d_id_j$.
        Directly, using Equation \eqref{mirror-spec}, \begin{align*}
        P_{-K}(a,q)=P_K(a^{-1}, q^{-1})=&\sum_{\mathbf{d}=(d_1,\dots,d_m)\in \N^m}(-1)^{R\cdot \mathbf{d}}a^{-A\cdot \mathbf{d}}q^{-S\cdot \mathbf{d}+\mathbf{d}^T\cdot (-Q)\cdot \mathbf{d}}\qmultinom{|\mathbf{d}|}{d_1,d_2,\dots,d_m}{q^{-2}}x^{|\mathbf{d}|}.
        \\
        \intertext{Now, using Lemma \ref{multinomswap} along with the identity above,}
        P_{-K}(a,q)=&\sum_{\mathbf{d}=(d_1,\dots,d_m)\in \N^m}(-1)^{-R\cdot \mathbf{d}}a^{-A\cdot \mathbf{d}}q^{-S\cdot \mathbf{d}+\mathbf{d}^T\cdot (I-1-Q)\cdot \mathbf{d}}\qmultinom{|\mathbf{d}|}{d_1,d_2,\dots,d_m}{q^2}x^{|\mathbf{d}|}.
        \end{align*}
    \end{proof}
\end{proposition}

\subsection{Skeleton of the Quiver}\label{subsec:skeleton} Let $K_{u/v}$ be a rational knot, with (geometric) quiver $Q_{K_{u/v}}$. As noted in Remark \ref{fixedquiverrmk}, we will assume that $Q_{K_{u/v}}$ is computed by Theorem \ref{geothm} with the $\kappa_i$'s ordered from left to right.

\begin{definition}
    We define the \textit{skeleton} of the $u\times u$ quiver $Q_{K_{u/v}}$ to be the $(u-1)\times (u-1)$ matrix $S_{{u/v}}$, where \[
{(S_{u/v})}_{i, j}={(Q_{K_{u/v}})}_{i+1, j+1}-{(Q_{K_{u/v}})}_{i+1, j}-{(Q_{K_{u/v}})}_{i, j+1}+{(Q_{K_{u/v}})}_{i, j}.
\]
\end{definition}

\begin{example}
    Using the quiver of $K_{5/2}$ given in \ref{ex:5/2}, we calculate:
    \[
    S_{5/2}=
    \begin{pmatrix}
          1 & -1 & 0 & 1\\
         -1 & 1 & 0 & 0\\
          0 & 0 & 1 & -1\\
          1 & 0 & -1 & 1\\
    \end{pmatrix}.
    \]
\end{example}

Now we prove an interesting property of the skeleton.

\begin{theorem}\label{skeletontheorem}
    For all $1\le i,j\le u-1$, $(S_{u/v})_{i,i}=1$, and $-1\le {(S_{u/v})}_{i,j}\le 1$. 
\end{theorem}

Let $\rho_{k,l}$ be the path from $\kappa_k$ to $\kappa_{l}$ along $\overline{\alpha_{u/v}}$.
The proof will give an explicit geometric method to calculate ${(S_{u/v})}_{i,j}$ by looking at intersections of $\kappa_i$ and $\kappa_{i+1}$ with $\rho_{j,j+1}$.

\begin{proof}

One can check, using Theorem \ref{geothm}, that
\begin{align*}
    Q_{i,j}-Q_{i,j+1}-Q_{i+1,j}+Q_{i+1,j+1} =&w_\Delta(\gamma_{j, i}) - w_\Delta(\gamma_{j+1, i})-w_\Delta(\gamma_{j, i+1})+w_\Delta(\gamma_{j+1, i+1}).
\end{align*}

Consider the term $w_\Delta(\gamma_{j, i})$. The loop $\gamma_{j, i}$ is the concatenation of two paths in $\op{Conf}_2(M)$: $\eta_{j,i}$, which traverses $\mathbb{L}^2_{u/v}$, and $\overline{\eta}_{i,j}$, which traverses $\mathbb{L}^2_\beta$, as described in Section \ref{GeoQuiverSection}. The winding number can be thought of as the writhe of an associated 2-strand braid, which is the concatenation of two braids associated with $\eta_{j,i}$ and $\overline{\eta}_{i,j}$. Let $$\op{wr}_{j,i}\coloneqq w_\Delta({\eta}_{j,i}) - w_\Delta({\eta}_{j+1,i})-w_\Delta({\eta}_{j,i+1})+w_\Delta({\eta}_{j+1,i+1})$$
where $w_\Delta(\eta_{j, i})$ refers to the writhe of the portion of the braid associated with $\eta_{j, i}$, and similarly let 
$$\overline{\op{wr}_{i,j}}\coloneqq  w_\Delta(\overline{\eta}_{i,j}) - w_\Delta(\overline{\eta}_{i,j+1})-w_\Delta(\overline{\eta}_{i+1,j})+w_\Delta(\overline{\eta}_{i+1,j+1})$$
Now ${(S_{u/v})}_{i,j} = \op{wr}_{j,i}+ \overline{\op{wr}_{i,j}}$.

\begin{figure}[!hb]
    \centering
    \begin{tikzpicture}
        \node at (0,2.5) {\includegraphics[height=0.5cm]{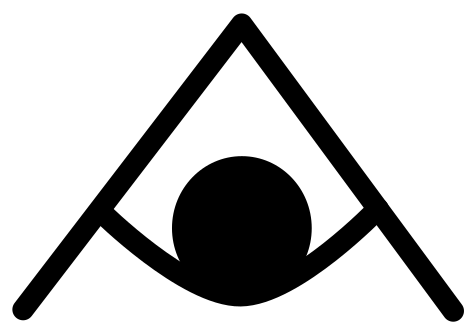}};
        \node at (0,-1) {\includegraphics[height=5cm]{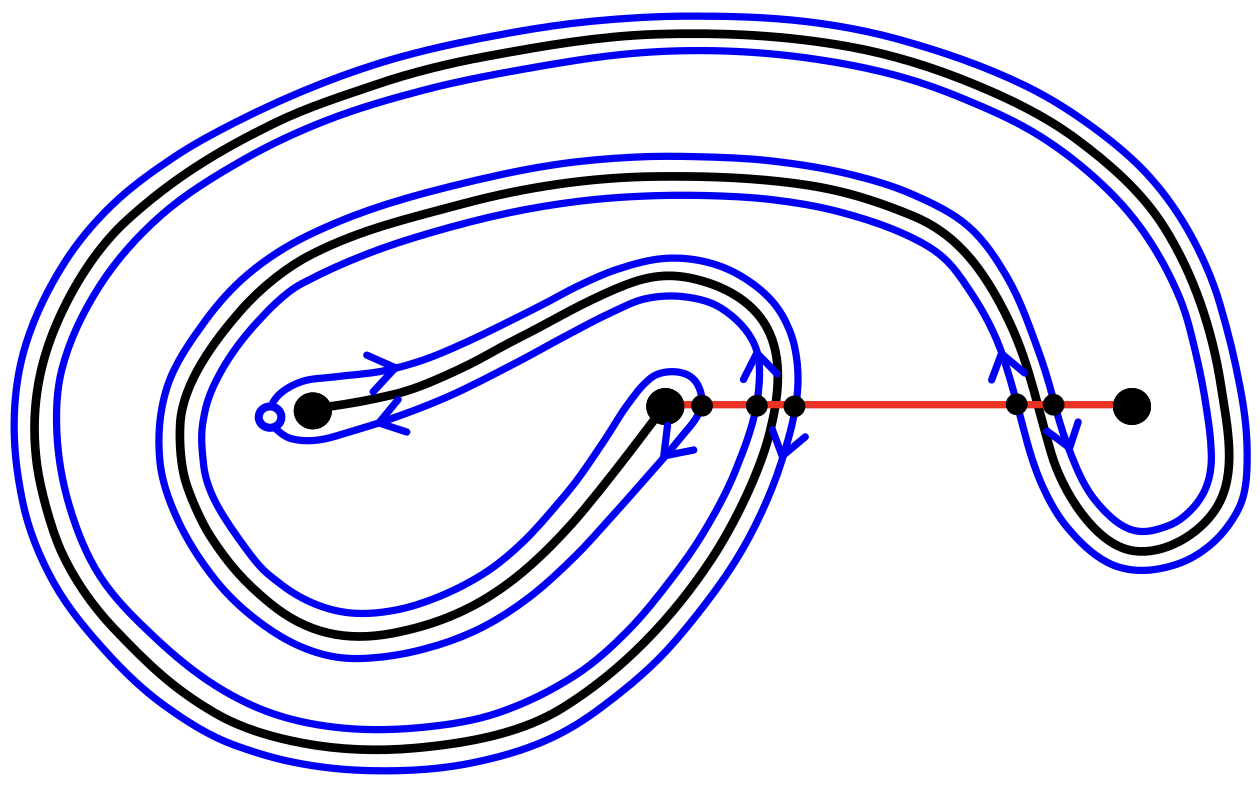}};
        \node at (8,0) {\includegraphics[height=7cm]{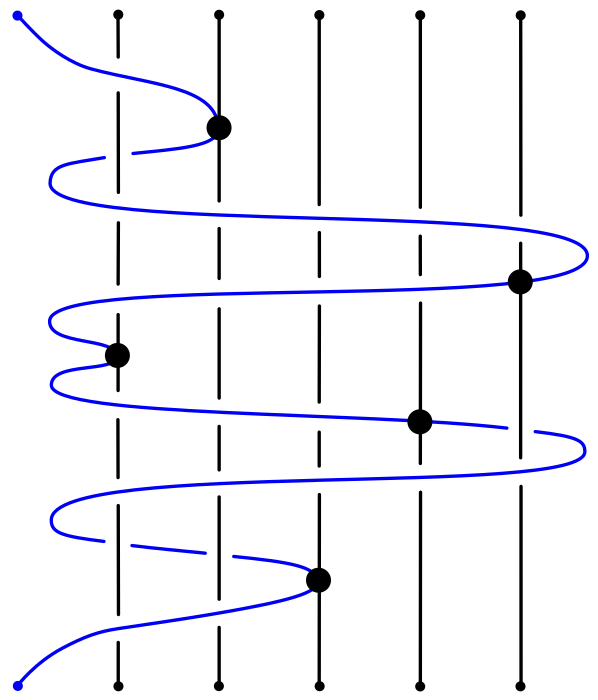}};
    \end{tikzpicture}
    \caption{Left: Lagrangians $\alpha$ and $\beta$ for $K_{5/2}$, with loop drawn for defining corresponding singular braid (and eyeball drawn to determine correct perspective). Right: Singular braid for $K_{5/2}$}
    \label{fig:singbraid}
\end{figure}

\begin{lemma}\label{writhe1}
    $\op{wr}_{j,i}$ can be expressed as a sum of terms $wr(p)$ ($p=i,i+1$) corresponding to the intersections of $\kappa_i$ and $\kappa_{i+1}$ with $\rho_{j,j+1}$. In other words, we add the term $wr(i)$ corresponding to $\kappa_i$ if and only if it lies on $\rho_{j,j+1}$, and similarly with $\kappa_{i+1}$. The terms $wr(p)$ are given below, with the direction being the direction $\rho_{j,j+1}$ passes through $\beta$ at $\kappa_p$. 
    \smallskip\\
    \resizebox{\textwidth}{!}{$
\begin{array}{c}
    p \ne j,j+1 \\[6pt]
    \begin{array}{|c|c|c|c|}
        \hline
        \mathrm{Point}\ \kappa_p & p\le v & \mathrm{Direction} & wr(p) \\ \hline\hline
        \kappa_i & \mathrm{True} & \mathrm{Down} & 0 \\ \hline
        \kappa_i & \mathrm{True} & \mathrm{Up} & 0 \\ \hline
        \kappa_i & \mathrm{False} & \mathrm{Down} & 1 \\ \hline
        \kappa_i & \mathrm{False} & \mathrm{Up} & -1 \\ \hline
        \kappa_{i+1} & \mathrm{True} & \mathrm{Down} & 2 \\ \hline
        \kappa_{i+1} & \mathrm{True} & \mathrm{Up} & -2 \\ \hline
        \kappa_{i+1} & \mathrm{False} & \mathrm{Down} & 1 \\ \hline
        \kappa_{i+1} & \mathrm{False} & \mathrm{Up} & -1 \\ \hline
    \end{array}
\end{array}
\qquad
\begin{array}{c}
    p=j \\[6pt]
    \begin{array}{|c|c|c|c|}
        \hline
        \mathrm{Point}\ \kappa_p & p\le v & \mathrm{Direction} & wr(p) \\ \hline\hline
        \kappa_i & \mathrm{True} & \mathrm{Down} & 0 \\ \hline
        \kappa_i & \mathrm{True} & \mathrm{Up} & 0 \\ \hline
        \kappa_i & \mathrm{False} & \mathrm{Down} & 1 \\ \hline
        \kappa_i & \mathrm{False} & \mathrm{Up} & 0 \\ \hline
        \kappa_{i+1} & \mathrm{True} & \mathrm{Down} & 1 \\ \hline
        \kappa_{i+1} & \mathrm{True} & \mathrm{Up} & -1 \\ \hline
        \kappa_{i+1} & \mathrm{False} & \mathrm{Down} & 0 \\ \hline
        \kappa_{i+1} & \mathrm{False} & \mathrm{Up} & -1 \\ \hline
    \end{array}
\end{array}
\qquad
\begin{array}{c}
    p=j+1 \\[6pt]
    \begin{array}{|c|c|c|c|}
        \hline
        \mathrm{Point} &p \le v & \mathrm{Direction} & wr(p) \\ \hline\hline
        \kappa_i & \mathrm{True} & \mathrm{Down} & 0 \\ \hline
        \kappa_i & \mathrm{True} & \mathrm{Up} & 0 \\ \hline
        \kappa_i & \mathrm{False} & \mathrm{Down} & 0 \\ \hline
        \kappa_i & \mathrm{False} & \mathrm{Up} & -1 \\ \hline
        \kappa_{i+1} & \mathrm{True} & \mathrm{Down} & 1 \\ \hline
        \kappa_{i+1} & \mathrm{True} & \mathrm{Up} & -1 \\ \hline
        \kappa_{i+1} & \mathrm{False} & \mathrm{Down} & 1 \\ \hline
        \kappa_{i+1} & \mathrm{False} & \mathrm{Up} & 0 \\ \hline
    \end{array}
\end{array}
$}
\end{lemma}

\begin{proof}

A \textit{singular braid} is an extension of classical braids where intersections of strands are allowed, and the intersections do not contribute to the writhe. Just as a $k$-strand braid is related to a loop in $\op{Conf}_k(\R^2)$, a $k$-strand singular braid is related to a loop in $\op{Sym}^k(\R^2)$.

We will abbreviate $\overline{\alpha_{u/v}}$ as $\alpha$. Define a singular braid associated with $\alpha$ by looking from above, as shown in Figure \ref{fig:singbraid}, with the strands being constant paths at each of the intersection points $\kappa_i$, and the path $\overline{\alpha_{u/v}}$, indicated in blue. Define $w_{i,j;k}$ as the writhe of the subbraid associated with the strands $\alpha$ and $\kappa_k$ restricted to the portion along $\rho_{i,j}$. This definition has the useful properties that $w_{i, j;k}+w_{j, l;k}=w_{i, l;k}$ and $w_\Delta(\eta_{j, i})=w_{j, i;i}$. Now, \begin{equation}\label{eq:wrsum}
\mathrm{wr}_{j,i} = w_{j, i;i} -w_{j+1, i;i}-w_{j, i+1;i+1}+w_{j+1, i+1;i+1}=w_{j,j+1;i}-w_{j,j+1;i+1}
\end{equation}

Because $i$ and $i+1$ are consecutive, if $\alpha$ passes above or below $\kappa_i$, it will pass above or below $\kappa_{i+1}$ as well, unless it intersects it. Likewise, if it passes above or below $\kappa_{i+1}$, it will pass above or below $\kappa_i$ as well, unless it intersects it. The contributions to $\op{wr}_{j,i}$ therefore come from the intersections of $\rho_{j,j+1}$ with $\kappa_i$ and $\kappa_{i+1}$. Let \[
wr(p) = \begin{cases} w_{j,j+1;i} & p=i\\ -w_{j,j+1;i+1} & p=i+1
\end{cases}
\]
Now the contributions of each intersection can be determined using images such as those in Figure \ref{fig:writhes}. 
\begin{figure}[htbp]
    \centering
    \begin{subfigure}{0.2\textwidth}
        \centering
        \includegraphics[height=2cm]{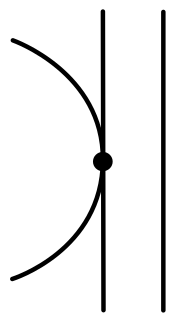}
        \caption*{$p=i,i\le v$}
    \end{subfigure}
    \hfill
    \begin{subfigure}{0.2\textwidth}
        \centering
        \includegraphics[height=2cm]{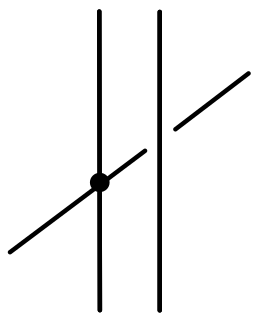}
        \caption*{$p={i}$, $i > v$, downwards}
    \end{subfigure}
    \hfill
    \begin{subfigure}{0.2\textwidth}
        \centering
        \includegraphics[height=2cm]{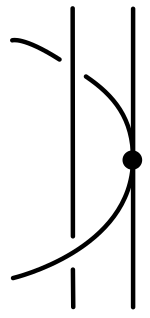}
        \caption*{$p={i+1}$, $i+1\le v$, downwards}
    \end{subfigure}
    \hfill
    \begin{subfigure}{0.2\textwidth}
        \centering
        \includegraphics[height=2cm]{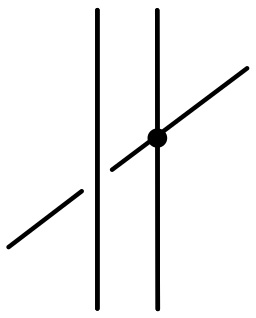}
        \caption*{$p={i+1}$, $i+1 > v$, upwards}
    \end{subfigure}

    \par\vspace{2em}
    
    \centering
    \begin{subfigure}{0.2\textwidth}
        \centering
        \includegraphics[height=2cm]{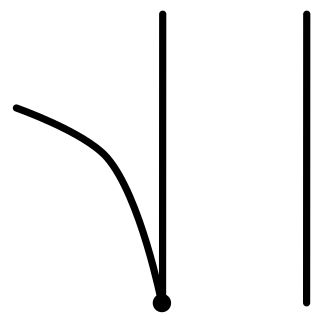}
        \caption*{$p=i=j, p\le v$, either direction}
    \end{subfigure}
    \hfill
    \begin{subfigure}{0.2\textwidth}
        \centering
        \includegraphics[height=2cm]{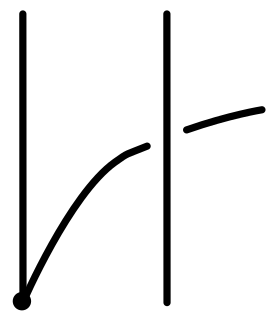}
        \caption*{$p=i+1=j$, $ p\le  v$, downwards}
    \end{subfigure}
    \hfill
    \begin{subfigure}{0.2\textwidth}
        \centering
        \includegraphics[height=2cm]{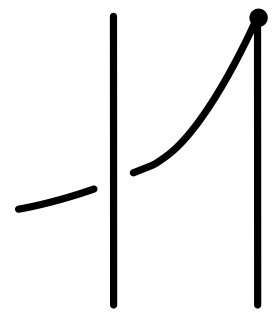}
        \caption*{$p=i+1=j+1$, $p\le v$, upwards}
    \end{subfigure}
    \hfill
    \begin{subfigure}{0.2\textwidth}
        \centering
        \includegraphics[height=2cm]{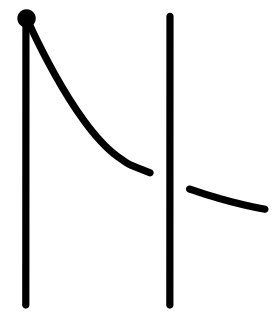}
        \caption*{$p=i=j+1$, $p > v$, upwards}
    \end{subfigure}
    \caption{Diagrams of braids near intersection points}
    \label{fig:writhes}
\end{figure}
\end{proof}

\begin{lemma}\label{writhe2}
    $\overline{\op{wr}_{i,j}}$ can be expressed as a sum of terms $\overline{wr}(p)$ ($p=i,i+1$) corresponding to the intersections of $\kappa_i$ and $\kappa_{i+1}$ with $\rho_{j,j+1}$. The terms $\overline{wr}(p)$ are given below, with the direction being the direction $\rho_{j,j+1}$ passes through $\beta$ at $\kappa_p$.
\smallskip\\
\resizebox{\textwidth}{!}{$
\begin{array}{c}
    p \ne j,j+1 \\[6pt]
    \begin{array}{|c|c|c|c|}
        \hline
        \mathrm{Point}\ \kappa_p & p \le v & \mathrm{Direction} & \overline{wr}(p) \\ \hline\hline
        \kappa_i& \mathrm{True} & \mathrm{Down} & 1 \\ \hline
        \kappa_i& \mathrm{True} & \mathrm{Up} & -1 \\ \hline
        \kappa_{i+1}& \mathrm{True} & \mathrm{Down} & -1 \\ \hline
        \kappa_{i+1}& \mathrm{True} & \mathrm{Up} & 1 \\ \hline
        \mathrm{Any}&\mathrm{False} & \mathrm{Any} & 0 \\ \hline
    \end{array}
\end{array}
\qquad
\begin{array}{c}
    p=j \\[6pt]
    \begin{array}{|c|c|c|c|}
        \hline
        \mathrm{Point}\ \kappa_p & p \le v & \mathrm{Direction} & \overline{wr}(p) \\ \hline\hline
        \kappa_i & \mathrm{True} & \mathrm{Down} & 1 \\ \hline
        \kappa_i & \mathrm{True} & \mathrm{Up} & 0 \\ \hline
        \kappa_{i+1} & \mathrm{True} & \mathrm{Down} & -1 \\ \hline
        \kappa_{i+1} & \mathrm{True} & \mathrm{Up} & 0 \\ \hline
        \mathrm{Any} & \mathrm{False} & \mathrm{Any} & 0 \\ \hline
    \end{array}
\end{array}
\qquad
\begin{array}{c}
    p=j+1 \\[6pt]
    \begin{array}{|c|c|c|c|}
        \hline
        \mathrm{Point}\ \kappa_p & p \le v & \mathrm{Direction} & \overline{wr}(p) \\ \hline\hline
        \kappa_i & \mathrm{True} & \mathrm{Down} & 0 \\ \hline
        \kappa_i & \mathrm{True} & \mathrm{Up} & -1 \\ \hline
        \kappa_{i+1} & \mathrm{True} & \mathrm{Down} & 0 \\ \hline
        \kappa_{i+1} & \mathrm{True} & \mathrm{Up} & 1 \\ \hline
        \mathrm{Any} & \mathrm{False} & \mathrm{Any} & 0 \\ \hline
    \end{array}
\end{array}
$}
\end{lemma}

\begin{proof}

Consider the term $w_\Delta(\overline{\eta}_{i,j})$. There are 8 possible cases, depending on whether or not $j<i$, the direction $\eta_{j,i}$ passes through $\beta$ at $\kappa_i$, and whether or not $i\le v$. The 8 cases, as well as the corresponding braids, are given in Figure \ref{fig:braids}.

\begin{figure}[htbp]
    \centering
    \begin{subfigure}{0.25\textwidth}
        \centering
        \includegraphics[height=2cm]{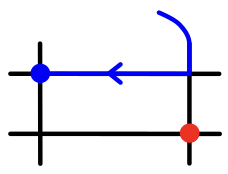}%
        \includegraphics[height=2cm]{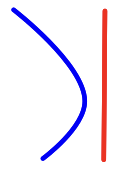}
        \caption*{$j<i$, downwards}
    \end{subfigure}
    \hfill
    \begin{subfigure}{0.25\textwidth}
        \centering
        \includegraphics[height=2cm]{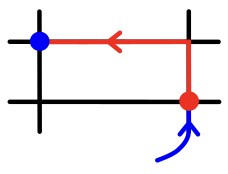}%
        \includegraphics[height=2cm]{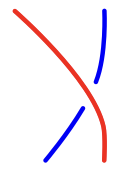}
        \caption*{$j<i$, $i\le v$, upwards}
    \end{subfigure}
    \hfill
    \begin{subfigure}{0.25\textwidth}
        \centering
        \includegraphics[height=2cm]{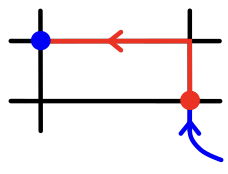}%
        \includegraphics[height=2cm]{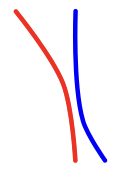}
        \caption*{$j<i$, $i \ge v$, upwards}
    \end{subfigure}
    
    \par\vspace{2em}
    
    \centering
    \begin{subfigure}{0.25\textwidth}
        \centering
        \includegraphics[height=2cm]{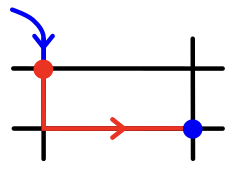}%
        \includegraphics[height=2cm]{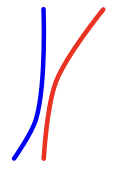}
        \caption*{$j>i$, downwards}
    \end{subfigure}
    \hfill
    \begin{subfigure}{0.25\textwidth}
        \centering
        \includegraphics[height=2cm]{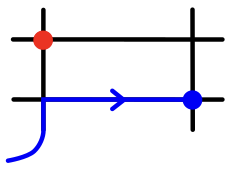}%
        \includegraphics[height=2cm]{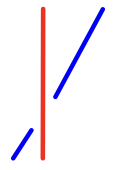}
        \caption*{$j>i$, $i\le v$, upwards}
    \end{subfigure}
    \hfill
    \begin{subfigure}{0.25\textwidth}
        \centering
        \includegraphics[height=2cm]{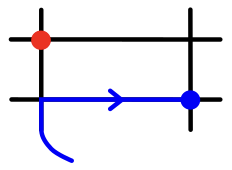}%
        \includegraphics[height=2cm]{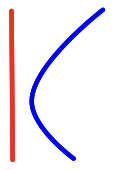}
        \caption*{$j>i$, $i \ge v$, upwards}
    \end{subfigure}
    \caption{Paths of $\eta_{j,i}$ and corresponding braids}
    \label{fig:braids}
\end{figure}
From the braids, the writhe contributed by $\overline{\eta}_{i,j}$ is as follows:
\[
\begin{array}{|c|c|c|c|}
    \hline i \le v & \trm{Direction} & \trm{Writhe}
    \\ 
     \hline \hline
    \trm{True} & \trm{Down} & 0 \\ \hline
    \trm{True}& \trm{Up} & -1 \\ \hline
     \trm{False} & \trm{Down} & 0 \\ \hline
    \trm{False} & \trm{Up} & 0 \\ \hline
 
\end{array}
\]

Consider the sum $\overline{wr}(i)= w_\Delta (\overline{\eta}_{i,j})-w_{\Delta}(\overline{\eta}_{i,j+1})$. If $\rho_{j,j+1}$ does not intersect $\kappa_i$, then the direction of $\eta_{j,i}$ through $\beta$ at $\kappa_i$ is the same as the direction of $\eta_{j+1,i}$ through $\beta$ at $\kappa_i$, otherwise the directions will be opposite to each other. If the directions are the same, the terms must cancel, and so there is no contribution to $\overline{wr_{i,j}}$.

Also, if $i=j$ or $i=j+1$, one of the terms will vanish, so the terms may not cancel. We can still consider the direction which the path from $\kappa_j$ to $\kappa_{j+1}$ passes through $\beta$ at $\kappa_i$, and can then determine $\overline{wr}(i) $ by the intersection as follows:

\[
\begin{array}{c}
    i \ne j,j+1 \\[6pt]
    \begin{array}{|c|c|c|}
        \hline
        i \le v & \mathrm{Direction} & \overline{wr}(i) \\ \hline\hline
        \mathrm{True} & \mathrm{Down} & 1 \\ \hline
        \mathrm{True} & \mathrm{Up} & -1 \\ \hline
        \mathrm{False} & \mathrm{Any} & 0 \\ \hline
    \end{array}
\end{array}
\qquad
\begin{array}{c}
    i=j \\[6pt]
    \begin{array}{|c|c|c|}
        \hline
        i \le v & \mathrm{Direction} & \overline{wr}(i) \\ \hline\hline
        \mathrm{True} & \mathrm{Down} & 1 \\ \hline
        \mathrm{True} & \mathrm{Up} & 0 \\ \hline
        \mathrm{False} & \mathrm{Any} & 0 \\ \hline
    \end{array}
\end{array}
\qquad
\begin{array}{c}
    i=j+1 \\[6pt]
    \begin{array}{|c|c|c|}
        \hline
        i \le v & \mathrm{Direction} &\overline{wr}(i) \\ \hline\hline
        \mathrm{True} & \mathrm{Down} & 0 \\ \hline
        \mathrm{True} & \mathrm{Up} & -1 \\ \hline
        \mathrm{False} & \mathrm{Any} & 0 \\ \hline
    \end{array}
\end{array}
\]

The value of $\overline{wr}(i)$ is therefore determined by the intersection of $\rho_{j,j+1}$ with $\kappa_i$. Let $\overline{wr}(i+1)= -\big(w_\Delta (\overline{\eta}_{i+1,j})-w_{\Delta}(\overline{\eta}_{i+1,j+1})\big)$, so that \[
\overline{\op{wr}_{i,j}} = \overline{wr}(i) + \overline{wr}(i+1).
\]
Because of the negative sign in the definition, $\overline{wr}(i+1)$ is given by the table above, but with the signs reversed.
\end{proof}

Now we can prove Theorem \ref{skeletontheorem}. The tables from Lemmas \ref{writhe1} and \ref{writhe2} can be combined and simplified. If $p=i,i+1$, the total writhe contributed to $w_\Delta(\gamma_{j, i}) - w_\Delta(\gamma_{j+1, i})-w_\Delta(\gamma_{j, i+1})+w_\Delta(\gamma_{j+1, i+1})$ by the intersection of $\kappa_p$ with $\rho_{j,j+1}$ is given by the following tables.\\
\resizebox{\textwidth}{!}{$
\begin{array}{c}
    p\ne j,j+1 \\[6pt]
    \begin{array}{|c|c|}
        \hline
        \mathrm{Direction} & wr(p) + \overline{wr}(p) \\ \hline\hline
        \mathrm{Down} & 1 \\ \hline
        \mathrm{Up} & -1 \\ \hline
    \end{array}
\end{array}
\quad
\begin{array}{c}
    p=j \\[6pt]
    \begin{array}{|c|c|c|}
        \hline
        \mathrm{Point}\ \kappa_p & \mathrm{Direction} & wr(p) + \overline{wr}(p) \\ \hline\hline
        \kappa_i & \mathrm{Down} & 1 \\ \hline
        \kappa_i & \mathrm{Up} & 0 \\ \hline
        \kappa_{i+1} & \mathrm{Down} & 0 \\ \hline
        \kappa_{i+1} & \mathrm{Up} & -1 \\ \hline
    \end{array}
\end{array}
\quad
\begin{array}{c}
    p=j+1 \\[6pt]
    \begin{array}{|c|c|c|}
        \hline
        \mathrm{Point}\ \kappa_p & \mathrm{Direction} & wr(p) + \overline{wr}(p) \\ \hline\hline
        \kappa_i & \mathrm{Down} & 0 \\ \hline
        \kappa_i & \mathrm{Up} & -1 \\ \hline
        \kappa_{i+1} & \mathrm{Down} & 1 \\ \hline
        \kappa_{i+1} & \mathrm{Up} & 0 \\ \hline
    \end{array}
\end{array}
$}
\smallskip

Observe that $\alpha$ is a loop around an arc associated with the tangle $\tau_{u/v}$ (the black arc in Figure \ref{fig:singbraid}). Orient $\alpha$ as in Figure \ref{fig:singbraid}, with the tangle to the right. As the tangle always stays to the right of $\alpha$, $\alpha$ passes through $\kappa_i$ and $\kappa_{i+1}$ in opposite directions. This remains true if we reverse the orientation of $\alpha$. 

Now, if $i=j$, then $i+1=j+1$. Suppose the direction $\rho_{j,j+1}$ passes through $\beta$ at $\kappa_i$ is upwards. Then necessarily the direction $\rho_{j,j+1}$ passes through $\beta$ at $\kappa_{i+1}$ is downwards, so we have ${(S_{u/v})}_{i,i}=0+1=1$. If the direction $\rho_{j,j+1}$ passes through $\beta$ at $\kappa_i$ is downwards, then necessarily the direction $\rho_{j,j+1}$ passes through $\beta$ at $\kappa_{i+1}$ is upwards, so we have ${(S_{u/v})}_{i,i}=1+0=1$. Therefore, the diagonal of $S_{u/v}$ is identically $1$.

Now we prove that $-1\le {(S_{u/v})}_{i,j}\le 1$ in general. All the terms appearing in the tables are between \(-1\) and \(1\), so the only way this could fail to be true is if $\rho_{j,j+1}$ intersects both $\kappa_i$ and $\kappa_{i+1}$. But if $\rho_{j,j+1}$ intersects both $\kappa_i$ and $\kappa_{i+1}$, the direction of intersection at $\kappa_i$ and $\kappa_{i+1}$ will be opposite, so any nonzero terms will cancel out.
\end{proof}

\subsection{Shape of the Tail}\label{subsec:shape}

When computing and visualizing the colored HOMFLY-PT polynomials, the shape of the plotted $a$- and $q$-powers seems to be bounded by a collection of curves, as in Figure \ref{fig:K5/2} (see Conjecture \ref{limshapeconj} for a precise statement). In fact, when considering just the tail or head of the colored HOMFLY-PT polynomial of a knot in class $\mathcal Q$, the shape is bounded by a quadratic curve. 
\begin{proposition}\label{prop:shape}
    Let $K$ be a knot of class $\mathcal{Q}$, with quiver $Q$ and vectors $A,S,R$, and let \[\Psi_K(a,q)=\sum\limits_{m\in\Z, n\ge 0}c_{m,n}a^{m}q^{n}\] be the tail of the HOMFLY-PT polynomial of $K$. Then there exists a quadratic function $f(m)$ such that $n\ge f(m)$ for all $n$ with $c_{m,n}\ne 0$.
\end{proposition}
\begin{proof} Recall that for a knot of class $\mathcal{Q}$, $Q^t$ is nonnegative and nonzero along the diagonal. If $Q^t$ is empty, then $\Psi_K(a,q)=1$, so suppose $Q^t$ is a $k\times k$ matrix with $k>0$. Let $\mu=\min Q^t_{i,i}>0$ and $\lambda=\max|A^t_i|$. If $A^t = 0$, then $c_{m,n}=0$ whenever $m\neq0$, so suppose $A^t\neq0$. Then $\lambda >0$. By \eqref{homftail}, up to an overall shift we have \[\Psi_K(a,q) = \sum_{v\in\N^k}(-1)^{R^t\cdot v}a^{A^t\cdot v}q^{v^T\cdot Q^t\cdot v+S\cdot v}\inv{\prod_{i=1}^k \poch{q^2}{q^2}{v_i}}. \]

Suppose $v\in \N^k$. As $\inv{\prod_{i=1}^k \poch{q^2}{q^2}{v_i}}$ is a power series with constant coefficient $1$ (which follows from Lemma \ref{multinomlim}), the lowest $q$-degree in the corresponding summand is then $n=v^T\cdot Q^t\cdot v+S^t\cdot v$. The corresponding $a$-degree in the summand is $m=A^t\cdot v$. Since $K$ is of class $\mathcal{Q}$, $S^t=\text{diag}(Q^t)+2A^t$, we can write $n=v^T\cdot Q^t\cdot v+\text{diag}(Q^t)\cdot v+2m$. Then
\begin{align*}
    n&=\sum_{i j}Q^t_{ij}v_iv_j+\sum_{i=1}^{k}Q_{ii}^tv_i+2m\\
    &\ge \mu\sum_{i=1}^{k}v_i^2+\mu\sum_{i=1}^{k}v_i+2m.
\end{align*}
The Cauchy-Schwarz inequality gives
$$m^2=(A^t\cdot v)^2\leq\|A^t\|_2^2\|v\|_2^2\leq k\lambda^2\sum_{i=1}^{k}v_i^2.$$
Also, $m=A^t\cdot v\geq -\lambda \sum_{i=1}^{k}v_i$, so
$$n\geq \frac{\mu}{k\lambda^2}m^2+(2-\frac{\mu}{\lambda})m.$$
Therefore, $c_{m,n}=0$ for all $n<\frac{\mu}{k\lambda^2}m^2+(2-\frac{\mu}{\lambda})m$.
\end{proof}

Note that this quadratic bound is on the unnormalized tail, so a bound on the normalized tail will require a substitution $m\to m-\alpha$, $n\to n-\beta$, where $\alpha,\beta$ are such that $a^\alpha q^\beta$ is the monomial in the unnormalized tail with the lowest $q$-power and highest $a$-power. 
\begin{example}
We compute the quadratic bound of Proposition \ref{prop:shape} for the knot $K_{5/2}$. Recall in Example \ref{ex:5/2}, we computed the subquiver data: 
$$Q^t=(1),\ A^t=(-2),\ S^t=(-3).$$
Then $k=1$, $\mu=1$, and $\lambda=2.$ Substituting these values into the bound of Proposition \ref{prop:shape} gives 
$$n\geq \frac{1}{4}m^2+\frac{3}{2}m.$$
Hence, the coefficients $c_{m,n}$ of $\Psi_{K_{5/2}}$ vanish whenever $n<\frac{1}{4}m^2+\frac{3}{2}m$. We visualize this in Figure \ref{fig:bound}, where we can see that the bound is tight for this knot.  
\begin{figure}
    \centering
    \includegraphics[width=0.5\linewidth]{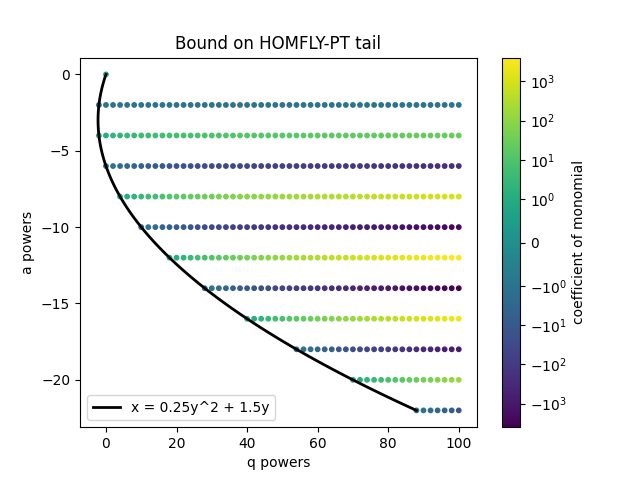}
    \caption{Quadratic bound on $\Psi_{K_{5/2}}$}
    \label{fig:bound}
\end{figure}
\end{example}

\section{Rational Knots are of class \texorpdfstring{$\mathcal{Q}$}{Q}}\label{sec:Q}

In this section, we will prove the following theorem, which motivates the definition of class $\mathcal{Q}$:
\begin{theorem}\label{rationalsareq}
    All rational knots are of class $\mathcal{Q}$.
\end{theorem}

Our strategy relies on the inductive construction of the quiver in the proof of Theorem 4.1 of  Sto\v{s}i\'c--Wedrich \cite{Sto_i__2019}. The Knots-Quivers Correspondence has an extension to tangles, where the individual colored HOMFLY-PT invariants are no longer polynomials, but skein module evaluations written as linear combinations of simple basis webs. Importantly, this extension to tangles can still be understood in terms of linear and quadratic forms analogous to the ones used for knots. In particular, if $\langle \tau_{u/v}\rangle_j$ is the $j$-colored skein module evaluation for the rational tangle $\tau_{u/v}$, then the Tangles-Quivers Correspondence gives that
\begin{multline}\label{eq:tangle_polynomial}
    \sum_{j\geq 0}\langle\tau_{u/v}\rangle_j=\sum_{\textbf{d}=(d_1,...,d_{k+l})\in\mathbb{N}^{k+l}}(-q)^{S\cdot \textbf{d}}a^{A\cdot \textbf{d}}q^{\textbf{d}^t\cdot Q\cdot \textbf{d}}(q^2;q^2)_{K\cdot \textbf{d}}\qmultinom{d_1+...+d_k}{d_1,\dots,d_k}{q^2}\qmultinom{d_{k+1}+...+d_{k+l}}{d_{k+1},\dots,d_{k+l}}{q^2}\\
    \times X[d_1+...+d_{k+l},d_1+...+d_k],
\end{multline}
where $S,A,$ and $Q$ are linear and quadratic forms as in the case for knots, $K$ is an additional linear form of $0$'s and $1$'s, and $X[j,k]$ denotes the basis web ($X$ is the orientation, $j$ is the color, and $k$ is called the \textit{weight} of the web).

An important distinction from the Knots-Quivers Correspondence is that the free $\mathbb{Z}$-module on which the linear and quadratic forms are defined splits into \textit{active} and \textit{inactive} blocks. Algebraically, summation indices coming from the active block, or the \textit{active indices}, contribute to the weight of the basis webs, whereas the indices coming from the inactive block, or the \textit{inactive indices}, do not. Geometrically, the active and inactive indices come from two different types of Lagrangian intersections in \cite{higgins1}. We will not elaborate too much on the Tangles-Quivers Correspondence here, but the interested reader should refer to \cite{Sto_i__2019} or \cite{higgins1} for more details. We will only need the fact that the quivers can be computed for rational tangles via the sequence of twists used to build them, and if $K_{u/v}=Cl(\tau_{u/v})$, then one can obtain quivers for the knot $K_{u/v}$ from quivers for $\tau_{u/v}$.

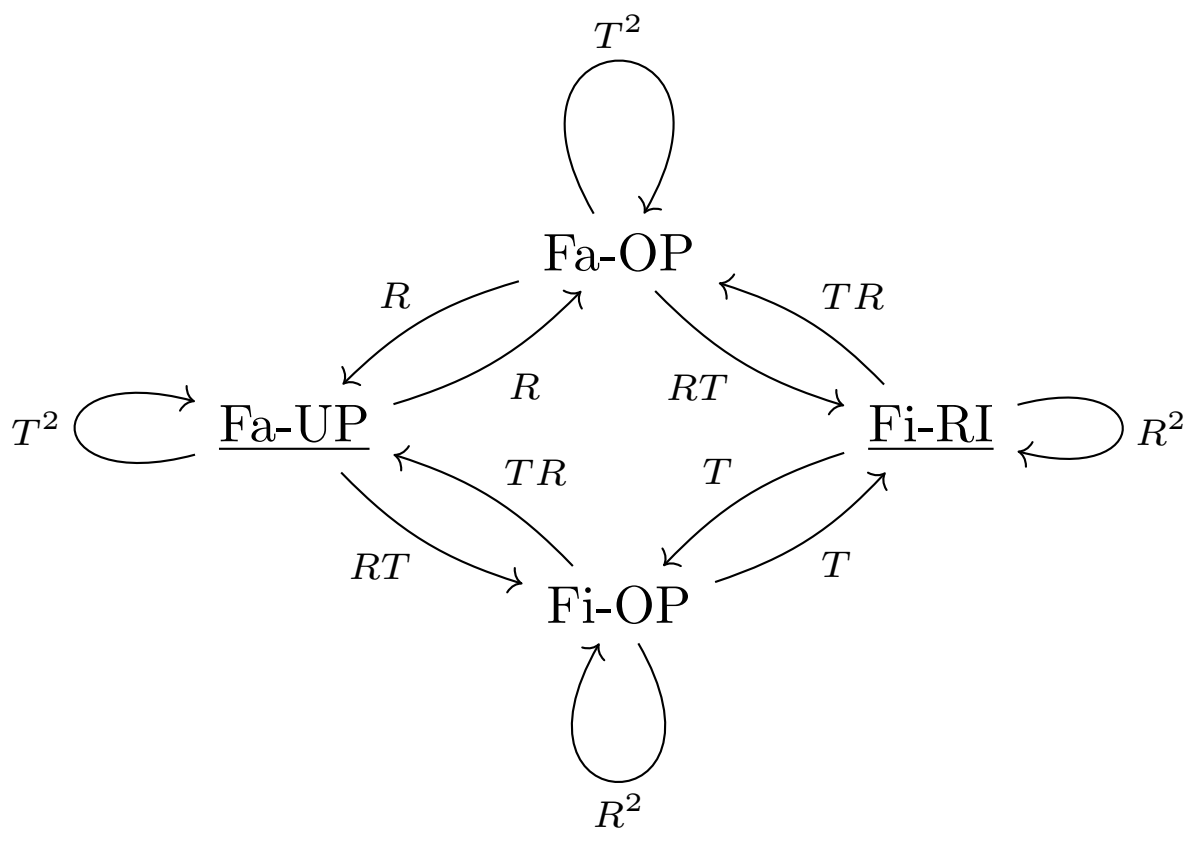
\begin{figure}
    \centering
    \begin{tikzcd}
    & \text{Fa-OP} \arrow[dl,bend right=0.5cm,"R"'] \arrow[dr, bend right=0.5cm,"RT"']\arrow[loop above,looseness=12, out=120, in=60,"T^2"] & \\
    \underline{\text{Fa-UP}} \arrow[ur, bend right=0.5cm,"R"']\arrow[loop left,"T^2"]\arrow[dr, bend right=0.5cm,"RT"'] && \underline{\text{Fi-RI}} \arrow[ul, bend right=0.5cm,"TR"']\arrow[loop right,"R^2"]\arrow[dl, bend right=0.5cm,"T"']\\
    & \text{Fi-OP} \arrow[ul, bend right=0.5cm,"TR"']\arrow[ur,bend right=0.5cm,"T"']\arrow[loop below,looseness=14, out=-60, in=-120,"R^2"] &\\
\end{tikzcd}
    \caption{Graph depicting the four tangle types and legal twist operations between them. Types where the $Cl(T-)$ operation can be applied are underlined.}
    \label{tangletypes}
\end{figure}

Let $u/v>1$ be reduced and suppose that $u$ is odd, so that the numerator closure of $\tau_{u/v}$ (if allowed by its orientation) is a rational knot $K_{u/v}$.  Its continued-fraction expansion determines a word in top and right twists, denoted by $T$ and $R$, ending with the final top twist used in the knot construction.  First, observe that $K_{u/v}=Cl(\tau_{u/v})=Cl(T\tau_{(u-v)/v})$; following Proposition~4.4 of \cite{Sto_i__2019}, we do not process this word for $\tau_{(u-v)/v}$ one letter at a time.  Instead, 
 Sto\v{s}i\'c--Wedrich prove that $\tau_{(u-v)/v}$ can be built inductively using only four intermediate types of tangles, which we call Fa-UP, Fa-OP, Fi-OP, and Fi-RI, where the $a$ and $i$ refer to a property involving the linear form $K$ involving the active and inactive indices, respectively, and $UP$, $OP$, and $RI$ refer to the orientations of the intermediate tangles. At each step, the following table shows which twists can be performed to stay within one of the four types. See Figure \ref{tangletypes} for a graphical depiction of the tangle types and the twist operations between them.
\[
\begin{array}{c|c}
\text{current type} & \text{legal next operation}\\ \hline
\Fa\text{-UP or }\Fa\text{-OP} & T^2,\ RT,\ \text{or }R\\
\Fi\text{-OP or }\Fi\text{-RI} & R^2,\ TR,\ \text{or }T
\end{array}
\]
Crucially, for the numerator closure of $\tau_{u/v}$ to be a knot, $\tau_{(u-v)/v}$ must be of type $\Fa\text{-UP or }\Fi\text{-RI}$ by Proposition 4.4 of \cite{Sto_i__2019}. These two tangle types are underlined in Figure \ref{tangletypes}.

\begin{remark}
    \cite{Sto_i__2019} defines $\mathcal{F}_a(\alpha,\beta,X)$ and $\mathcal{F}_i(\alpha,\beta,X)$, which denote sets of generating functions of the form of \eqref{eq:tangle_polynomial}. For our purposes, we only care that these states give us rules for how to compute the intermediate quivers up to $\tau_{(u-v)/v}$. Hence, we abbreviate this notation by only writing $\Fa\text{-X}$ and $\Fi\text{-X}$.
\end{remark}

At every intermediate stage, the quiver of a rational tangle can be written in block matrix form:
\[
Q=
\left(
\begin{array}{c||c}
Q_{++} & Q_{+-}\\ \hline \hline
Q_{-+} & Q_{--}
\end{array}
\right),
\qquad
Q_{-+}=Q_{+-}^{\mathsf T},
\]
where $Q_{++}$ denotes the restriction of $Q$ to the active indices, $Q_{--}$ to the inactive indices, and the other two blocks to mixed indices. These blocks can be computed by inductive rules based on the sequence of top and right twists used to construct the tangle, with base case given by the quiver for the trivial tangle $\tau_{0/1}$, or
\[Q_{\tau_{0/1}} = (0)\]
taken as the $Q_{--}$ block (with all others being empty). The inductive rules for computing the quiver of the tangle will be given in the proofs of the upcoming lemmas.

To prove Theorem \ref{rationalsareq}, we also need formulas for the quivers we obtain for $K_{u/v}$ by applying $Cl(T-)$ to $\tau_{(u-v)/v}$ when $\tau_{(u-v)/v}$ is of type Fa-UP or Fi-RI. These are given by the following lemma, which was originally proved in \cite{Sto_i__2019}, but we follow the framing conventions of \cite{higgins2}.

\begin{lemma}\label{lem:closure}
    Suppose $K_{u/v}=Cl(T\tau_{(u-v)/v})$. If $\tau_{(u-v)/v}$ is of type $\Fa\text{-UP}$, then the $Cl(T-)$ operation is given on matrix blocks by
    \begin{equation}
      Q_{\tau_{(u-v)/v}}=
\left(
\begin{array}{c||c}
Q_{++} & Q_{+-}\\ \hline \hline
Q_{-+} & Q_{--}
\end{array}
\right) \xrightarrow{\quad Cl(T-)\quad} 
Q_{K_{u/v}}=\begin{pmatrix}
    Q_{++}+2 & Q_{++}+L & Q_{+-}\\
    Q_{++}+U & Q_{++}-1 & Q_{+-}-1\\
    Q_{-+} & Q_{-+}-1 & Q_{--}
    \end{pmatrix}.
    \end{equation}
   If $\tau_{(u-v)/v}$ is of type $\Fi$-RI, then the $Cl(T-)$ operation is given by
   \begin{equation}
      Q_{\tau_{(u-v)/v}}=
\left(
\begin{array}{c||c}
Q_{++} & Q_{+-}\\ \hline \hline
Q_{-+} & Q_{--}
\end{array}
\right) \xrightarrow{\quad Cl(T-)\quad} 
Q_{K_{u/v}}=\begin{pmatrix}
    Q_{++} & Q_{+-} & Q_{+-}-1\\
    Q_{-+} & Q_{--}+2 & Q_{--}+L\\
    Q_{-+}-1 & Q_{--}+U & Q_{--}-1
    \end{pmatrix},
    \end{equation}
    where $U$ and $L$ are the strictly upper/lower triangular matrices of $1$s.
\end{lemma}

Now, we state the key lemmas that will be used to prove Theorem \ref{rationalsareq}. All matrix inequalities below are defined entry-wise; for example, $Q_{+-}\leq a$ means that every entry of $Q_{+-}$ is at most $a$. We say that a matrix has a \textit{unique diagonal maximum} if the diagonal of the matrix has a unique maximum which is also the unique maximum of the matrix, and that it has a \textit{diagonal minimum} if the diagonal has a unique minimum which is also the minimum of the matrix.

\begin{lemma}\label{maxineqs}
Every intermediate tangle satisfies the following inequalities:
\begin{itemize}
    \item 
If the tangle is of type $\Fa$-UP, then either $Q_{++}$ has a unique diagonal maximum $\lambda$, and $Q_{--}$ has a unique diagonal maximum of value $\eta$, such that
\begin{equation}\label{faup1}
    Q_{+-}\leq \lambda+1,
\qquad
\eta\leq \lambda+1.
\end{equation}
or $Q_{--}$ has a unique diagonal maximum $\lambda$, and
\begin{equation}\label{faup2}
Q_{++}\leq \lambda-3,
\qquad
Q_{+-}\leq \lambda-1.
\end{equation}

\item
If the tangle is of type $\Fa$-OP, then $Q_{--}$ has a unique diagonal maximum $\lambda$, and
\begin{equation}\label{faopmax}
Q_{++}\leq \lambda-2,
\qquad
Q_{+-}\leq \lambda-1.
\end{equation}

\item 
If the tangle is of type $\Fi$-OP, then $Q_{--}$ has a unique diagonal maximum $\lambda$, and
\begin{equation}\label{fiopmax}
Q_{++}\leq \lambda,
\qquad
Q_{+-}\leq \lambda.
\end{equation}
\item 
If the tangle is of type $\Fi$-RI, then $Q_{++}$ has a unique diagonal maximum of value $\lambda$, $Q_{--}$ has a unique diagonal maximum of value $\lambda-1$, and
\begin{equation}\label{firimax}
Q_{+-}\leq \lambda.
\end{equation}
\end{itemize}
\end{lemma}

\begin{lemma}\label{minineqs}
    Every intermediate tangle satisfies the following inequalities:
    \begin{itemize}
        \item If the tangle is of type $\Fa$-UP, then $Q_{--}$ has a diagonal minimum $\lambda$, $Q_{++}$ has the same diagonal minimum $\lambda$ (if nonempty), and 
        \begin{equation}\label{faupmin}
        Q_{+-}\ge \lambda
        \end{equation}
        \item If the tangle is of type $\Fa$-OP or $\Fi$-OP, then $Q_{++}$ has a diagonal minimum $\lambda$ and $Q_{--}$ has a diagonal minimum $\eta$ such that 
        \begin{equation}\label{faopmin}
        Q_{+-}\ge \lambda+1, \qquad \eta\ge\lambda+1
        \end{equation}
        \item If the tangle is of type $\Fi$-RI, then $Q_{++}$ has a diagonal minimum $\lambda$ and $Q_{--}$ has a diagonal minimum $\eta$ such that either 
        \begin{equation}\label{firi1}
Q_{+-}\ge \lambda+1, \qquad \eta\ge\lambda+2
        \end{equation}
        or
        \begin{equation}\label{firi2}
        Q_{+-}\ge \eta,\qquad \eta\le\lambda
        \end{equation}
    \end{itemize}
\end{lemma}

From the lemmas, the proof of Theorem \ref{rationalsareq} proceeds as follows:

\begin{proof}[Proof of Theorem \ref{rationalsareq}]
    Let $K_{u/v}$ be a rational knot. Then $K_{u/v} = Cl(T\tau_{(u-v)/v})$ where $\tau_{(u-v)/v}$ must be of type $\Fa$-UP or $\Fi$-RI.
    Additionally, $\tau_{(u-v)/v}$ has a quiver of the form:
    \[
    Q_{\tau_{(u-v)/v}}=
    \left(
    \begin{array}{c||c}
    Q_{++} & Q_{+-}\\ \hline \hline
    Q_{-+} & Q_{--}
    \end{array}
    \right).
    \]
    
    If $\tau_{(u-v)/v}$ is of type $\Fa$-UP, then by Lemma $\ref{lem:closure}$, the quiver of $K_{u/v}$ is: 
    \[
    Q_{K_{u/v}}=
    \begin{pmatrix}
    Q_{++}+2 & Q_{++}+L & Q_{+-}\\
    Q_{++}+U & Q_{++}-1 & Q_{+-}-1\\
    Q_{-+} & Q_{-+}-1 & Q_{--}
    \end{pmatrix}.
    \]
    By Lemma \ref{minineqs}, $Q_{K_{u/v}}$ has a diagonal minimum coming from the $Q_{++}-1$ block.
    
    By Lemma \ref{maxineqs}, if $Q_{\tau_{(u-v)/v}}$ satisfies \eqref{faup1}, then $Q_{K_{u/v}}$ has a unique diagonal maximum coming from the $Q_{++}+2$ block, and if $Q_{\tau_{(u-v)/v}}$ satisfies \eqref{faup2}, then $Q_{K_{u/v}}$ has a unique diagonal maximum coming from the $Q_{--}$ block.

    If $\tau_{(u-v)/v}$ is of type $\Fi$-RI, then by Lemma $\ref{lem:closure}$, the quiver of $K_{u/v}$ is:
    \[
    Q_{K_{u/v}}=
    \begin{pmatrix}
    Q_{++} & Q_{+-} & Q_{+-}-1\\
    Q_{-+} & Q_{--}+2 & Q_{--}+L\\
    Q_{-+}-1 & Q_{--}+U & Q_{--}-1
    \end{pmatrix}.
    \]
    By Lemma \ref{minineqs}, if $Q_{\tau_{(u-v)/v}}$ satisfies \eqref{firi1}, $Q_{K_{u/v}}$ has a diagonal minimum coming from the $Q_{++}$ block, and if $Q_{\tau_{(u-v)/v}}$ satisfies \eqref{firi2}, $Q_{K_{u/v}}$ has a diagonal minimum coming from the $Q_{--} - 1$ block.
    
    By Lemma \ref{maxineqs}, we have that $Q_{K_{u/v}}$ has a unique diagonal maximum coming from the $Q_{--}+2$ block. 
    
    Therefore, for all rational knots $K_{u/v}$ with $u/v > 1$, $Q_{K_{u/v}}$ has a diagonal minimum and a unique diagonal maximum. By Proposition \ref{mirrorquiv}, $Q_{-K_{u/v}}$ also has a diagonal minimum.

    Here, $Q_{K_{u/v}}$ and $Q_{-K_{u/v}}$ are algebraic quivers, and since they are related to the geometric quivers by permutations of off-diagonal entries \cite{JKLN21}, they also satisfy the diagonal minimum conditions. Lastly, for both the algebraic and geometric quivers of a rational knot, we have $\text{diag}(Q) = S-2A$ up to a constant shift by Theorem 1.4 of \cite{Sto_i__2019} and by Theorem \ref{geothm}.
    Therefore $K$ and $-K$ are in class $\mathcal{Q}$.
\end{proof}

Now we prove Lemma \ref{maxineqs}. The proof is simply to check that the inequalities are preserved through all the 12 possible transformations. We will write the first 2 cases in detail, and the remaining cases will follow similarly. 

\begin{proof}[Proof of Lemma \ref{maxineqs}]
We first consider the base case $\tau_{0/1}$, which is of type $\Fa$-UP and holds trivially. We now proceed by induction.

Suppose $\tau$ has tangle quiver 
\[Q=
    \left(
    \begin{array}{c||c}
    Q_{++} & Q_{+-}\\ \hline \hline
    Q_{-+} & Q_{--}
    \end{array}
    \right).
\]

\textbf{Case 1:} $T^2(\Fa$-UP)

If $\tau$ is of type $\Fa$-UP, then $T^2\tau$ is of type $\Fa$-UP, and \[
T^2(Q) =
\left(
\begin{array}{cc||c}
Q_{++}+2 & Q_{+-} & Q_{+-}\\
Q_{-+} & Q_{--} & Q_{--}+L\\ \hline \hline
Q_{-+} & Q_{--}+U & Q_{--}
\end{array}
\right).
\]
If $\tau$ satisfies \eqref{faup1} with $\lambda$ the $Q_{++}$ unique diagonal maximum and $\eta$ the $Q_{--}$ unique diagonal maximum, then $T^2(Q)_{++}$ will have a unique diagonal maximum of $\lambda'=\lambda + 2$ coming from the $Q_{++}+2$ block. Then $Q_{+-} \le \lambda'+1$ and $Q_{--}+L \le \lambda' + 1$, so $T^2(Q)_{+-} \le \lambda'+1$. Moreover $T^2(Q)_{--}=Q_{--}$, so it has a unique diagonal maximum $\eta'=\eta$ with $\eta' \le \lambda'+1$.
Hence $T^2\tau$ satisfies \eqref{faup1}.

If $\tau$ satisfies \eqref{faup2} with $\lambda$ the $Q_{--}$ unique diagonal maximum, then $T^2(Q_{++})$ will have a unique diagonal maximum of $\lambda'=\lambda$ coming from the $Q_{--}$ block. Then $Q_{+-} \le \lambda'+1$ and $Q_{--}+L \le \lambda' + 1$, so $T^2(Q)_{+-} \le \lambda'+1$.  Moreover $T^2(Q)_{--}=Q_{--}$, so it has a unique diagonal maximum $\eta'=\lambda$ with $\eta' \le \lambda'+1$.
Hence $T^2\tau$ satisfies \eqref{faup1} as well.

\textbf{Case 2:} $RT(\Fa$-UP)

If $\tau$ is of type $\Fa$-UP, then $RT\tau$ is of type $\Fi$-OP, and
\[
RT(Q) =
\left(
\begin{array}{ccc||cc}
Q_{++}+1 & Q_{++}+1+U & Q_{+-} & Q_{++}+2+L & Q_{+-}+1\\
Q_{++}+1+L & Q_{++}+2 & Q_{+-} & Q_{++}+3+L & Q_{+-}+1\\
Q_{-+} & Q_{-+} & Q_{--} & Q_{-+}+2 & Q_{--}+1+L\\ \hline \hline
Q_{++}+2+U & Q_{++}+3+U & Q_{+-}+2 & Q_{++}+3 & Q_{+-}+1\\
Q_{-+}+1 & Q_{-+}+1 & Q_{--}+1+U & Q_{-+}+1 & Q_{--}+1
\end{array}
\right).
\]
If $\tau$ satisfies \eqref{faup1} with $\lambda$ the $Q_{++}$ unique diagonal maximum, then $RT(Q)_{--}$ will have a unique diagonal maximum of $\lambda'=\lambda + 3$ coming from the $Q_{++}+3$ block. Then $Q_{++} + 2 \le \lambda'$, $Q_{-+} \le \lambda'$, and $Q_{--} \le \lambda'$ so $RT(Q)_{++} \le \lambda'$. Moreover, $Q_{++} + 3 + U \le \lambda'$, $Q_{+-}+2 \le \lambda'$ and $Q_{--}+1+U\le \lambda'$, so $RT(Q)_{+-} \le \lambda'$.
Hence $RT\tau$ satisfies \eqref{fiopmax}.

If $\tau$ satisfies \eqref{faup2} with $\lambda$ the $Q_{--}$ unique diagonal maximum, then $RT(Q)_{--}$ will have a unique diagonal maximum of $\lambda'=\lambda+1$ coming from the $Q_{--}+1$ block. Then $Q_{++} + 2 \le \lambda'$, $Q_{-+} \le \lambda'$, and $Q_{--} \le \lambda'$ so $RT(Q)_{++} \le \lambda'$. Moreover, $Q_{++} + 3 + U \le \lambda'$, $Q_{+-}+2 \le \lambda'$ and $Q_{--}+1+U\le \lambda'$, so $RT(Q)_{+-} \le \lambda'$.
Hence $RT\tau$ satisfies \eqref{fiopmax} as well.

\textbf{Case 3:} $R(\Fa$-UP)

If $\tau$ is of type $\Fa$-UP, then $R\tau$ is of type $\Fa$-OP, and \[
R(Q)=
\left(
\begin{array}{c||ccc}
Q_{++} & Q_{++}+2+L & Q_{++}+1+L & Q_{+-}+1\\ \hline \hline
Q_{++}+2+U & Q_{++}+3 & Q_{++}+2+L & Q_{+-}+1\\
Q_{++}+1+U & Q_{++}+2+U & Q_{++}+2 & Q_{+-}+1\\
Q_{-+}+1 & Q_{-+}+1 & Q_{-+}+1 & Q_{--}+1
\end{array}
\right).
\]
If $\tau$ satisfies \eqref{faup1}, then the unique diagonal maximum of $R(Q)_{--}$ comes from the $Q_{++}+3$ block, and if $\tau$ satisfies \eqref{faup2}, then it comes from the $Q_{--}+1$ block. In either case the inequalities in \eqref{faopmax} hold.

\textbf{Case 4:} $T^2(\Fa$-OP)

If $\tau$ is of type $\Fa$-OP, then $T^2\tau$ is of type $\Fa$-OP, and \[
T^2(Q)= 
\left(
\begin{array}{cc||c}
Q_{++}-2 & Q_{+-}-3 & Q_{+-}-2\\ 
Q_{-+}-3 & Q_{--}-2 & Q_{--}-1+L\\ \hline \hline
Q_{-+}-2 & Q_{--}-1+U & Q_{--}
\end{array}
\right).
\]
and so the inequalities in \eqref{faopmax} are preserved. Note that if $Q_{--}$ has a unique diagonal maximum $\lambda$, then $Q_{--}-1+L<\lambda$.

\textbf{Case 5:} $RT(\Fa$-OP)

If $\tau$ is of type $\Fa$-OP, then $RT\tau$ is of type $\Fi$-RI, and \[
RT(Q)=
\left(
\begin{array}{ccc||cc}
Q_{++}-1 & Q_{++}-1+U & Q_{+-}-1 & Q_{++}-1+L & Q_{+-}-1\\
Q_{++}-1+L & Q_{++} & Q_{+-}-1 & Q_{++}+L & Q_{+-}-1\\
Q_{-+}-1 & Q_{-+}-1 & Q_{--} & Q_{-+} & Q_{--}+L\\ \hline \hline
Q_{++}-1+U & Q_{++}+U & Q_{+-} & Q_{++}-1 & Q_{+-}-2\\
Q_{-+}-1 & Q_{-+}-1 & Q_{--}+U & Q_{-+}-2 & Q_{--}-1
\end{array}
\right).
\]

Now $RT(Q)_{++}$ has a unique diagonal maximum in the $Q_{--}$ block, and $RT(Q)_{--}$ has a unique diagonal maximum in the $Q_{--}-1$ block, and the off-diagonal inequality in \eqref{firimax} holds.

\textbf{Case 6:} $R(\Fa$-OP)

If $\tau$ is of type $\Fa$-OP, then $R\tau$ is of type $\Fa$-UP, and \[
R(Q)=
\left(
\begin{array}{c||ccc}
Q_{++} & Q_{++}+1+L & Q_{++}+L & Q_{+-}+1\\ \hline \hline
Q_{++}+1+U & Q_{++}+1 & Q_{++}+L & Q_{+-}\\
Q_{++}+U & Q_{++}+U & Q_{++} & Q_{+-}\\
Q_{-+}+1 & Q_{-+} & Q_{-+} & Q_{--}+1
\end{array}
\right).
\]
$R(Q)_{--}$ has a unique diagonal maximum in the $Q_{--}+1$ block. Hence $R(Q)$ satisfies \eqref{faup2}.

\textbf{Case 7:} $R^2(\Fi$-OP)

If $\tau$ is of type $\Fi$-OP, then $R^2\tau$ is of type $\Fi$-OP, and \[
R^2(Q)=
\left(
\begin{array}{c||cc}
Q_{++} & Q_{++}+1+L & Q_{+-}+2\\ \hline \hline
Q_{++}+1+U & Q_{++}+1 & Q_{+-}+1\\
Q_{-+}+2 & Q_{-+}+1 & Q_{--}+2
\end{array}
\right).
\]
We have that $R^2(Q)_{--}$ has a unique diagonal maximum coming from the $Q_{--}+2$ block, and the inequalities in \eqref{fiopmax} are preserved.

\textbf{Case 8:} $TR(\Fi$-OP)

If $\tau$ is of type $\Fi$-OP, then $TR\tau$ is of type $\Fa$-UP, and \[
TR(Q)=
\left(
\begin{array}{cc||ccc}
Q_{++} & Q_{+-} & Q_{++}+L & Q_{+-} & Q_{+-}+1\\
Q_{-+} & Q_{--}+1 & Q_{-+}+1 & Q_{--}+1+L & Q_{--}+2+L\\ \hline \hline
Q_{++}+U & Q_{+-}+1 & Q_{++} & Q_{+-} & Q_{+-}\\
Q_{-+} & Q_{--}+1+U & Q_{-+} & Q_{--}+1 & Q_{--}+1+U\\
Q_{-+}+1 & Q_{--}+2+U & Q_{-+} & Q_{--}+1+L & Q_{--}+2
\end{array}
\right).
\]
We now have that $TR(Q)_{++}$ has a unique diagonal maximum in the $Q_{--}+1$-block, $TR(Q)_{--}$ has a unique diagonal maximum in the $Q_{--}+2$ block, and \eqref{faup1} is preserved.

\textbf{Case 9:} $T(\Fi$-OP)

If $\tau$ is of type $\Fi$-OP, then $T\tau$ is of type $\Fi$-RI, and \[
T(Q)=
\left(
\begin{array}{ccc||c}
Q_{++}-1 & Q_{+-}-1 & Q_{+-}-1 & Q_{+-}-1\\
Q_{-+}-1 & Q_{--} & Q_{--}+U & Q_{--}+L\\
Q_{-+}-1 & Q_{--}+L & Q_{--}+1 & Q_{--}+1+L\\ \hline \hline
Q_{-+}-1 & Q_{--}+U & Q_{--}+U+1 & Q_{--}
\end{array}
\right).
\]
We have that $T(Q)_{++}$ has a unique diagonal maximum coming from the $Q_{--}+1$ block, $T(Q)_{--}$ has a unique diagonal maximum coming from the $Q_{--}$ block, and the $T(Q)_{+-}$ inequality in \eqref{firimax} holds.

Lastly, we consider the case where $\tau$ is of type $\Fi$-RI.

\textbf{Case 10:} $R^2(\Fi$-RI)

If $\tau$ is of type $\Fi$-RI, then $R^2\tau$ is of type $\Fi$-RI, and \[
R^2(Q)=
\left(
\begin{array}{c||cc}
Q_{++} & Q_{++}+L & Q_{+-}\\ \hline \hline
Q_{++}+U & Q_{++}-1 & Q_{+-}-2\\
Q_{-+} & Q_{-+}-2 & Q_{--}-2
\end{array}
\right).
\]
$R^2(Q)_{++}$ has a unique diagonal maximum coming from the $Q_{++}$ block, $R^2(Q)_{--}$ has a unique diagonal maximum coming from the $Q_{++}-1$ block, and the $R^2(Q)_{+-}$ inequality in \eqref{firimax} holds.

\textbf{Case 11:} $TR(\Fi$-RI)

If $\tau$ is of type $\Fi$-RI, then $TR\tau$ is of type $\Fa$-OP, and \[
TR(Q)=
\left(
\begin{array}{cc||ccc}
Q_{++}-2 & Q_{+-}-3 & Q_{++}-1+L & Q_{+-}-2 & Q_{+-}-1\\
Q_{-+}-3 & Q_{--}-3 & Q_{-+}-1 & Q_{--}-2+L & Q_{--}-1+L\\ \hline \hline
Q_{++}-1+U & Q_{+-}-1 & Q_{++} & Q_{+-}-1 & Q_{+-}-1\\
Q_{-+}-2 & Q_{--}-2+U & Q_{-+}-1 & Q_{--}-1 & Q_{--}-1+U\\
Q_{-+}-1 & Q_{--}-1+U & Q_{-+}-1 & Q_{--}-1+L & Q_{--}
\end{array}
\right).
\]
$TR(Q)_{--}$ has a unique diagonal maximum coming from the $Q_{++}$ block, and the inequalities in \eqref{faopmax} hold.

\textbf{Case 12:} $T(\Fi$-RI)

Lastly, if $\tau$ is of type $\Fi$-RI, then $T\tau$ is of type $\Fi$-OP, and \[
T(Q)=
\left(
\begin{array}{ccc||c}
Q_{++}-1 & Q_{+-}-2 & Q_{+-}-2 & Q_{+-}-1\\
Q_{-+}-2 & Q_{--}-2 & Q_{--}-2+U & Q_{--}-1+L\\
Q_{-+}-2 & Q_{--}-2+L & Q_{--}-1 & Q_{--}+L\\ \hline \hline
Q_{-+}-1 & Q_{--}-1+U & Q_{--}+U & Q_{--}
\end{array}
\right).
\]
$T(Q)_{--}$ has a unique diagonal maximum, and the inequalities hold in \eqref{fiopmax}.

Therefore, the inequalities are preserved for all of the intermediate tangles.
\end{proof}

The proof of Lemma \ref{minineqs} is similar and is left to the reader.

\section{Examples}\label{sec:ex}

In this section, we compute $\Psi_K(a,q)$ for two families of knots as well as for various examples of small knots.
By the work of Osburn and Storzer \citep{osburn2025}, the Jones tail of a rational knot can be explicitly given as a product of theta/false theta functions $h_b$. Recall that these functions are defined as

\[
h_b(q^2)=\sum_{n\in\mathbb{Z}}\epsilon_b(n)q^{bn(n+1)-2n},
\]
where 
\[
\epsilon_b(n)=\begin{cases}
    (-1)^n &\text{if $b$ is odd},\\
    1 &\text{if $b$ is even and $n \ge 0$},\\
    -1 &\text{if $b$ is even and $n < 0$}.
\end{cases}
\]

In contrast, the form that the HOMFLY-PT tails of rational knots take is more complicated. That being said, there are clear patterns for some families of knots. Moreover, by looking at the torus and twist knots, we will finally prove at the end of Section \ref{TwistKSection} that the HOMFLY-PT tail is a stronger invariant than the Jones tail. 

For the following examples we will write $\Psi_{u/v}(a,q)$ to denote $\Psi_{K_{u/v}}(a,q)$. Additionally, when computing the head of a knot $K$ (the tail of the mirror knot $-K$), we will write the subquiver data as $S^h,A^h,Q^h$.

\subsection{Torus Knots}
\label{TorusKsection}

First, we compute the HOMFLY-PT tail for the torus knots $K_{(2n+1)/1}=T(2,2n+1)$. These have relatively simple quiver data, and give us a nice formula for the tail.

\begin{theorem}
    For all $n \ge 1$, $$\Psi_{{(2n+1)/1}}(a,q)=\frac{(a^{-2}q^{-2};q^2)_\infty}{(q^2;q^2)_\infty} \cdot h_{2n+1}(q^2).$$
\end{theorem}

\begin{proof}
    The quiver data of $K_{(2n+1)/1}$ is given in \citep{higgins2}, and it is easy to see that the subquiver data is:
 \[
S^t=\begin{pmatrix}
    -3\\2\\4\\6\\8\\\vdots
\end{pmatrix},\ 
A^t = \begin{pmatrix}
    -2\\0\\0\\0\\0\\\vdots
\end{pmatrix},\ 
Q^t = \begin{pmatrix}
    1&0&0&0&\dots&0\\
    0&2&2&2&\dots&2\\
    0&2&4&4&\dots & 4\\
    0&2&4&6&\dots&6\\
    \vdots&\vdots&\vdots&\vdots&\ddots&\vdots\\
    0&2&4&6&\dots&2n-2
\end{pmatrix} .
\]
Note that $Q^t$ is block-diagonal. Using Lemma \ref{prop:block}, the first block gives $\poch{a^{-2}q^{-2}}{q^2}{\infty}$, while the second gives \begin{align*}
&\sum_{v_1,\dots,v_{n-1}\ge 0}\frac{q^{\sum_{k=1}^{n-1}(\sum_{i=k}^{n-1} 2v_i+\sum_{i,j=k}^{n-1}2v_iv_j)}}{\prod\poch{q^2}{q^2}{v_i}}\\
=&\sum_{v_1, \dots,v_{n-1}\ge 0}\frac{q^{2\lt(\sum_{i=1}^{n-1}V_i^2+V_i\rt)}}{\prod_{i=1}^{n-1}(q^2;q^2)_{v_i}},\\
\intertext{where $V_i=v_i+v_{i+1}+\dots + v_{n-1}$,}
=& \frac{\lt(q^{2};q^{4n+2}\rt)_\infty\lt(q^{4n};q^{4n+2}\rt)_\infty\lt(q^{4n+2};q^{4n+2}\rt)_\infty}{\lt(q^2;q^2\rt)_\infty} & \trm{by \eqref{ag}}\\
= &\frac{1}{\poch{q^2}{q^2}{\infty}}h_{2n+1}(q^2) & \trm{by \eqref{jtp}}.
\end{align*}
Therefore, \[
\Psi_{{(2n+1)/1}}(a,q)=\frac{(a^{-2}q^{-2};q^2)_\infty}{(q^2;q^2)_\infty} \cdot h_{2n+1}(q^2).
\]
\end{proof}

The Jones head of the torus knots is $h_{2n+1}(q^2)$, so specializing to $\Psi_{K_{2n+1/1}}(q^{-2},q)$ gives the correct Jones head. This particular factor appears in many of the HOMFLY-PT tails, so we define the following $a$-deformed theta function.

\begin{definition}
    For $b\ge3$, define $$t_{b}(a^2,q^2):= \frac{(a^{-2}q^{-2};q^2)_\infty}{(q^2;q^2)_\infty} \cdot h_{b}(q^2).$$
\end{definition}

The HOMFLY-PT head of these torus knots also takes a simple form, being $\Psi_{-(2n+1)/1}(a,q)=1$ for all $n \ge 1$. Consider the geometric quiver $Q$ of $K_{(2n+1)/1}$, as calculated in \cite{higgins2}. One can see that the maximum on the diagonal of $Q$ is at least $2$ greater than every other element. Thus, the quiver of the mirror knot $Q'=I - 1 - Q$ will have a diagonal entry which is the unique minimum across the entire quiver. Therefore, the subquiver will be empty, giving us a head of $1$. This correctly aligns with the Jones tail, as they are also equal to $1$ for each of these torus knots.

\subsection{Twist Knots}
\label{TwistKSection}

In this section, we will prove a formula for the HOMFLT-PT head of the twist knots. The tails of the twist knots are all $t_3$ as the tail subquivers all take the same form as the subquiver of $K_{5/2}$. This can be seen from looking at the quivers of the twist knots, which are computed in Appendix \ref{TwistKAppendix}. Before proving a general form of the heads for the twist knots, we work out the heads for $K_{7/3}$ and $K_{9/4}$.

\begin{figure}
\centering
  \begin{tikzpicture}
      \node at (0,0) {\includegraphics[height=2.5cm]{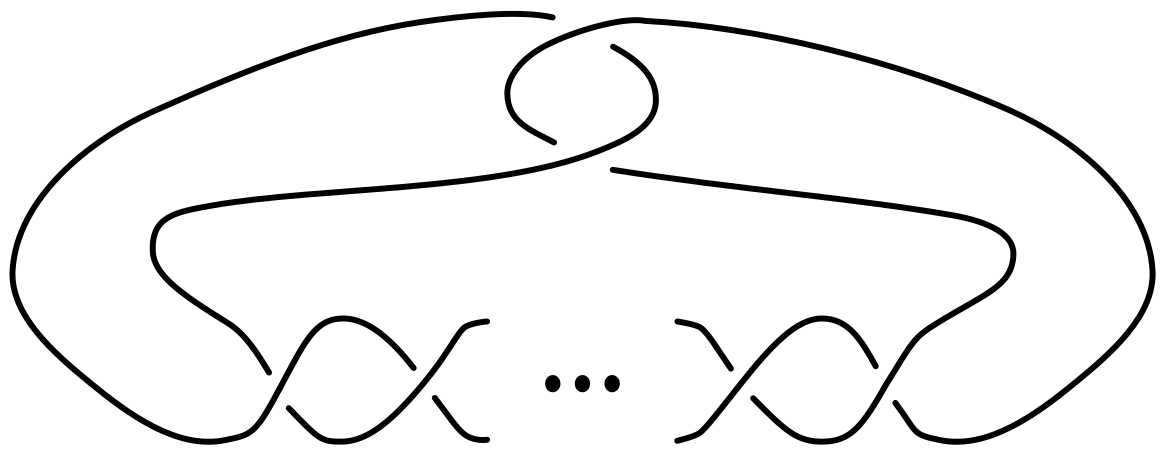}};
      \node at (0,-1.6) {$\underbrace{\hspace{3.7cm}}_{\textstyle n \,\text{half twists}}$};   
  \end{tikzpicture}
    \caption{The twist knot $K_{(2n+1)/n}$}
    \label{twistknotsfig}
\end{figure}

\begin{example}
    First, we compute the head of $K_{7/3}$. Recall that for the head, we want to look at the quiver $I - 1 - Q$ corresponding to the mirror knot.

We compute
\[
S= \begin{pmatrix}
    0\\-1\\-2\\-3\\0\\1\\2
\end{pmatrix},\
A = \begin{pmatrix}
    0\\-2\\-2\\-4\\-2\\0\\0
\end{pmatrix},\
Q=\begin{pmatrix}
    0&1&0&1&1&0&1\\1&3&2&3&3&2&2\\0&2&2&3&3&1&1\\1&3&3&5&4&2&2\\1&3&3&4&4&2&2\\0&2&1&2&2&1&2\\1&2&1&2&2&1&2
\end{pmatrix}.
\]

Thus, the subquiver data we get are $$S^h=\begin{pmatrix}-2\\1\end{pmatrix},\ A^h= \begin{pmatrix}-2\\0\end{pmatrix},\ Q^h = \begin{pmatrix}
    2&1\\1&1
\end{pmatrix}.$$

Now, we compute the head.

\begin{align*}
    \Psi_{-7/3}(a,q)&=\sum_{v_1,v_2\ge0}(-1)^{v_2}a^{-2v_1}q^{2v_1^2+v_2^2+2v_1v_2-2v_1+v_2}\frac{1}{(q^2;q^2)_{v_1}(q^2;q^2)_{v_2}}\\
    &=\sum_{v_1,v_2\ge0}q^{2v_1^2-2v_1}a^{-2v_1}\frac{(-1)^{v_2}q^{v_2^2-v_2}(a^2q^{4-2v_1})^{v_2}}{(q^2;q^2)_{v_1-v_2}(q^2;q^2)_{v_2}}&\text{shifting $v_1\to v_1-v_2$}\\
    &=\sum_{v_1\ge 0}q^{2v_1^2-2v_1}a^{-2v_1}\frac{(a^2q^{4-2v_1};q^2)_{v_1}}{(q^2;q^2)_{v_1}}&\text{by \eqref{eq:t/q} with $t=a^2q^{4-2v_1}$}\\
    &=\sum_{n\ge0}(-1)^nq^{n^2+n}\frac{(a^{-2}q^{-2};q^2)_n}{(q^2;q^2)_n}&\text{by \eqref{eq:qpochinverse}}.
\end{align*}

    Again, by specializing $a\to q^{-2}$, we recover the tail of the colored Jones polynomial for $K_{7/3}$. Thus $\Phi_{{7/3}} = \sum_{n\ge0}(-1)^nq^{n^2+n} = h_4(q^2)$.
    
\end{example}

\begin{example}\label{ex:9/4}

Computation for the head of $K_{9/4}$.

The subquiver data we have is 
$$S^h=\begin{pmatrix}-3\\-2\\-5\end{pmatrix},\ A^h= \begin{pmatrix}-2\\-2\\-4\end{pmatrix},\ Q^h = \begin{pmatrix}
    1&1&1\\1&2&2\\1&2&3
\end{pmatrix}.$$

Thus, we get:

\begin{align*}
    \Psi_{-9/4}(a,q)&=\sum_{v_1,v_2,v_3\ge0}(-1)^{v_1+v_3}a^{-2v_1-2v_2-4v_3} \frac{q^{v_1^2+2v_2^2+3v_3^2+2v_1v_2+4v_2v_3+2v_1v_3-3v_1-2v_2-5v_3}}{(q^2;q^2)_{v_1}(q^2;q^2)_{v_2}(q^2;q^2)_{v_3}}\\
    &=\sum_{v_2,v_3\ge0}(a^{-2}q^{2v_2+2v_3-2};q^2)_\infty (-1)^{v_3}a^{-2v_2-4v_3}\frac{q^{2v_2^2+3v_3^2+4v_2v_3-2v_2-5v_3}}{(q^2;q^2)_{v_2}(q^2;q^2)_{v_3}}&\text{solve } v_1 \text{ sum with } \eqref{eq:qbinom3}\\
    &=\sum_{v_2,v_3\ge0}(a^{-2}q^{2v_2-2};q^2)_\infty (-1)^{v_3}a^{-2v_2}q^{2v_2^2-2v_2} \frac{q^{v_3^2-v_3}(aq)^{-2v_3}}{(q^2;q^2)_{v_2-v_3}(q^2;q^2)_{v_3}}&\text{substitute } v_2 \to v_2 - v_3\\
    &=\sum_{v_2\ge0}(a^{-2}q^{2v_2-2};q^2)_\infty \frac{(a^{-2}q^{-2};q^2)_{v_2}}{(q^2;q^2)_{v_2}}a^{-2v_2}q^{2v_2^2-2v_2}&\text{solve } v_3 \text{ sum with \eqref{eq:t/q}}\\
    &=(a^{-2}q^{-2};q^2)_\infty \sum_{n\ge0}\frac{a^{-2n}q^{2n^2-2n}}{(q^2;q^2)_{n}}&\text{by \eqref{eq:qinfty}}.
\end{align*}

Lastly, specializing $a\to q^{-2}$, we get $\Phi_{{9/4}}(q) = (q^2;q^2)_\infty \sum_{n\ge0}\frac{q^{2n^2+2n}}{(q^2;q^2)_{n}}=h_5(q^2)$ by \eqref{eq:rr2} and \eqref{jtp}.
    
\end{example}

These previous 2 examples give us an idea about what the general form of the HOMFLY-PT tail for the twist knots looks like. We have derived a general formula whose statement we split into two cases, corresponding to the positive and negative twist knots. For $p>0$, the positive twist knot $K(p)$ is the rational knot $K_{(2n+1)/n}$ for $n=2p-1$, and the negative twist knot $K(-p)$ is the rational knot $K_{(2n+1)/n}$ for $n=2p$.

\begin{theorem}
    Let $p > 1 $ and consider the positive twist knots $K(p)$. Then 
$$
    \Psi_{-K(p)}(a,q)=\sum_{n_1,\dots, n_{p-1}\ge 0}(-1)^{N_1}a^{-2\sum_{i=2}^{p-1}N_i}q^{N_{1}^2+N_1+2\sum_{i=2}^{p-1}N_{i}(N_i-1)}\times\frac{\poch{a^{-2}q^{-2}}{q^2}{N_1}}{\prod_{i=1}^{p-1}\poch{q^2}{q^2}{n_i}}
$$
where $N_j=\sum_{i=j}^{p-1}n_i$.  
\end{theorem}

\begin{proof}
Before we calculate the head, it is easy to check
\begin{equation}\label{reindex_twistknot}
    \sum_{j=1}^{n}i_j^2-i_j-2i_j\sum_{k=j}^ni_k=-\lt(\sum_{j=1}^ni_j\rt)^2-\sum_{j=1}^ni_j.
\end{equation}

Let $p > 1$ and $K=-K(p)$. The subquiver data for $K$ is given by \[
S^h_i = \begin{cases}
    2-i & i\trm{ odd}\\
    -i & i\trm{ even}
\end{cases},\qquad
A^h_i = \begin{cases}
    1-i & i\trm{ odd}\\
    -i & i\trm{ even}
\end{cases},\qquad
Q^h_{i, j} = \trm{min}(i, j)
\]
for $1\le i, j\le 2|p|-2$.
Then we can write \begin{multline*}
\Psi_{K}=\sum_{v_1,\dots,v_{2p-2}\ge 0}(-1)^{\sum_{i=1}^{p-1}v_{2i-1}}a^{\sum_{i=1}^{p-1}(2-2i)v_{2i-1}-2iv_{2i}}q^{\sum_{i=1}^{p-1}(3-2i)v_{2i-1}-2iv_{2i}}\\q^{\sum_{i=1}^{2p-2}iv_i^2+\sum_{i=1}^{2p-2}\sum_{j=i+1}^{2p-2}2iv_iv_j}\frac{1}{\prod_{i=1}^{2p-2}\poch{q^2}{q^2}{v_i}}.
\end{multline*}
Now, shift $v_{2i}\to v_{2i}-v_{2i-1}$.
\begin{multline*}
    \Psi_{K}=\sum_{v_1,\dots v_{2p-2}\ge 0}(-1)^{\sum_{i=1}^{p-1}v_{2i-1}}a^{\sum_{i=1}^{p-1}2v_{2i-1}-2iv_{2i}}q^{\sum_{i=1}^{p-1}3v_{2i-1}-2iv_{2i}}\\
    q^{\sum_{i=1}^{p-1}2iv_{2i}^2-2v_{2i-1}v_{2i}+v_{2i-1}^2+\sum_{i=1}^{p-1}\sum_{j=i+1}^{p-1}4iv_{2i}v_{2j}-2v_{2i-1}v_{2j}}\inv{\prod_{i=1}^{p-1}\poch{q^2}{q^2}{v_{2i-1}}\prod_{i=1}^{p-1}\poch{q^2}{q^2}{v_{2i}-v_{2i-1}}}.
\end{multline*}
Next, we solve all $v_{2i-1}$ terms using \eqref{eq:t/q}: 
\begin{multline*}
    \Psi_{K} = \sum_{v_2,v_4,\dots, v_{2p-2}\ge 0}a^{\sum_{i=1}^{p-1}-2iv_{2i}}q^{\sum_{i=1}^{p-1}2iv_{2i}^2-2iv_{2i}+\sum_{i=1}^{p-1}\sum_{j=i+1}^{p-1}4iv_{2i}v_{2j}}\\
    \times \prod_{i=1}^{p-1}\frac{\poch{a^{2}q^{4-2\sum_{j=i}^{p-1}v_{2j}}}{q^2}{v_{2i}}}{\poch{q^2}{q^2}{v_{2i}}}.
\end{multline*}
Applying \eqref{eq:qpochinverse}, one obtains \begin{align*}
\prod_{i=1}^{p-1}\poch{a^{2}q^{4-2\sum_{j=i}^{p-1}v_{2j}}}{q^2}{v_{2i}}&=\prod_{i=1}^{p-1} (-1)^{v_{2i}} \lt(a^2q^4\rt)^{v_{2i}} q^{v_{2i}^2-v_{2_i}-2v_{2i}\sum_{j=i}^{p-1}v_{2j}} \poch{a^{-2}q^{-2+2\sum_{j=i+1}^{p-1}v_{2j}}} {q^2} {v_{2i}}\\
&=(-1)^{\sum_{i=1}^{p-1}v_{2i}}a^{2\sum_{i=1}^{p-1}v_{2i}}q^{-\lt(\sum_{i=1}^{p-1}v_{2i}\rt)^2+3\sum_{i=1}^{p-1}v_{2i}}\poch{a^{-2}q^{-2}}{q^2}{\sum_{i=1}^{p-1}v_{2i}}
\end{align*}
by \eqref{eq:tele} and \eqref{reindex_twistknot}. Combining it back, 
\begin{multline*}
\Psi_{K} = \sum_{v_2,v_4,\dots,v_{2p-2}\ge 0}(-1)^{\sum_{i=1}^{p-1}v_{2i}}a^{-2\sum_{i=2}^{p-1}\sum_{j=i}^{p-1}v_{2j}}q^{\lt(\sum_{i=1}^{p-1}v_{2i}\rt) \lt(\sum_{i=1}^{p-1}v_{2i}+1\rt)+2\sum_{i=2}^{p-1}\lt(\sum_{j=i}^{p-1}v_{2j}\rt)\lt(\sum_{j=i}^{p-1}v_{2j}-1\rt)}\\
\times \frac{\poch{a^{-2}q^{-2}}{q^2}{\sum_{i=1}^{p-1}v_{2i}}}{\prod_{i=1}^{p-1}\poch{q^2}{q^2}{v_{2i}}}.
\end{multline*}
Rewriting $v_{2i}=n_i$, $N_i=\sum_{j=i}^{p-1}n_j$, this looks like
\begin{equation}\label{pos_twist_tail}
    \Psi_{-K(p)}=\sum_{n_1,\dots, n_{p-1}\ge 0}(-1)^{N_1}a^{-2\sum_{i=2}^{p-1}N_i}q^{N_{1}^2+N_1+2\sum_{i=2}^{p-1}N_{i}(N_i-1)}\times\frac{\poch{a^{-2}q^{-2}}{q^2}{N_1}}{\prod_{i=1}^{p-1}\poch{q^2}{q^2}{n_i}}.
\end{equation}

\end{proof}

Next, we compute the colored HOMFLY-PT heads of negative twist knots.

\begin{theorem}
    Let $p > 1$ and consider the negative twist knots $K(-p)$. 
    Then $$
    \Psi_{-K(-p)}(a,q)=(a^{-2}q^{-2};q^2)_\infty\sum_{n_1,...,n_{p-1}\ge0}a^{-2\sum_{i=1}^{p-1}N_i}q^{-2\sum_{i=1}^{p-1}N_i+2\sum_{i=1}^{p-1}N_i^2}\frac{1}{\prod_{i=1}^{p-1}(q^2;q^2)_{n_i}}
    $$

    where $N_j=\sum_{i=j}^{p-1}n_i$.
\end{theorem}

\begin{proof}
    Let $p > 1$ and $K = -K(-p)$. The subquiver data for $K$ are given by 
    $$S^h_i = \begin{cases} 
      -i-2 & i \text{ odd} \\
      -i & i \text{ even}
   \end{cases},\qquad 
   A^h_i = \begin{cases}
       -i-1 & i \text{ odd} \\
       -i & i \text{ even} \\
   \end{cases},\qquad 
   Q^h_{ij}= \min(i,j)$$
   for $1 \le i,j \le 2p-1$.

   Thus, we can write 

\begin{multline*}
    \Psi_K(a,q)=\sum_{v_1,...,v_{2p-1}\ge0}(-1)^{\sum_{i=1}^{p}v_{2i-1}}a^{\sum_{i=1}^{p}(-2i)v_{2i-1}+\sum_{i=1}^{p-1}(-2i)v_{2i}} \\
    q^{\sum_{i=1}^{p}(-2i-1)v_{2i-1}+\sum_{i=1}^{p-1}(-2i)v_{2i}+\sum_{i=1}^{2p-1}iv_i^2+\sum_{i=1}^{2p-1}\sum_{j=i+1}^{2p-1}2iv_iv_j}\frac{1}{\prod_{i=1}^{2p-1}(q^2;q^2)_{v_i}}.
\end{multline*}

Note that from here, we proceed similarly as in \cite{Keilthy_2016}.

First, we solve the $v_1$ sum using \eqref{eq:qbinom3} to get

\begin{multline*}
    \Psi_K(a,q)=\sum_{v_2,...,v_{2p-1}\ge0} (a^{-2}q^{2(\sum_{i=2}^{2p-1}v_i-1)};q^{2})_\infty (-1)^{\sum_{i=1}^{p-1}v_{2i+1}}
    a^{\sum_{i=1}^{p-1} (-2i-2)v_{2i+1}+\sum_{i=1}^{p-1}(-2i)v_{2i}} \\
    q^{\sum_{i=1}^{p-1}(-2i-3)v_{2i+1}+\sum_{i=1}^{p-1}(-2i)v_{2i}+\sum_{i=2}^{2p-1}iv_i^2+\sum_{i=2}^{2p-1}\sum_{j=i+1}^{2p-1}2iv_iv_j}\frac{1}{\prod_{i=2}^{2p-1}(q^2;q^2)_{v_i}}.
\end{multline*}

Now, we shift $v_{2i}\to v_{2i} - v_{2i+1}$ for all $1\le i \le p-1 $, which gives

\begin{multline*}
    \Psi_K(a,q)=\sum_{v_2,...,v_{2p-1}\ge0} (a^{-2}q^{2(\sum_{i=1}^{p-1}v_{2i}-1)};q^{2})_\infty (-1)^{\sum_{i=1}^{p-1}v_{2i+1}}
    a^{\sum_{i=1}^{p-1}(-2)v_{2i+1}+\sum_{i=1}^{p-1}(-2i)v_{2i}} \\
    q^{\sum_{i=1}^{p-1}(-3)v_{2i+1}-
    \sum_{i=1}^{p-1}2iv_{2i}+
    \sum_{i=1}^{p-1}2iv_{2i}^2+
    \sum_{i=1}^{p-1}v_{2i+1}^2+    
    \sum_{i=1}^{p-1}\sum_{j=i+1}^{p-1}4iv_{2i}v_{2j}+  \sum_{i=1}^{p-1}\sum_{j=i+1}^{p-1}2v_{2i+1}v_{2j}}\\
    \frac{1}{\prod_{i=1}^{p-1}(q^2;q^2)_{v_{2i}-v_{2i+1}}\prod_{i=1}^{p-1}(q^2;q^2)_{v_{2i+1}}}.
\end{multline*}

Next, we solve all $v_{2i+1}$ sums for $1 \le i \le p-1$ using \eqref{eq:t/q}:

\begin{multline*}
    \Psi_K(a,q)=\sum_{v_2,v_4,...,v_{2p-2}\ge0} (a^{-2}q^{2(\sum_{i=1}^{p-1}v_{2i}-1)};q^{2})_\infty 
    \prod_{i=1}^{p-1}(a^{-2}q^{2(\sum_{j=i+1}^{p-1} v_{2j}-1)};q^2)_{v_{2i}}
    a^{-2\sum_{i=1}^{p-1}iv_{2i}}\\
    q^{-2\sum_{i=1}^{p-1}iv_{2i}+
    2\sum_{i=1}^{p-1}iv_{2i}^2+    \sum_{i=1}^{p-1}\sum_{j=i+1}^{p-1}4iv_{2i}v_{2j}}
    \frac{1}{\prod_{i=1}^{p-1}(q^2;q^2)_{v_{2i}}}.
\end{multline*}

By repeatedly applying \eqref{eq:qinfty} to the Pochhammer symbols, we have  $$(a^{-2}q^{-2};q^2)_\infty=(a^{-2}q^{2(\sum_{i=1}^{p-1}v_{2i}-1)};q^{2})_\infty \prod_{i=1}^{p-1}(a^{-2}q^{2(\sum_{j=i+1}^{p-1} v_{2j}-1)};q^2)_{v_{2i}}.$$

Hence, we have

\[
    \Psi_K(a,q)=(a^{-2}q^{-2};q^{2})_\infty \sum_{v_2,v_4,\dots,v_{2p-2}\ge0}    
    \frac{a^{-2\sum_{i=1}^{p-1}iv_{2i}}
    q^{-2\sum_{i=1}^{p-1}iv_{2i}+
    2\sum_{i=1}^{p-1}iv_{2i}^2+    \sum_{i=1}^{p-1}\sum_{j=i+1}^{p-1}4iv_{2i}v_{2j}}}{\prod_{i=1}^{p-1}(q^2;q^2)_{v_{2i}}}.
\]

Lastly, we reindex $n_i=v_{2i}$ for $1 \le i \le p-1$, and set $N_j = \sum_{i=j}^{p-1}n_i$ to get
\[
    \Psi_{-K(-p)}(a,q)=(a^{-2}q^{-2};q^2)_\infty\sum_{n_1,...,n_{p-1}\ge0}    \frac{a^{-2\sum_{i=1}^{p-1}N_i}q^{-2\sum_{i=1}^{p-1}N_i+2\sum_{i=1}^{p-1}N_i^2}}{\prod_{i=1}^{p-1}(q^2;q^2)_{n_i}}.
\]

\end{proof}

One may check that specializing $a\to q^{-2}$ gives us the odd $h_b$ functions for the negative twist knots and the even $h_b$ functions for the positive twist knots, which agrees with the Jones tail. Like with the torus knot case, we will define another $a$-deformed function for these families of knots.

\begin{definition}
    For $b \ge 3$, define 
    $$e_b(a^2,q^2):= 
    \begin{cases}
        \Psi_{-K(b/2)}(a,q) &\text{$b$ even},\\
        \Psi_{-K(-(b-1)/2)}(a,q) &\text{$b$ odd}.
    \end{cases}$$
\end{definition}

Note that $e_3=t_3$ from these definitions, so we will by convention write $t_3$ for $(a^{-2}q^{-2};q^2)_\infty$.

\begin{proposition}\label{diff_tails}
    For $p > 1$, $e_{2p+1} \neq t_{2p+1}$.
\end{proposition}

\begin{proof}
    Let $p > 1$. Then 
    \begin{align*}
        e_{2p+1} &= (a^{-2}q^{-2};q^2)_\infty\sum_{n_1,...,n_{p-1}\ge0}    \frac{a^{-2\sum_{i=1}^{p-1}N_i}q^{-2\sum_{i=1}^{p-1}N_i+2\sum_{i=1}^{p-1}N_i^2}}{\prod_{i=1}^{p-1}(q^2;q^2)_{n_i}}\\
        t_{2p+1}&= \frac{(a^{-2}q^{-2};q^2)_\infty}{(q^2;q^2)_\infty}h_{2p+1}(q^2).
    \end{align*}

    Suppose that they are equal. Then, after canceling the $(a^{-2}q^{-2};q^2)_\infty$ factor, we get

    $$\sum_{n_1,...,n_{p-1}\ge0}    \frac{a^{-2\sum_{i=1}^{p-1}N_i}q^{-2\sum_{i=1}^{p-1}N_i+2\sum_{i=1}^{p-1}N_i^2}}{\prod_{i=1}^{p-1}(q^2;q^2)_{n_i}} = \frac{1}{(q^2;q^2)_\infty}h_{2p+1}(q^2).$$ 
    Looking at the left hand side, since $N_i \ge N_{i+1}$ for all $i$, a term with $a^{-2}$ can only appear when $n=(1,0,...,0)$ in the summation. Then the coefficient of $a^{-2}$ is $\frac{1}{1-q^2}$, but the right hand side is independent of $a$. Therefore, they cannot be equal. 
\end{proof}

\begin{remark}
    Recall that given the HOMFLY-PT tail $\Psi_K(a,q)$ of a knot $K$, the Jones head of $K$ is $\Phi_{-K}(q)=\Psi_{K}(q^{-2},q)$. The HOMFLY-PT tail specializing to the Jones head is due to the colored HOMFLY-PT polynomials being $\bigwedge^j$-colored invariants and the colored Jones polynomials being $Sym^j$-colored invariants. Hence, to discuss the strength of these invariants, we compare the HOMFLY-PT tails to the Jones heads.
\end{remark}

Combining Proposition \ref{diff_tails} with the computations for torus and twist knots gives the following result.

\begin{corollary}
    For all $n > 1$, consider the torus knot $K_{(2n+1)/1}$ and the mirror negative twist knot $K_{-(4n+1)/2n}$. Then we have $\Psi_{(2n+1)/1}(q^{-2},q)=\Psi_{-(4n+1)/2n}(q^{-2},q)$, but $\Psi_{(2n+1)/1}(a,q)\neq\Psi_{-(4n+1)/2n}(a,q)$. That is, they have the same Jones head yet different HOMFLY-PT tails.
\end{corollary}

We now have computed the HOMFLY-PT tails of two families of knots which give the same Jones heads yet different HOMFLY-PT tails. Thus, the tail of the colored HOMFLY-PT polynomial is a stronger invariant than the tail of the colored Jones polynomial. The Jones tail has a correspondence with reduced Tait graphs \cite{garoufalidis2012nahmsumsstabilitycolored}, so one may wonder what underlying structure of the knot the HOMFLY-PT tail is related to.

\subsection{Other Examples}

From the previous families of knots we have computed, it seems like it is very complicated to compute these tails. However, in some cases, the subquivers break up into block diagonal matrices consisting of blocks that we already know the tails of. We will also give the Alexander-Briggs notation for the following examples.

\begin{example}
    For an example of a knot that is not a torus or twist knot, we compute the tail of $K_{29/6}$. In Alexander-Briggs notation, this is the knot $10_8$.

    For the subquiver data, we have \[
    S^t = 
    \begin{pmatrix}
         -3\\-2\\-5\\1\\2\\3\\4
    \end{pmatrix},\ 
    A^t = \begin{pmatrix}
        -2\\-2\\-4\\0\\0\\0\\0
    \end{pmatrix},\ 
    Q^t = \begin{pmatrix}
        1&1&1&0&0&0&0\\
        1&2&2&0&0&0&0\\
        1&2&3&0&0&0&0\\
        0&0&0&1&1&1&1\\
        0&0&0&1&2&2&2\\
        0&0&0&1&2&3&3\\
        0&0&0&1&2&3&4\\
    \end{pmatrix}.
    \]
    By the block-diagonal format of $Q^t$, we reduce to computations we have already done. These are simply the subquivers for the twist knots $K(-2)$ and $K(3)$. Note that the 2nd block has 0 in every A entry, which corresponds to a specialization of $a \to q^{-2}$, meaning this gives the Jones head of $K(3)$. 
    Thus by Proposition \ref{prop:block}, \[
    \Psi_{-29/6}(a,q)=\lt((a^{-2}q^{-2};q^2)_\infty \sum_{n\ge0}\frac{a^{-2n}q^{2n^2-2n}}{(q^2;q^2)_{n}}\rt)h_6(q^2)=e_5h_6
    \]
    and
    \[
    \Phi_{29/6}(q)=h_5h_6.
    \]
\end{example}

\begin{example}
    Computation of the head of $K_{13/8}=6_3$.
        For the subquiver data, we have \[
    S^h = 
    \begin{pmatrix}
        1\\-3\\-2
    \end{pmatrix},\ 
    A^h = \begin{pmatrix}
        0\\-2\\-2
    \end{pmatrix},\ 
    Q^h = \begin{pmatrix}
        1&1&1\\1&1&1\\1&1&2
    \end{pmatrix}.
    \]

    \begin{align*}
    \Psi_{-13/8}(a,q)&=\sum_{v_1,v_2,v_3\ge0}(-1)^{v_1+v_2}a^{-2v_2-2v_3} \frac{q^{v_1^2+v_2^2+2v_3^2+v_1-3v_2-2v_3+2v_1v_2+2v_1v_3+2v_2v_3}}{(q^2;q^2)_{v_1}(q^2;q^2)_{v_2}(q^2;q^2)_{v_3}}\\
    &=\sum_{v_2,v_3\ge0}(q^{2v_2+2v_3+2};q^2)_\infty (-1)^{v_2}a^{-2v_2-2v_3}\frac{q^{v_2^2+2v_3^2+2v_2v_3-3v_2-2v_3}}{(q^2;q^2)_{v_2}(q^2;q^2)_{v_3}}&\text{solve } v_1 \text{ sum with } \eqref{eq:qbinom3}\\
    &=\sum_{v_2,v_3\ge0}(q^{2v_2+2};q^2)_\infty (-1)^{v_2+v_3}a^{-2v_2}q^{v_2^2-3v_2} \frac{q^{v_3^2-v_3}q^{2v_3}}{(q^2;q^2)_{v_2-v_3}(q^2;q^2)_{v_3}}&\text{substitute } v_2 \to v_2 - v_3\\
    &=\sum_{v_2\ge0}(q^{2v_2+2};q^2)_\infty \frac{(q^2;q^2)_{v_2}}{(q^2;q^2)_{v_2}}(-1)^{v_2}a^{-2v_2}q^{v_2^2-3v_2}&\text{solve } v_3 \text{ sum with \eqref{eq:t/q}}\\
    &=(a^{-2}q^{-2};q^2)_\infty h_3(q^2) &\text{by \eqref{eq:qinfty} and \eqref{eq:qbinom3}}\\
    &=t_3h_3 .
\end{align*}
    
\end{example}

\begin{example}
    Computation of the head of $K_{27/10}=8_{11}$.

    For the subquiver data, we have \[
    S^h = 
    \begin{pmatrix}
        -3\\1\\-2\\-2\\-5
    \end{pmatrix},\ 
    A^h = \begin{pmatrix}
        -2\\0\\-2\\-2\\-4
    \end{pmatrix},\ 
    Q^h = \begin{pmatrix}
        1&1&1&1&1\\1&1&1&1&2\\1&1&2&1&2\\1&1&1&2&2\\
        1&2&2&2&3
    \end{pmatrix}.
    \]

    Direct computation shows that $\Psi_{-27/10}(a,q)=(a^{-2}q^{-2};q^2)_\infty\sum_{n\ge0}(-1)^nq^{n^2+n}\frac{(a^{-2}q^{-2};q^2)_n}{(q^2;q^2)_n}=t_3e_4$.
\end{example}

\begin{example}\label{ex:29/12}
    Computation of the tail of $K_{29/12}=8_{12}$.

    For the subquiver data, we have \[
    S^t = 
    \begin{pmatrix}
        -3\\-3\\-2\\-6\\-5
    \end{pmatrix},\ 
    A^t = \begin{pmatrix}
        -2\\-2\\-2\\-4\\-4
    \end{pmatrix},\ 
    Q^t = \begin{pmatrix}
        1&1&1&1&1\\1&1&1&1&1\\1&1&2&2&2\\1&1&2&2&2\\
        1&1&2&2&3
    \end{pmatrix}.
    \]

    Direct computation shows that $\Psi_{29/12}(a,q)=(a^{-2}q^{-2};q^2)_\infty\sum_{n\ge0}(-1)^n(aq)^{-2n}q^{n^2-n}$. The sum in this tail has not appeared in any previous tails, so we name it $o_4(a^2,q^2)=\sum_{n\ge0}(-1)^n(aq)^{-2n}q^{n^2-n}$, as it specializes to $h_4$. Then $\Psi_{29/12}(a,q)=t_3o_4$.
\end{example}

\begin{remark}
    It turns out that the previous two subquivers are related by a perturbation matrix which does not influence the calculation of the HOMFLY-PT tail. Thus, the reason for the differing tails is due to the $A$ vectors and not the subquivers. 
\end{remark}

\begin{conjecture}\label{conj:weird}
    Let $k> 1$ be such that $k \equiv 1 \pmod{4}$. Then, consider a knot $K$ with the recursively defined subquiver data 
    $$S^t_{1,2,3,4}=(-3, -3, -2, -6), \quad S^t_{n+4}=S^t_{n} - 2$$
    $$A^t_{1,2,3,4}=(-2, -2, -2, -4), \quad A^t_{n+4}=A^t_{n} -2$$
    $$Q^t_{i,j} = \min\left(\left\lceil\frac{i}{2} \right\rceil, \left\lceil\frac{j}{2} \right\rceil\right)$$
    for $1 \le n \le k-4$ and $1 \le i,j \le k$.
    
    Then, $\Phi_{-K}(q)=h_3h_{(k+3)/2}$.

    Furthermore, we can write $\Psi_K(a,q)=t_3o_{(k+3)/2}$ for some function $o_b(a^2,q^2)$ that specializes to $h_b(q^2)$.
\end{conjecture}

Consider $k$ and $K$ as above. Taking this conjecture to be true, without even explicitly calculating the HOMFLY-PT tail of $K$, we can give a definition of $o_b$ in terms of $\Psi_K(a,q)$. That is, we can write $$o_b(a^2,q^2) := \frac{\Psi_K(a,q)}{t_3}$$ for $b = (k+3)/2$. This definition of $o_b$ is motivated by Example \ref{ex:29/12}, which has exactly this subquiver data.

Computational data suggests that this subquiver data appears for the knots $K_{(7+12n)/(3+4n)}$ for $n > 1$. We calculate the continued fraction to be
\[
\frac{7+12n}{3+4n}=2+\frac{1}{1+\frac{1}{2n+\frac{1}{1+\frac{1}{1}}}},
\]
from which one can calculate their Jones tails to be $h_3h_{2n+2}$ via the methods of \cite{osburn2025}. In addition, we have that their HOMFLY-PT tails are $\Psi_K(a,q)=t_3o_{2n+2}$ by Conjecture \ref{conj:weird}. Hence, to prove this conjecture, one could try to show that this family of rational knots has the above subquiver data.

\begin{example}
    We can also work backwards to find out which quivers represent which Jones or HOMFLY-PT tails. Consider $K_{31/7}=9_9$. Its subquiver data is:

    \[
    S^t = 
    \begin{pmatrix}
        -3\\1\\2\\2\\3\\4
    \end{pmatrix},\ 
    A^t = \begin{pmatrix}
        -2\\0\\0\\0\\0\\0
    \end{pmatrix},\ 
    Q^t = \begin{pmatrix}
        1&0&0&0&0&0\\
        0&1&1&1&1&2\\
        0&1&2&1&1&2\\
        0&1&1&2&2&2\\
        0&1&1&2&3&3\\
        0&2&2&2&3&4
    \end{pmatrix}.
    \]
    From the continued fraction decomposition of $31/7$, we can compute its Jones tail to be $\Phi_{-31/7}(q)=h_4h_5$. Reading off the blocks of the subquiver data, we see the 2nd block has a $0$ A vector, so there is no $a$ contribution from that block. Then the first block contributes $h_3$ in the Jones case, meaning the second block must contribute $h_4h_5/h_3$. In the HOMFLY-PT case, the first block is $t_3$, so we can compute the HOMFLY-PT tail to be $\Psi_{31/7}(a,q)=t_3h_4h_5/h_3=h_4t_5$.
\end{example}

\subsection{Table of HOMFLY-PT tails}

So far we have seen 4 different families of subquiver data. It turns that up to numerator 30, all rational knots are built out of blocks of these subquivers. Thus we reduce the computations to recognizing the blocks that these subquivers are built out of. Recall the $t_b$, $e_b$, and $o_b$ functions we defined earlier:
\begin{align*}
    t_b(a^2,q^2) &:= \frac{(a^{-2}q^{-2};q^2)_\infty}{(q^2;q^2)_\infty}h_b(q^2) & \text{for $b \ge 3$}\\
    e_b(a^2,q^2) &:= 
\begin{cases}
    \Psi_{-K(b/2)}(a,q) &\text{$b$ even}\\
    \Psi_{-K(-(b-1)/2)}(a,q) &\text{$b$ odd}
\end{cases} & \text{for $b \ge 3$}\\
    o_{b}(a^2,q^2) &:= \frac{\Psi_{(6b-5)/(2b-1)}(a,q)}{t_3} &\text{for $b\ge 4$ even (conjectured)}.
\end{align*}

These functions are defined such that $t_b(q^{-4},q^2)=e_b(q^{-4},q^2)=o_b(q^{-4},q^2)=h_b(q^2)$. The HOMFLY-PT tails are given in Table \ref{tab1}. The subquiver data for each knot was computed using \cite{homfly-tails}, from which the tails were determined by matching the subquivers to previous examples. 

\begin{table}[!ht]
\begin{tabular}{|c|c|c||c|c|c||c|c|c|}
\hline
$K$ & $\Psi_{K}(a,q)$ & $\Psi_{-K}(a,q)$ & $K$ & $\Psi_{K}(a,q)$ & $\Psi_{-K}(a,q)$ & $K$ & $\Psi_{K}(a,q)$ & $\Psi_{-K}(a,q)$\\ \hline
$3/1$ & $ t_{3} $ & $ 1 $ & $17/4$ & $ e_{5} $ & $ e_{5} $ & $25/2$ & $ e_{13} $ & $ t_{3} $\\
$5/1$ & $ t_{5} $ & $ 1 $ & $17/5$ & $ t_{3}  h_{4} $ & $ e_{4} $ & $25/3$ & $ t_{9} $ & $ e_{4} $\\
$5/2$ & $ t_{3} $ & $ t_{3} $ & $17/6$ & $ t_{3}  h_{6} $ & $ t_{3} $ & $25/4$ & $ e_{7} $ & $ e_{5} $\\
$7/1$ & $ t_{7} $ & $ 1 $ & $19/1$ & $ t_{19} $ & $ 1 $ & $25/14$ & $ h_{3}  t_{5} $ & $ h_{3}  t_{3} $\\
$7/3$ & $ t_{3} $ & $ e_{4} $ & $19/3$ & $ t_{7} $ & $ e_{4} $ & $25/18$ & $ t_{3}  h_{4} $ & $ t_{3}  h_{4} $\\
$9/1$ & $ t_{9} $ & $ 1 $ & $19/4$ & $ h_{4}  e_{5} $ & $ t_{3} $ & $27/1$ & $ t_{27} $ & $ 1 $\\
$9/2$ & $ e_{5} $ & $ t_{3} $ & $19/7$ & $ t_{3}  o_{4} $ & $ h_{3}  t_{3} $ & $27/4$ & $ h_{4}  e_{7} $ & $ t_{3} $\\
$11/1$ & $ t_{11} $ & $ 1 $ & $19/9$ & $ t_{3} $ & $ e_{10} $ & $27/5$ & $ t_{3}  h_{6} $ & $ e_{4} $\\
$11/4$ & $ t_{3}  h_{4} $ & $ t_{3} $ & $21/1$ & $ t_{21} $ & $ 1 $ & $27/8$ & $ t_{3}  e_{4} $ & $ t_{3}  h_{4} $\\
$11/5$ & $ t_{3} $ & $ e_{6} $ & $21/2$ & $ e_{11} $ & $ t_{3} $ & $27/13$ & $ t_{3} $ & $ e_{14} $\\
$13/1$ & $ t_{13} $ & $ 1 $ & $21/5$ & $ t_{5} $ & $ e_{6} $ & $29/1$ & $ t_{29} $ & $ 1 $\\
$13/2$ & $ e_{7} $ & $ t_{3} $ & $21/13$ & $ t_{3}^2 $ & $ t_{3}  h_{3}^2 $ & $29/2$ & $ e_{15} $ & $ t_{3} $\\
$13/3$ & $ t_{5} $ & $ e_{4} $ & $23/1$ & $ t_{23} $ & $ 1 $ & $29/6$ & $ e_{5}  h_{6} $ & $ t_{3} $\\
$13/8$ & $ h_{3}  t_{3} $ & $ h_{3}  t_{3} $ & $23/8$ & $ t_{3}  h_{8} $ & $ t_{3} $ & $29/7$ & $ t_{5} $ & $ e_{8} $\\
$15/1$ & $ t_{15} $ & $ 1 $ & $23/10$ & $ t_{3}  h_{4} $ & $ e_{5} $ & $29/8$ & $ e_{4}  h_{3}  t_{3} $ & $ h_{3}  t_{3} $\\
$15/7$ & $ t_{3} $ & $ e_{8} $ & $23/11$ & $ t_{3} $ & $ e_{12} $ & $29/9$ & $ t_{3}  h_{4} $ & $ e_{6} $\\
$15/11$ & $ t_{3} $ & $ e_{4}^2 $ & $23/14$ & $ t_{5}  h_{3} $ & $ h_{3}  t_{3} $ & $29/10$ & $ t_{3}  h_{10} $ & $ t_{3} $\\
$17/1$ & $ t_{17} $ & $ 1 $ & $23/17$ & $ t_{3} $ & $ e_{4}  e_{6} $ & $29/12$ & $ t_{3}  o_{4} $ & $ t_{3}  o_{4} $\\
$17/2$ & $ e_{9} $ & $ t_{3} $ & $25/1$ & $ t_{25} $ & $ 1 $ & $31/1$ & $ t_{31}$ & $1$\\

\hline
\end{tabular}
\caption{Colored HOMFLY-PT heads and tails of rational knots with numerator up to 30}
\label{tab1}
\end{table}

\section{Questions and Conjectures}\label{sec:questions}

\subsection{Shape of the Colored HOMFLY-PT polynomials}

In Section \ref{subsec:shape}, we showed that there is a quadratic bound on the tail of the colored HOMFLY-PT polynomial. We believe that this may extend in some sense to the entire colored HOMFLY-PT polynomial. An interesting phenomenon appearing in our visualizations of the colored HOMFLY-PT
polynomials is that, for a fixed rational knot $K_{u/v}$, the set of $(a,q)$-monomials appearing in
$P^j_{K_{u/v}}(a,q)$ seems to approach a limiting shape as we increase the color $j$.

\begin{figure}
    \centering
    \includegraphics[width=1\linewidth]{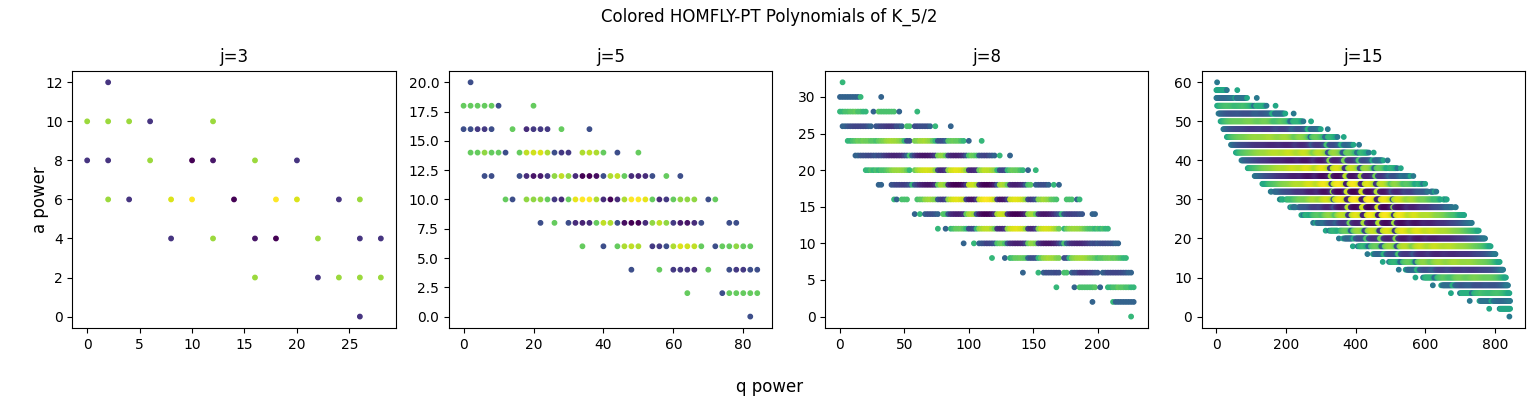}
    \hfill
    \includegraphics[width=1\linewidth]{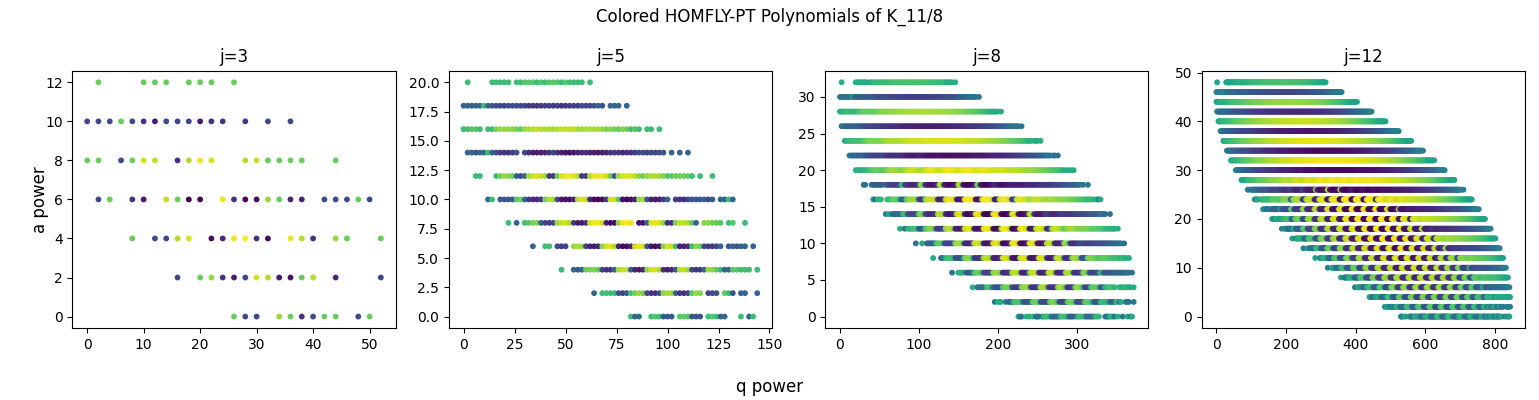}
    \hfill
    \includegraphics[width=1\linewidth]{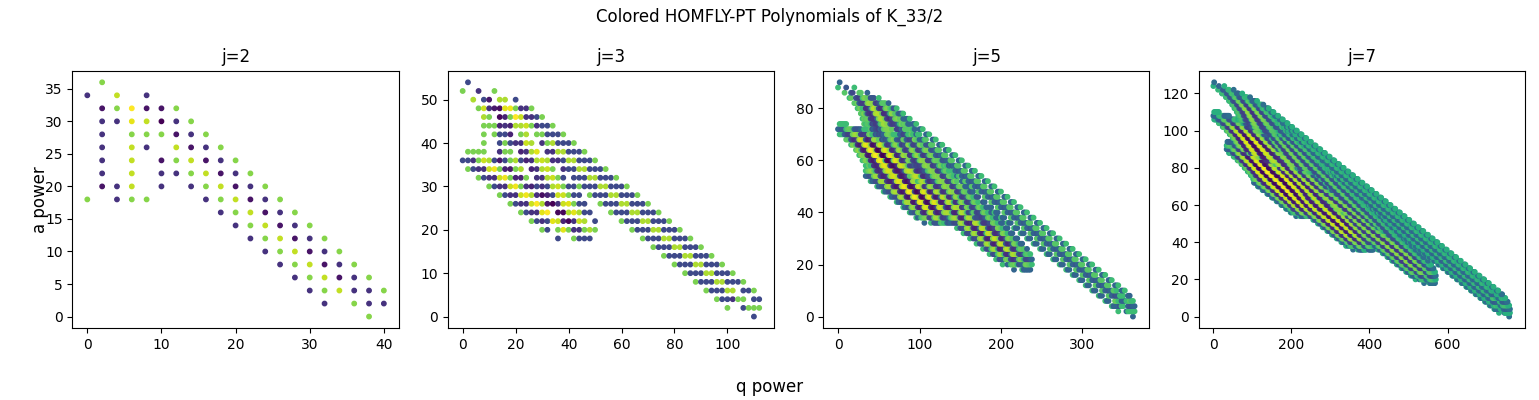}
    \caption{$j$th colored HOMFLY-PT polynomials of $K_{5/2}$, $K_{11/8}$, and $K_{33/2}$}
    \label{fig:homflyplots}
\end{figure}

Recall that we write
\[
{\widehat{P}^j_{K_{u/v}}}(a,q)
=
\sum_{s\in\mathbb{Z},\,r\geq 0} c^j_{r,s} a^{2s}q^{2r}.
\]
Since the $a$-degree is determined by $A\cdot d$ and the $q$-degree is determined by $d^T\cdot Q\cdot d+S\cdot d$, where $|d|=j$,
it is natural to rescale the $a$-degree linearly in $j$ and the $q$-degree quadratically in $j$. The plots in Figure \ref{fig:homflyplots}
suggest that for certain knots, the support for the plots of $a$- and $q$-powers approaches a fixed region.

\begin{conjecture}[Limiting Shape of the colored HOMFLY-PT polynomial]
\label{limshapeconj}
    Let $K_{u/v}$ be a rational knot and write $\hat{P}^j_{K_{u/v}}(a,q)=\sum_{s\in\mathbb{Z},r\geq0}c^j_{s,r}a^{2s}q^{2r}$. For each $j\geq1$, we define the rescaled support $$\mathcal{S}_j(K_{u/v})=\left\{ \left(\frac{s}{j},\frac{r}{j^2}\right)\;\middle|\;c^j_{s,r}\neq0\ \right\}.$$ Then there exists a compact region $$\mathcal{S}(K_{u/v})\subset\mathbb{R}^2$$ so that $$\lim_{j\longrightarrow{\infty}}d_H\left(\mathcal{S}_j(K_{u/v}),\mathcal{S}(K_{u/v})\right)=0$$ where $d_H$ denotes the Hausdorff metric. We call $\mathcal{S}(K_{u/v})$ the limiting shape of the colored HOMFLY-PT polynomial of $K_{u/v}$.
\end{conjecture}

For several rational knots, this region appears to be convex. In particular, the examples
$K_{5/2}$ and $K_{11/8}$ exhibit increasingly well-defined shapes as the color increases. The rational
torus knots $K_{(2n+1)/1}$ give a particularly simple example of the phenomenon,
coinciding with their relatively simple quiver data. In contrast, examples such as $K_{33/2}$
suggest that its limiting shape may be composed of several boundary pieces for other knots. For more examples, see \cite{homfly-tails}.

\subsection{Rate of Convergence of HOMFLY-PT Tail}

When the existence of the colored Jones tail for alternating knots was proved, the authors also showed that the $n$th coefficient was independent of the color $j$ for $j \ge n$ \cite{garoufalidis2012nahmsumsstabilitycolored}. One may also wonder about the convergence for the HOMFLY-PT case.

\begin{conjecture}
Let $K_{u/v}$ be a rational knot and consider $\hat{P}_{K_{u/v}}^j(a,q) = \sum_{s\in\mathbb{Z},r\ge0}c^j_{s,r}a^{2s}q^{2r}$. Fix $s \in \mathbb{Z}$. Then for all $n \in\mathbb{N}$, $c^j_{s,n}$ is independent of $j$ for all $j \ge n$.
\end{conjecture}

Note that it may be that there exists $k < n$ such that $c^j_{s,n}$ is independent of $j$ for all $j \ge k$. Additionally, it is unknown whether or not one can determine the first nonzero coefficient for each $c^j_{s,r}$ when fixing $s$. This conjecture has been verified for the low coefficients of rational knots up to numerator $15$. Visualizations of this phenomenon can be found in \cite{homfly-tails}. 

\subsection{Nonvanishing of Colored \texorpdfstring{$\mathfrak{sl}_N$}{sl(N)} Tails}
As noted in Remark \ref{rmk:slNtailnonvanishing}, $A^t\leq0$ ensures that the tails of the colored $\mathfrak{sl}_N$ polynomials are nonvanishing. Computational data suggests that this holds for all rational knots.
\begin{conjecture}
Let $K$ be of class $\mathcal{Q}$. Then $K$ admits quiver data $Q,R,A,S$ satisfying the conditions of class $\mathcal{Q}$ such that the associated subquiver vector $A^t$, normalized as in Definition \ref{def:subquiverdata}, satisfies  $$A^t\leq0,$$ 
where the inequality holds component-wise.
\end{conjecture}
This conjecture has been verified for all rational knots up to numerator $100$ \cite{homfly-tails}. If true, it follows that the tails of the $\mathfrak{sl}_N$ polynomials are nonvanishing for all class $\mathcal{Q}$ knots and all $N\geq2.$ Additionally, this condition guarantees that the HOMFLY-PT tail does not vanish without having to assume that $\text{diag}(Q) = S - 2A$. Hence, one could have defined class $\mathcal Q$ to have this condition instead. It is unknown if one of these conditions defines a more general class than the other.

\subsection{Class \texorpdfstring{$\mathcal{Q}$}{Q}}

The tail of the colored Jones polynomial exists for all adequate links, so one may wonder if the tail of the colored HOMFLY-PT polynomial exists for all adequate links as well. While the colored HOMFLY-PT and colored Jones polynomials share many similar properties for rational knots, there are certain classes of knots for which these properties differ. For example, for alternating knots, the colored Jones polynomial always has alternating signs of the coefficients, but this is not always true of the colored HOMFLY-PT polynomial. Thus, it is conceivable that the tail of the Jones polynomial could exist for a knot in which the tail of the colored HOMFLY-PT polynomial does not exist. 

Another question to explore is what the quivers look like for knots in which the tail does not exist, like the Manx knots. These knots will not be of class $\mathcal{Q}$, and it would be interesting to see what form the quiver takes. For the Manx knots, the first coefficient of the colored Jones polynomial grows proportional to the color \citep{Lee_2022}, so that might give insight into the form of the quiver. 

Finally, we have the following inclusions:
    $$\text{Rational knots} \subseteq \text{Class $\mathcal Q$ knots} \subseteq \text{Knots with a HOMFLY-PT tail}.$$

\begin{question}
    \
    \begin{enumerate}[label=\roman*.]
        \item Are there class $\mathcal Q$ knots that are not rational knots?
        \item Are there knots with a HOMFLY-PT tail that are not in class $\mathcal Q$?
    \end{enumerate}
\end{question}

To try to resolve this question, one could attempt to calculate the quiver of a non-rational knot for which the Jones tail exists, like the knot $8_{10}$.

\appendix
\section{Geometric Computation of the Quivers for Twist Knots}
\label{TwistKAppendix}

In this appendix, we compute quivers for the twist knots $K_{(2n+1)/n}$ to support the tail derivations in Section \ref{TwistKSection}. We do so following the geometric methods of the third author in \cite{higgins2}, which were briefly outlined in Section \ref{GeoQuiverSection}. We must handle the even and odd $n$ cases separately; first, we consider the even case.

\begin{lemma}
\label{evenTknot}
    The twist knot $K_{(2n+1)/n}$ with $n=2r$ admits a quiver $Q$ given by the $5\times 5$ block matrix
    \[
    Q=\left(
    \begin{array}{ c | c | c | c | c}
    n-1  & n-1 & n-2 & n-2 & n-1 \\
    \hline
    n-1 & B_r+2 & B_r+U & B_r+L+U & B_r+2+L\\
    \hline
    n-2  & B_r+L & B_r-1 & B_r-1+L & B_r+1+L-U \\
    \hline
    n-2 & B_r+U+L & B_r-1+U & B_r & B_r+1+U \\
    \hline
    n-1 & B_r+2+U & B_r+1+U-L & B_r+1+L & B_r+3
  \end{array}\right)
    \]
where $B_r$ is the $r\times r$ matrix
\[
\begin{tikzpicture}
\node at (0,0) {
$B_r=\left( \begin{array}{ccccc}
    2r-2 & 2r-3 & 2r-3 & \dots &2r-3  \\
    2r-3 & 2r-4 & 2r-5 &\dots &  2r-5 \\
    2r-3 & 2r-5 & 2r-6 & \dots & 2r-7 \\
    \vdots & \vdots & \vdots & \ddots & \vdots \\
    2r-3 & 2r-5 & 2r-7 & \dots & 0
\end{array}\right)$};
\node at (5,0) {or};
\node at (9,0) {$(B_r)_{i,j}=\begin{cases}
    2r-2i, & i=j \\
    2r-2\,\text{min}(i,j)-1, & i\neq j
\end{cases}.$};
\end{tikzpicture}
\]
The $(1,1)$-block is the $1\times 1$ matrix $(n-1)$, so the blocks in the first row and column are constant vectors of the appropriate shape.
\end{lemma}

\begin{figure}
    \centering
    \includegraphics[height=6cm]{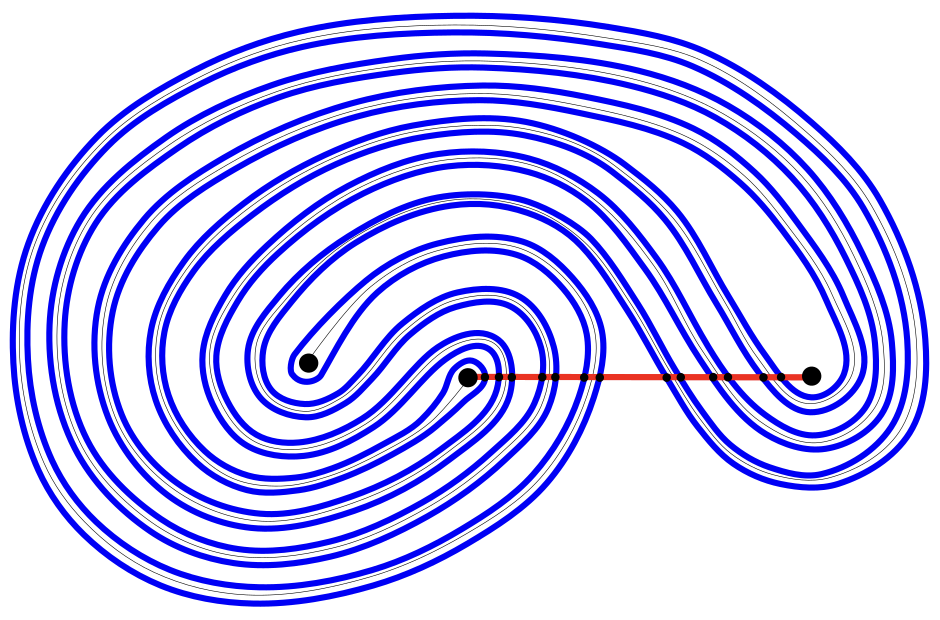}
    \caption{The two Lagrangians for $K_{13/6}$}
    \label{DK13_6}
\end{figure}

\begin{proof}
    As stated above, we compute the ``geometrically preferred'' quiver via the methods of \cite{higgins2}. In particular, the quiver will be computed as a quadratic form on the free $\mathbb{Z}$-module with basis given by Lagrangian intersections $\overline{\alpha_{(2n+1)/n}}\cap \beta=\{\kappa_1,...,\kappa_{2n+1}\}$. Furthermore, to compute the exact matrix in the lemma, this must be an ordered basis (i.e. we must fix an ordering on the intersection points).

    The Lagrangians $\overline{\alpha_{(2n+1)/n}}$ (or $\overline{\alpha_{(4r+1)/2r}}$) naturally follow a pattern where there are $2r+1$ intersections on the left and $2r$ to the right of the point where the two sides meet in the middle (see Figure \ref{DK13_6}). This can be made more explicit by applying the ``untwisting isotopy'' (or inverse Dehn twist $T^{-1}$) discussed in \cite{higgins2}. Furthermore, the punctures in $\mathcal{D}(K_{(4r+1)/2r})$ are always ordered $X^-,X^+$, and $Y$ from left to right. The leftmost intersection point will be the first point $\kappa_1$ in our ordered basis. Now, if we cut $\overline{\alpha_{(4r+1)/2r}}$ at the two ends near $X^+$ and $X^-$, we have two parallel copies of the arc $\alpha_{(4r+1)/2r}$. The second and fifth blocks come from the intersection points on one of the arcs, and the third and fourth blocks come from the other arc. Furthermore, the second and third blocks come from the remaining $2r$ intersection points on the left and the fourth and fifth blocks from the $2r$ intersections on the right.

    The precise ordering of the intersection points follows a simple pattern, so we will just use the example of $K_{13/6}$ to establish conventions. A neighborhood of $\overline{\alpha_{13/6}}\cap \beta$ looks like
    \[
    \begin{tikzpicture}
        \node at (0,0) {\includegraphics[height=1cm]{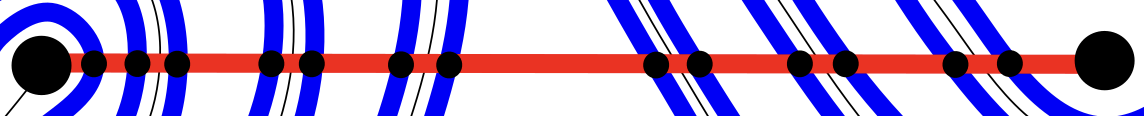}};
        \node at (-4.3,-0.7) {$\kappa_1$};
        \node at (-3.8,-0.7) {$\kappa_2$};
        \node at (-3.4,-0.7) {$\kappa_5$};
        \node at (-2.7,-0.7) {$\kappa_3$};
        \node at (-2.3,-0.7) {$\kappa_6$};
        \node at (-1.5,-0.7) {$\kappa_4$};
        \node at (-1.1,-0.7) {$\kappa_7$};
        \node at (0.8,-0.7) {$\kappa_{10}$};
        \node at (1.4,-0.7) {$\kappa_{13}$};
        \node at (2.1,-0.7) {$\kappa_9$};
        \node at (2.7,-0.7) {$\kappa_{12}$};
        \node at (3.5,-0.7) {$\kappa_8$};
        \node at (4.1,-0.7) {$\kappa_{11}$};
        \node at (-5.3,0) {$X^+$};
        \node at (5.3,0) {$Y$};
    \end{tikzpicture}
    \]
    with the intersection point labels written underneath. The same pattern used here applies to $K_{(4r+1)/2r}$.

    Observe that $\kappa_1=\kappa_\omega$, in the notation of \cite{higgins2}, so $Q_{11}=\mu_3(K_{(4r+1)/2r})=2r-1=n-1$. This determines the matrix $Q$ up to an overall shift. Now,
    \[
    Q_{1k}=n-1+w_{\Delta}(\gamma_{k,1})-2w_{X^+}(\gamma_{k,1}),
    \]
    where $w_{\Delta}(\gamma_{k,1})-2w_{X^+}(\gamma_{k,1})=-1$ if $r+1<k\le3r+1$ and $0$ otherwise, as one can easily check. Thus, it remains to prove the lemma for the lower $4\times 4$ blocks. By Lemma 4.15 of \cite{higgins2}, it suffices to prove that the diagonal blocks and the $(2,5)$-block $B_r+2+L$ have the correct form.
    
    The computation of the diagonal blocks is all essentially the same, so we will only explicitly discuss the $(2,2)$-block $B_r+2$. First, we check the diagonal of the block. We have 
    \[
    Q_{22}=Q_{11}+w_Y(\gamma_2)+w_{X^-}(\gamma_2)-3w_{X^+}(\gamma_2)=n
    \]
    since $w_Y(\gamma_2)=1$ and $w_{X^-}(\gamma_2)=w_{X^+}(\gamma_2)=0$. If $3\leq i+1\leq r+1$, we have $Q_{i+1\,i+1}=Q_{ii}-2$ because each loop $\gamma_{i+1,i}$ has the same basic shape such that $w_{X^+}(\gamma_{i+1,i})=w_Y(\gamma_{i+1,i})=1$ and $w_{X^-}(\gamma_{i+1,i})=0$. For the off-diagonal elements of the $(2,2)$-block, we can still use the regularity of the shapes of the relevant loops. In particular, for $2\leq i <j\leq r+1$, we have
    \[
    Q_{ji}=Q_{ij}=Q_{ii}+w_\Delta(\gamma_{j,i})-2w_{X^+}(\gamma_{j,i})=Q_{ii}-1
    \]
    because $w_{X^+}(\gamma_{j,i})=j-i$ and $w_\Delta(\gamma_{j,i})=2(j-i)-1$. Consequently, the diagonal block has the correct form.

    It remains to show that the $(2,5)$-block is of the form $B_r+2+L$. First, we check that the diagonal of the block is equal to the diagonal of the $(2,2)$-block. This holds because, for $1\leq i\leq r$, we have $w_\Delta(\gamma_{1+3r+i,1+i})=-2$ and $w_{X^+}(\gamma_{1+3r+i,1+i})=-1$. More generally, for $1\leq i,j\leq r$, one can check that $w_{X^+}(\gamma_{1+3r+j,1+i})=-1+j-i$ and $w_\Delta(\gamma_{1+3r+j,1+i})=-2$ if $i\geq j$ and $w_\Delta(\gamma_{1+3r+j,1+i})=-3+2(j-i)$ if $j>i$. It follows that, if we let $a_i=1+i$ and $b_j=1+3r+j$, then
    \[
Q_{a_ib_j}=Q_{a_ia_i}+w_\Delta(\gamma_{b_j,a_i)}-2w_{X^+}(\gamma_{b_j,a_i})= \begin{cases}
        n-2j+2, & j\leq i \\
        n-2i+1, & i<j,
    \end{cases}
    \]
    which is the $(i,j)$-entry of $B_r+2+L$.   
\end{proof}

\begin{remark}
    One consequence of Lemma \ref{evenTknot} is that the unique maximum in $Q$ occurs in the top-left entry of $(5,5)$-block (or the $(3r+2)$nd diagonal entry), and it has value $n+1$. For computing $Q^h$, we must then determine which entries in the $(3r+2)$nd column (or row) of $Q$ have value $n$. Lemma \ref{evenTknot} tells us that these are all the remaining entries in the $(5,5)$-block and all entries in the $(2,5)$-block in this column.
\end{remark}

\begin{figure}
    \centering
\includegraphics[height=6cm]{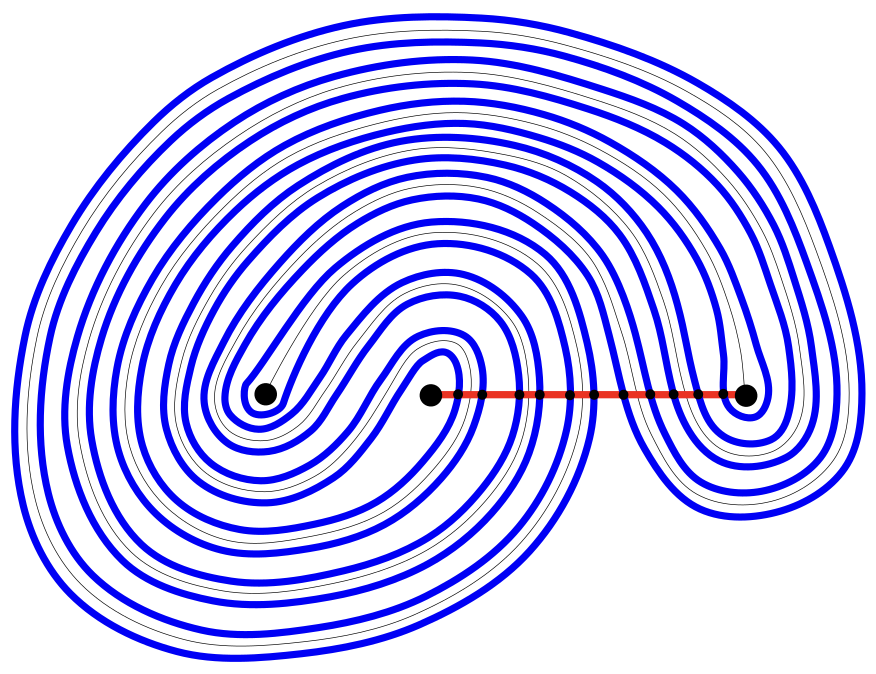}
    \caption{The two Lagrangians for $K_{11/5}$}
    \label{DK_11_5}
\end{figure}

Next, we consider the odd $n$ case.

\begin{lemma}
\label{oddTknots}
     The twist knot $K_{(2n+1)/n}$ with $n=2r-1$ admits a quiver $Q$ given by the $5\times 5$ block matrix
    \[
    Q=\left(
    \begin{array}{ c | c | c | c | c}
    B_r+2  & B_r+U & B_{r,r-1}+U-L & B_{r,r-1}+1+U & n-1 \\
    \hline
    B_r+L & B_r-1 & B_{r,r-1}-2+U & B_{r,r-1}-1+U+L & n-2\\
    \hline
    B_{r-1,r}+L-U  & B_{r-1,r}-2+L & B_{r-1} & B_{r-1}+1+U & n-2 \\
    \hline
    B_{r-1,r}+1+L & B_{r-1,r}-1+U+L & B_{r-1}+1+L & B_{r-1}+3 & n-1 \\
    \hline
    n-1 & n-2 & n-2 & n-1 & n-1
  \end{array}\right)
    \]
where $B_r$ is defined the same as in Lemma \ref{evenTknot}, $B_{r,r-1}$ is the $r\times (r-1)$ matrix obtained by deleting the $r$th column of $B_r$, and $B_{r-1,r}=B_{r,r-1}^T$.
The $(5,5)$-block is the $1\times 1$ matrix $(n-1)$, so the blocks in the fifth row and column are constant vectors of the appropriate shape.   
\end{lemma}

\begin{proof}
    As with Lemma \ref{evenTknot}, we must determine the ordering of the intersection points $\overline{\alpha_{(4r-1)/(2r-1)}}\cap \beta=\{\kappa_1,...,\kappa_{4r-1}\}$ used to obtain the matrix in the lemma statement. We take $\kappa_{4r-1}$ to be the rightmost intersection in $\mathcal{D}(K_{(2n+1)/n})$, which is the point called $\kappa_\omega$ in \cite{higgins2}. In $\mathcal{D}(K_{(2n+1)/n})$, there are now $2r$ intersection points to the left and $2(r-1)$ remaining intersection points to the right of the point in the middle where the two sides of $\overline{\alpha_{(4r-1)/(2r-1)}}$ come together, in the same sense as in Lemma \ref{evenTknot}. See Figure \ref{DK_11_5} for an example, where $r=3$. If we cut $\overline{\alpha_{(4r-1)/(2r-1)}}$ at both ends near $X^+$ and $X^-$, we have that the first and fourth blocks in $Q$ come from intersection points on one side, with the second and third blocks coming from the other. More precisely, the first and second blocks come from the $2r$ intersection points on the left and the third and fourth blocks come from the $2(r-1)$ intersection points on the right. As before, the precise ordering follows a simple pattern exemplified by $K_{11/5}$ as shown below.
    \[
    \begin{tikzpicture}
        \node at (0,0) {\includegraphics[height=0.9cm]{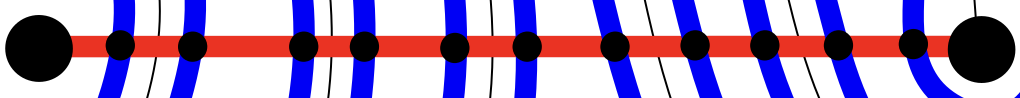}};
        \node at (-3.6,-0.65) {$\kappa_1$};
        \node at (-3,-0.65) {$\kappa_4$};
        \node at (-1.9,-0.65) {$\kappa_2$};
        \node at (-1.3,-0.65) {$\kappa_5$};
        \node at (-0.45,-0.65) {$\kappa_3$};
        \node at (0.15,-0.65) {$\kappa_6$};
        \node at (1.15,-0.65) {$\kappa_8$};
        \node at (1.9,-0.65) {$\kappa_{10}$};
        \node at (2.6,-0.65) {$\kappa_7$};
        \node at (3.2,-0.65) {$\kappa_9$};
        \node at (4,-0.65) {$\kappa_{11}$};
        \node at (-4.9,0) {$Y$};
        \node at (5,0) {$X^+$};
    \end{tikzpicture}
    \]
    Observe that the punctures in $\mathcal{D}(K_{(4r-1)/(2r-1)})$ are $X^-,Y,$ and $X^+$, from left to right. The geometric computations of the blocks in $Q$ are similar to the ones used in the proof of Lemma \ref{evenTknot}, so the remaining details are omitted.
\end{proof}

\bibliographystyle{plain}
\bibliography{refs}

\end{document}